\documentclass[a4paper,11pt]{article}
\usepackage{amsmath,amssymb,amsfonts}
\usepackage{mathrsfs}
\usepackage{latexsym}
\usepackage{CJK}
\usepackage{algorithm,algorithmic,multirow}     
\usepackage{clrscode}
\usepackage{graphicx,subfigure}
\usepackage{url}
\usepackage[T1]{fontenc}
\usepackage[utf8]{inputenc}
\usepackage{authblk}

\title{An extended finite element method with smooth nodal stress}
\author[1]{X. Peng}
\author[1]{S. Kulasegaram}
\author[1]{S. P. A. Bordas}
\author[2]{S. C. Wu}
\affil[1]{Institute of Mechanics and Advanced Materials (iMAM), Cardiff School of Engineering, Cardiff
         University, Cardiff, CF24 3AA, UK}
\affil[2]{State Key Laboratory of Traction Power, Southwest Jiaotong University, Chengdu, 610031, China}
\begin{document}

\maketitle

\begin{abstract}
The enrichment formulation of double-interpolation finite element method (DFEM) is developed in this paper. DFEM is first proposed by Zheng \emph{et al} (2011) and it requires two stages of interpolation to construct the trial function. The first stage of interpolation is the same as the standard finite element interpolation. Then the interpolation is reproduced by an additional procedure using the nodal values and nodal gradients which are derived from the first stage as interpolants. The re-constructed trial functions are now able to produce continuous nodal gradients, smooth nodal stress without post-processing and higher order basis without increasing the total degrees of freedom. Several benchmark numerical examples are performed to investigate accuracy and efficiency of DFEM and enriched DFEM. When compared with standard FEM, super-convergence rate and better accuracy are obtained by DFEM. For the numerical simulation of crack propagation, better accuracy is obtained in the evaluation of displacement norm, energy norm and the stress intensity factor.
 \\
$\boldsymbol{Keywords}$: Double-interpolation approximation; higher order element; smooth nodal stress; extended finite element method; crack propagation
\end{abstract}

\section{Introduction}

The extended finite element method (XFEM) \cite{Moes1999} is a powerful numerical tool to model cracks and discontinuity in a continuum. The essence of XFEM is to extend the finite element approximation basis by a set of enrichment functions that carry information about the character of the solution such as discontinuities and/or singularities in the context of partition of unity method (PUM) \cite{Melenk1996289}. The salient advantage of XFEM is to analyze cracks using the original underlying mesh. This is in great contrast to the traditional finite element method, which requires that the element topologies must conform to the crack geometries and the remeshing procedure needs to be performed on cracks when they evolve. The XFEM has achieved substantial developments over past decades. M\"oes and Gravouil setup the frame of 3D non-planar crack propagation \cite{NME:NME429}\cite{NME:NME430} using XFEM in conjunction with level sets method. Sukumar \emph{et al} \cite{NME:NME2344} coupled the fast marching method in order to improve the efficiency of updating the level sets functions. The XFEM for fatigue crack growth has also been incorporated into commercial software packages \cite{Eric2007}\cite{Shi20102840} with accurate error control for many cases in real-world applications \cite{NME:NME2332}\cite{Bordas20073381}.

However, the XFEM encounters a number of obstacles in becoming a robust method in spite of its successful development. Significant efforts have been devoted towards improving the accuracy and robustness of this method. In the standard XFEM, the local partition of unity is adopted, which means only certain nodes are enriched. This will result in some elements (the blending element) consisting of both regular and enriched nodes which do not conform to partition of unity. The existence of blending elements would lead to low accuracy. Chessa \emph{et al} \cite{NME:NME777} developed an enhanced strain formulation which led to a good performance of local partition of unity enrichments \cite{NME:NME2402}. Gracie \emph{et al} \cite{NME:NME2217} proposed the discontinuous Galerkin method aimed at eliminating the source of error in blending elements. More attempts can be referred from \cite{Fries2008}\cite{NME:NME2387}. In terms of integration, the additional non-polynomial enriched functions in the approximation space make the quadrature of the stiffness matrix of enriched element and blending element more delicate. Meanwhile the singularity problem in the crack tip enrichment adds the complexity for numerical integration. The traditional procedure to perform the integration is to sub-divide the enriched element and blending element into quadrature subcells to ensure the precision \cite{Moes1999}. Ventura \cite{NME:NME1570} proposed an approach to eliminate the quadrature subcells via replacing non-polynomial functions by 'equivalent' polynomials. But this method is only exact for triangular and tetrahedral elements. Another efficient integration scheme was proposed by transforming the domain integration into contour integration in \cite{NME:NME2387}. Natarajan \emph{et al} \cite{NME:NME2798} developed a new numerical integration for arbitrary polygonal domains in 2D to compute stiffness matrix. In this method, each part of the elements that are cut or intersected by a discontinuity is conformally mapped onto a unit disc using Schwarz-Christoffel mapping. Thus sub-dividing procedure is avoided. More attempts were seen in smoothed XFEM, in which interior integration is transformed into boundary integration. Then sub-dividing become unnecessary \cite{Chen2012}. Laborde \emph{et al} \cite{NME:NME1370} adopted the almost polar integration within the crack tip enriched element, which apparently improves the convergence rate. Another issue observed in standard XFEM is the non-optimal convergence rate. Many attempts have been performed in order to obtain the optimal convergence rate. One of the improvements is to use geometrical enrichment\cite{NME:NME1370}, i.e., the enrichment of a set of nodes within a radial domain of the crack tip. Nevertheless, this approach greatly deteriorates the condition number of the stiffness matrix, which limits its application to 3D or complex geometries. In order to reduce the condition number, effective preconditioners are proposed by B\'echet \cite{NME:NME1386} and Menk \emph{et al} \cite{NME:NME3032}.

Although XFEM presents great potential in industrial applications, traditional approach using finite element with remeshing technique \cite{Paluszny2011} still draw attentions of many researchers and feature in some commercial codes appearing such as FRANC2D\&3D \cite{Bittencourt1996} and Zencrack \cite{Maligno201094}. The boundary element method, which only requires the surface mesh discretization, provides a natural and convenient way to remesh crack surfaces \cite{Yan2007}\cite{Cisilino1999}.

Numerous interpolation methods have been developed in order to improve the efficiency of standard FE simulation. In terms of higher continuous FEM, Papanicolopulos and Zervos \cite{Papanicolopulos20121437}\cite{Papanicolopulos201353} created a series of triangular elements for $C^1$ continuous interpolation. Fischer \emph{et al} \cite{NME:NME2802} compared the performance of $C^1$ finite elements and $C^1$ natural element method (NEM) applied in non-linear gradient elasticity. Various meshfree methods are introduced and used widely in engineerings problems. The element free Galerkin method (EFG) \cite{Belytschko1994} adopts the moving least-squares interpolants to construct trial and test functions which can easily obtain higher continuous in both variable and its gradient. One similar method is the reproducing kernel particle method (RKPM)\cite{FLD:FLD1650200824}. The meshfree radial basis functions method (RBFs) utilizes the radial basis function to interpolate scattered nodal data and is employed with point interpolation method (PIM) by Liu \emph{et al}. The radial PIM (RPIM)\cite{Liew2004}\cite{Wang20022611}\cite{Liu2005}, consists of both radial basis and polynomial basis in the approximation, which can avoid the singularity of the moment matrix in polynomial basis. The maximum-entropy method (MAXENT) proposed by Arroyo and Oritiz bridges Delaunay triangulation and maximum-entropy statistical interface \cite{NME:NME1534}. And its second order form can be applied to optimize the support size in meshfree methods \cite{NME:NME2793}. Liu \emph{et al} developed a smoothed FEM (SFEM) based on the generalization of the strain smoothing technique to functions in the G space \cite{Liu2008}. According to different smoothing integration domain, the SFEM can be classified as node-based smoothed FEM (NSFEM) and edge-based smoothed FEM (ESFEM). Many researchers subsequently investigated the new methods into modelling discontinuities using the enriched form by applying partition of unity into these new methods, such as extended SFEM \cite{Chen2012}\cite{Bordas2010}\cite{NME:NME2868}\cite{Vu-Bac2011}\cite{Jiang2013}\cite{Bordas2011}or extended EFG \cite{Bordas2008943}.

In this paper, a new approximation, which is based on the simplex mesh discretization and has $C^1$ continuity on most nodes, is investigated. This approximation procedure can be named the extended double-interpolation FEM (XDFEM) due to the fact that two consecutive stages of interpolation are adopted in the construction of this approximation. The first stage of interpolation is performed by standard FEM to obtain nodal variables and nodal gradients. The problem field is reproduced in the latter interpolation using the nodal values and gradients derived from the previous interpolation. The re-constructed trial functions will maintain the $C^1$ continuity at nodes \cite{Zheng2009}. Cubic polynomial basis are contained in trial functions without increasing total degrees of freedom (DOFs). This feature enhances the adaptivity of DFEM simplex mesh to simulate high gradient field near the crack tip \cite{doi:10.1142/S0219876212500557} and improves the accuracy without increasing computational cost. Analogous to meshfree methods, the nodal stress can be calculated straightforwardly by interpolation in DFEM without any post-processing.

 The paper is organised as follows. In section $2$, the DFEM formulation for 1D and 2D are systematically demonstrated with a example of 1D bar simulation. Section 3 presents the discretized formulation of XDFEM for linear elastic fracture mechanics. Several numerical examples are presented to illustrate the advantages and probable limitations of DFEM and XDFEM in section $4$. Finally, in section 5, concluding remarks are made with pointers to possible future work.

\section{The double-interpolation approximation}

\subsection{1D approximation by double-interpolation}

\begin{figure}
\begin{center}
\includegraphics[width=0.5\textwidth]{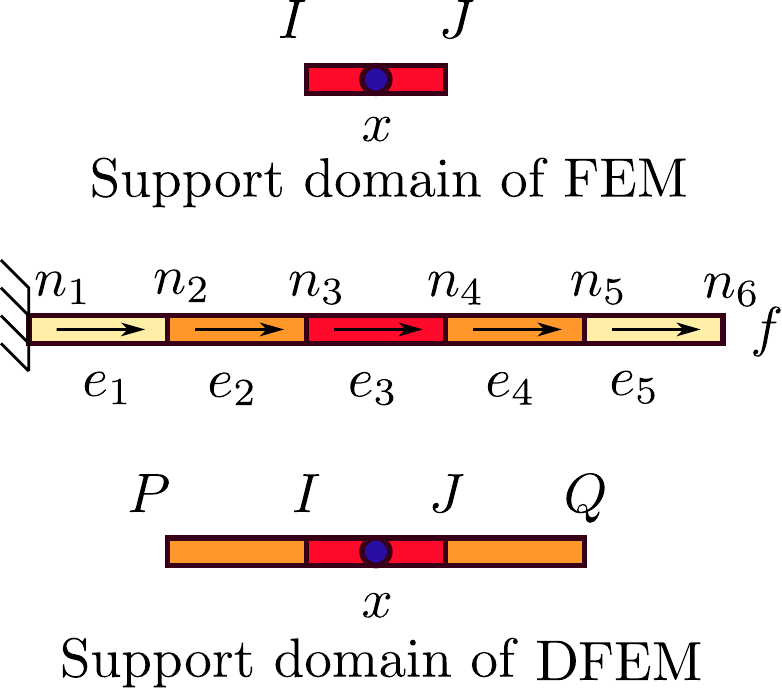}
\end{center}
\caption{Discretization of the 1D domain and the element domain of FEM and DFEM} \label{domain1D}
\end{figure}

The basic idea of the double-interpolation approximation is to interpolate the unknown fields, using both the primary nodal values and nodal gradients, which are generated by the finite element interpolation in simplex mesh discretization. The proposed 1D double-interpolation is comparable to the Hermite interpolation. Fig.\ref{domain1D} shows a 1D domain which is discretized by six 1D elements. For the point of interest $x$ in element $e_3$, the numerical value of the displacement can be interpolated by

\begin{equation}  \label{eq:00}
u^h(x)=\phi_I(x)u^I+\psi_I(x)u_{,x}^I+\phi_J(x)u^J+\psi_J(x)u_{,x}^J,\forall x\in[0,l]
\end{equation}
where $u^I, u_{,x}^I$ denote the nodal displacement and nodal derivative of the displacement field at node $I$, respectively. $l$ is the length of 1D domain. $\phi_I, \psi_I, \phi_J, \psi_J$ are the cubic Hermite basis polynomials given by;

\begin{subequations}
\begin{equation}
\phi_I(x)=\left(1+2\left(\frac{x-x_I}{x_J-x_I}\right)\right)\left(\frac{x-x_J}{x_J-x_I}\right)^2
\end{equation}
\begin{equation}
\psi_I(x)=(x-x_I)\left(\frac{x-x_J}{x_J-x_I}\right)^2
\end{equation}
\begin{equation}
\phi_J(x)=\left(1-2\left(\frac{x-x_J}{x_J-x_I}\right)\right)\left(\frac{x-x_I}{x_J-x_I}\right)^2
\end{equation}
\begin{equation}
\psi_J(x)=(x-x_J)\left(\frac{x-x_I}{x_J-x_I}\right)^2
\end{equation}
\end{subequations}

If we define the local coordinates as follows,

\begin{equation}
L_I(x)=\frac{x-x_J}{x_I-x_J},\quad L_J(x)=\frac{x-x_I}{x_J-x_I}
\end{equation}
and $c=x_J-x_I$, then the Hermite basis polynomials can be written as:

\begin{subequations}
\begin{equation}
\phi_I(x)=L_I(x)+\left(L_I(x)\right)^2L_J(x)-L_I(x)\left(L_J(x)\right)^2
\end{equation}
\begin{equation}
\psi_I(x)=cL_J(x)\left(L_I(x)\right)^2
\end{equation}
\begin{equation}
\phi_J(x)=L_J(x)+\left(L_J(x)\right)^2L_I(x)-L_J(x)\left(L_I(x)\right)^2
\end{equation}
\begin{equation}
\psi_J(x)=-cL_I(x)\left(L_J(x)\right)^2
\end{equation}
\end{subequations}

we use the "average" nodal gradients ($\bar{u}_{,x}^I$,$\bar{u}_{,x}^J$) derived from finite element interpolation at each node to replace the gradients ($u_{,x}^I,u_{,x}^J$) in equation (\ref{eq:00}). The calculation of the average nodal gradients at node $I$ will be presented as follows. We name $\Lambda_I$ as the support domain set of node $I$. For instance, in the support domain of DFEM in Fig.\ref{domain1D}, the adjacent elements of node $I$ are $e_2, e_3$, which are the only members of $\Lambda_I$ here. All the support nodes for point of interest $x$ are denoted by $\mathscr{N}_S$. Adopting the above notations, equation (\ref{eq:00}) can be rewritten as:

\begin{equation} \label{eq:01}
u^h(x)=\phi_I(x)u^I+\psi_I(x)\bar{u}_{,x}^I+\phi_J(x)u^J+\psi_J(x)\bar{u}_{,x}^J
\end{equation}
where
\begin{equation}\label{eq:02}
u^I=u(x_I)=N_I(x_I)u^I+N_J(x_I)u^J
\end{equation}
\begin{equation}\label{eq:03}
\bar{u}_{,x}^I=\bar{u}_{,x}(x_I)={\omega_{e_{2,I}}u_{,x}^{e_{2}}(x_I)}+{\omega_{e_{3,I}}u_{,x}^{e_{3}}(x_I)}
\end{equation}
in which $N_I$, $N_J$ are FEM shape functions. It is important to point out here that  the notation of FEM shape function "$N$" without superscript implies that the shape function is associated with the finite element where the point of interest is located.  $u_{,x}^{e_{2}}(x_I)$ is the nodal derivative in element $e_{2}$, which belongs to the $\Lambda_I$. $\omega_{e_{2,I}}$ denotes the weight of element $e_2$ in $\Lambda_I$. These parameters are calculated by:

\begin{equation} \label{eq:04}
u_{,x}^{e_{2}}(x_I)=N_{P,x}^{e_{2}}(x_I)u^P+N_{I,x}^{e_{2}}(x_I)u^I
\end{equation}
\begin{equation} \label{eq:05}
\omega_{e_{2,I}}=\frac{{meas}(e_{2,I})}
{{meas}(e_{2,I})+{meas}(e_{3,I})}
\end{equation}
where $N_P^{e_{2}}(x), N_{P,x}^{e_{2}}(x),N_I^{e_{2}}(x), N_{I,x}^{e_{2}}(x)$ are the linear finite element basis functions and their derivatives, in element $e_{2}$. ${meas(\cdot)}$ denotes the measure of characteristic dimension of the specified shape. In the 1D case, it means the length of the element.

Substituting equations (\ref{eq:05}) and (\ref{eq:04}), into equations (\ref{eq:03}) yields:

\begin{equation} \label{eq:06}
\begin{aligned}
\bar{u}_{,x}^I=\bar{u}_{,x}(x_I)=&\omega_{e_{2,I}}\left(N_{P,x}^{e_{2}}(x_I)u^P+N_{I,x}^{e_{2}}(x_I)u^I\right)+\\
   &             \omega_{e_{3,I}}\left(N_{I,x}^{e_{3}}(x_I)u^I+N_{J,x}^{e_{3}}(x_I)u^{J}\right)\\
\end{aligned}
\end{equation}
which can be rewritten as:
\begin{equation}
\begin{aligned}
\bar{u}_{,x}^I=
\left[\begin{array}{cccc}
\omega_{e_{2,I}}N_{P,x}^{e_{2}} &
\omega_{e_{2,I}}N_{I,x}^{e_{2}}+\omega_{e_{3,I}}N_{I,x}^{e_{3}} &
\omega_{e_{3,I}}N_{J,x}^{e_{3}}
\end{array}\right]
\left[\begin{array}{cccc}
u^P\\
u^I\\
u^J
\end{array}\right]
\end{aligned}
\end{equation}

By defining,
\begin{equation}
\bar{N}_{L,x}(x_I)=\sum_{e_i\in\Lambda_I}\omega_{e_{i,I}}N_{L,x}^{e_i}(x_I),\quad L\in \mathscr{N}_S
\end{equation}
the averaged derivative at node $I$ can be written as
\begin{equation} \label{eq:07}
\bar{u}_{,x}^I=\bar{u}_{,x}(x_I)=\bar{N}_{P,x}(x_I)u^P+\bar{N}_{I,x}(x_I)u^I+\bar{N}_{J,x}(x_I)u^J
\end{equation}

Now, substituting equations (\ref{eq:02}) and (\ref{eq:07}) into equations (\ref{eq:01}) results in:

\begin{equation}
\begin{aligned} \label{eq:08}
u^h(x)=&\phi_I(x)\left(N_I(x_I)u^I+N_J(x_I)u^J\right)+\\
&
       \psi_I(x)\left(\bar{N}_{P,x}(x_I)u^P+\bar{N}_{I,x}(x_I)u^I+\bar{N}_{J,x}(x_I)u^J\right)+ \\
&
       \phi_J(x)\left(N_I(x_J)u^I+N_J(x_J)u^J\right)+\\
&
       \psi_J(x)\left(\bar{N}_{I,x}(x_J)u^I+\bar{N}_{J,x}(x_J)u^J+\bar{N}_{Q,x}(x_J)u^Q\right) \\
=&\psi_I(x)\bar{N}_{P,x}(x_I)u^P+\\
&
\left(\phi_I(x)N_I(x_I)+\psi_I(x)\bar{N}_{I,x}(x_I)+
\phi_J(x)N_I(x_J)+\psi_J(x)\bar{N}_{I,x}(x_J)\right)u^I+\\
&
\left(\phi_I(x)N_J(x_I)+\psi_I(x)\bar{N}_{J,x}(x_I)+
\phi_J(x)N_J(x_J)+\psi_J(x)\bar{N}_{J,x}(x_J)\right)u^J+\\
&
\psi_J(x)\bar{N}_{Q,x}(x_J)u^Q
\end{aligned}
\end{equation}

Hence, by defining,

\begin{equation}
\hat{N}_L(x)=\phi_I(x)N_L(x_I)+\psi_I(x)\bar{N}_{L,x}(x_I)+
             \phi_J(x)N_L(x_J)+\psi_J(x)\bar{N}_{L,x}(x_J)
\end{equation}
the final form for the double-interpolation approximation can be obtained as:

\begin{equation}
u^h(x)=\sum_{L\in \mathscr{N}_S}\hat{N}_L(x) u^L
\end{equation}

\begin{figure}[htbp]
  \centering
  \subfigure[]{\includegraphics[width=1.\textwidth]{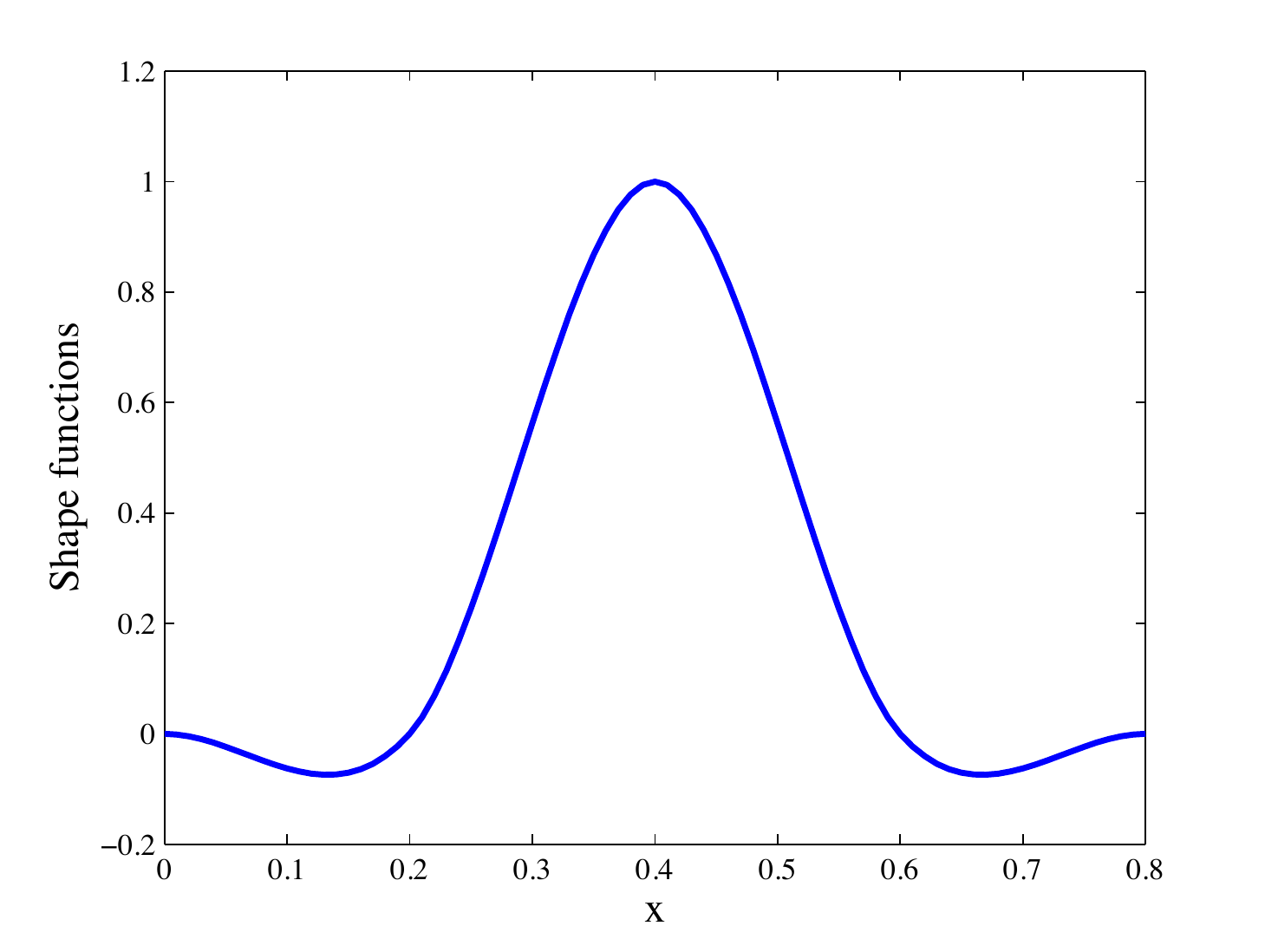}}
  \subfigure[]{\includegraphics[width=1.\textwidth]{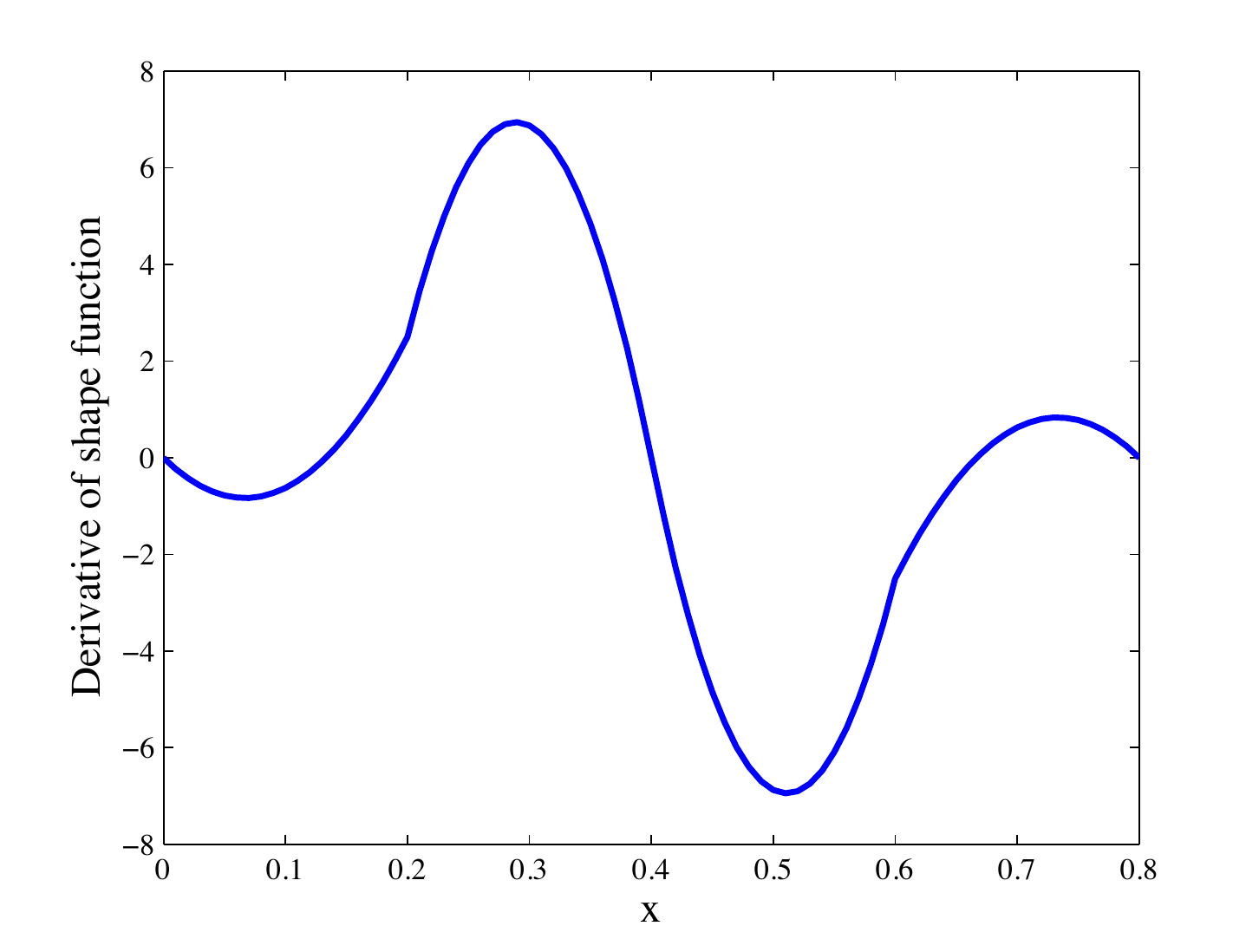}}
  \caption{The 1D DFEM shape function and its derivative at node 3}
  \label{1dshap}
\end{figure}

Fig.\ref{1dshap} shows the shape function and derivative for 1D.

To demonstrate the accuracy and characteristics of DFEM a numerical example is considered in the following discussion. For this purpose a 1D bar (as illustrated in Fig.\ref{domain1D}) problem is solved using both DFEM and FEM. The equilibrium equation for the 1D problem is defined as:
\begin{subequations}
\begin{equation}
EA\frac{\text{d}^2u}{\text{d}x^2}+f=0
\end{equation}
\begin{equation}
u|_{x=0}=0
\end{equation}
\end{subequations}
where $f$ is a uniform body force applied to the 1D bar. The exact solution for displacement and stress are given by:
\begin{subequations}
\begin{equation}
u(x)=\frac{fL^2}{EA}\left(\frac{x}{L}-\frac{1}{2}\left(\frac{x}{L}\right)^2\right)
\end{equation}
\begin{equation}
\sigma(x)=\frac{fL}{A}\left(1-\frac{x}{L}\right)
\end{equation}
\end{subequations}
where $L$, $A$ and $E$ are the total length, area of the cross section and Young's Modulus respectively. For simplicity, all these parameters are assumed to have unit value in the simulation.

\begin{figure}[htbp]
  \centering
  \subfigure[]{\includegraphics[width=1.\textwidth]{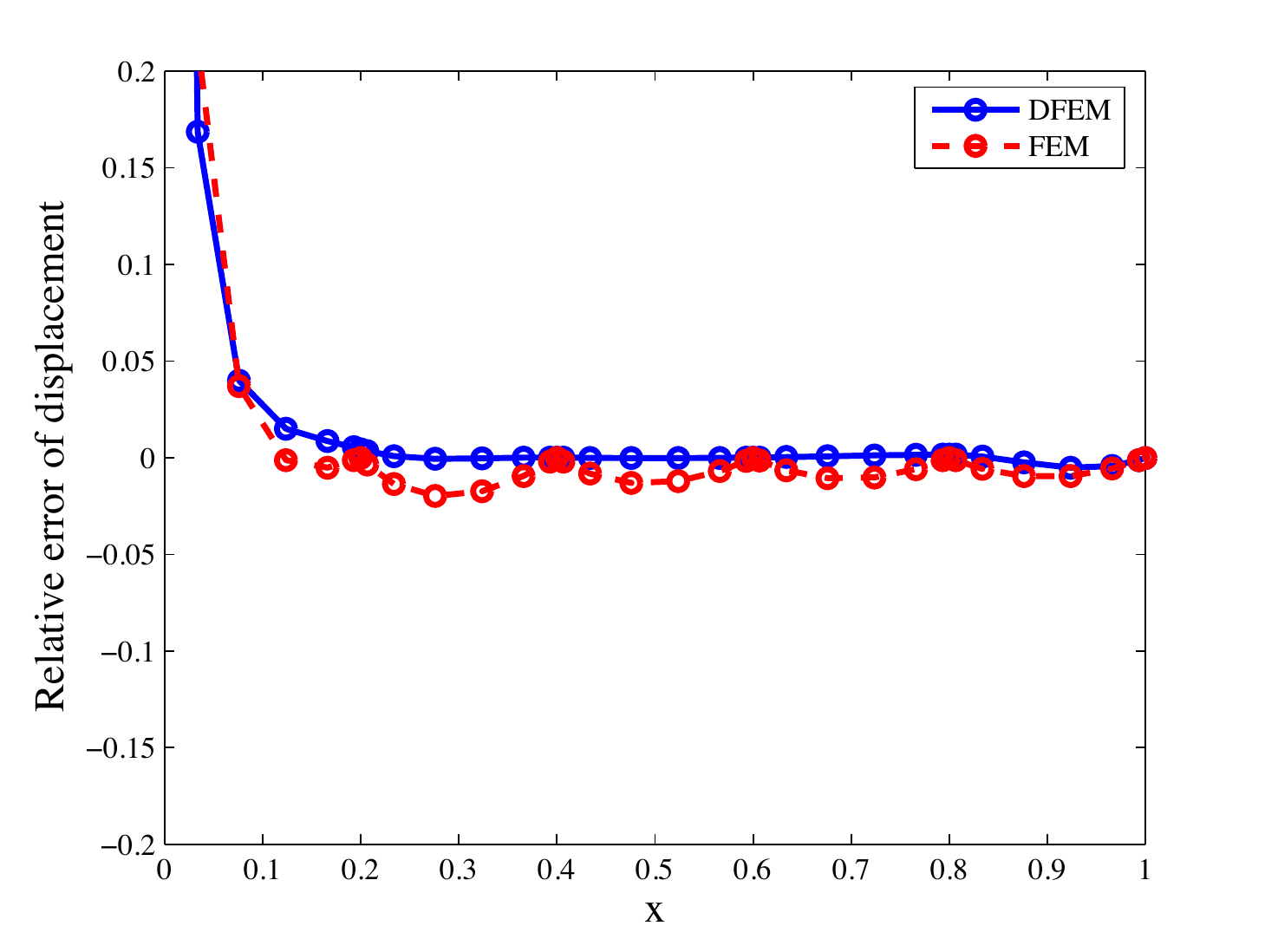}}
  \subfigure[]{\includegraphics[width=1.\textwidth]{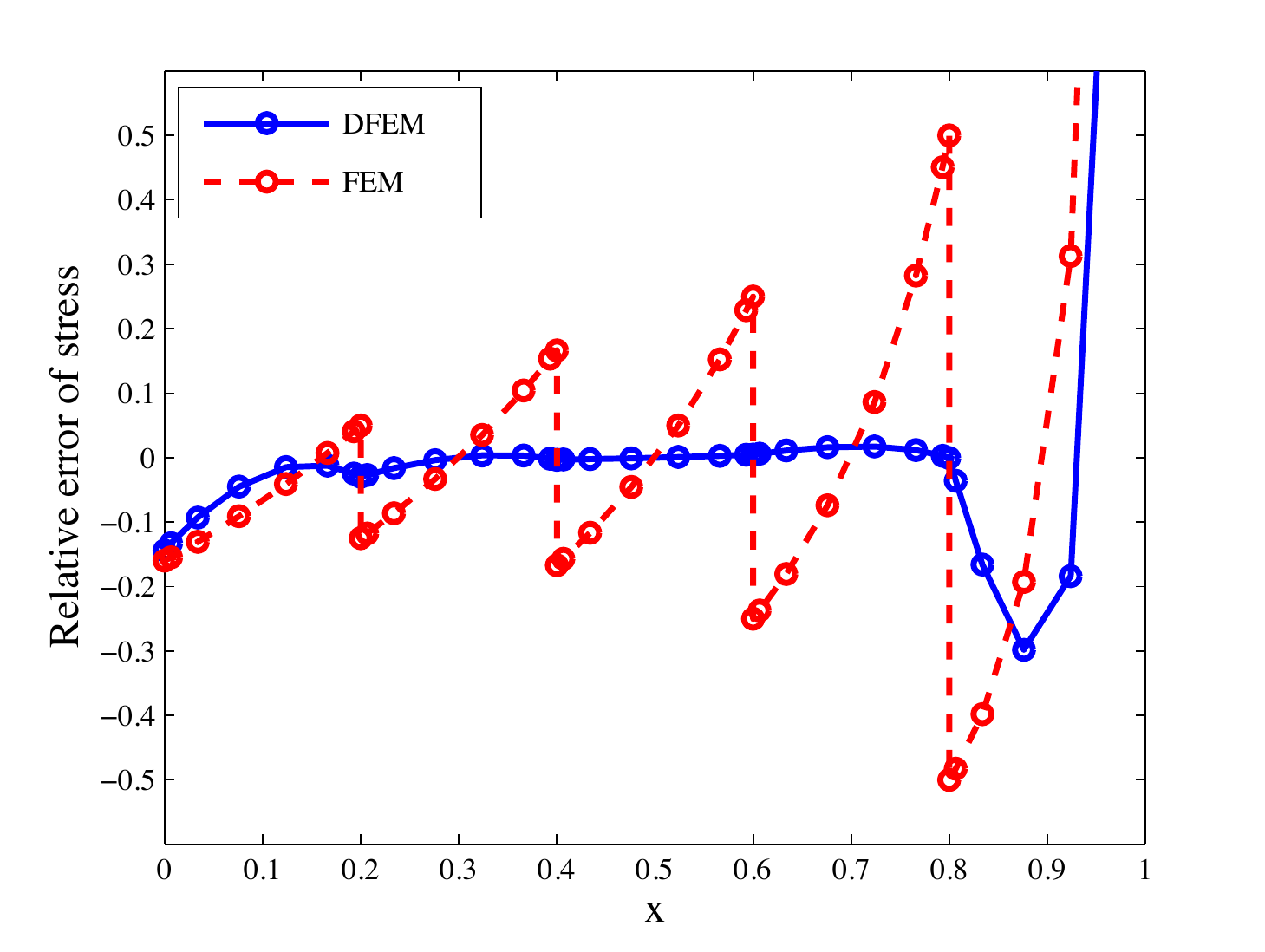}}
  \caption{The comparison of FEM and DFEM results for 1D bar:(a)displacement;(b)stress}  \label{1dresult}
\end{figure}

\begin{figure}[htbp]
  \centering
  \subfigure{\includegraphics[width=1.\textwidth]{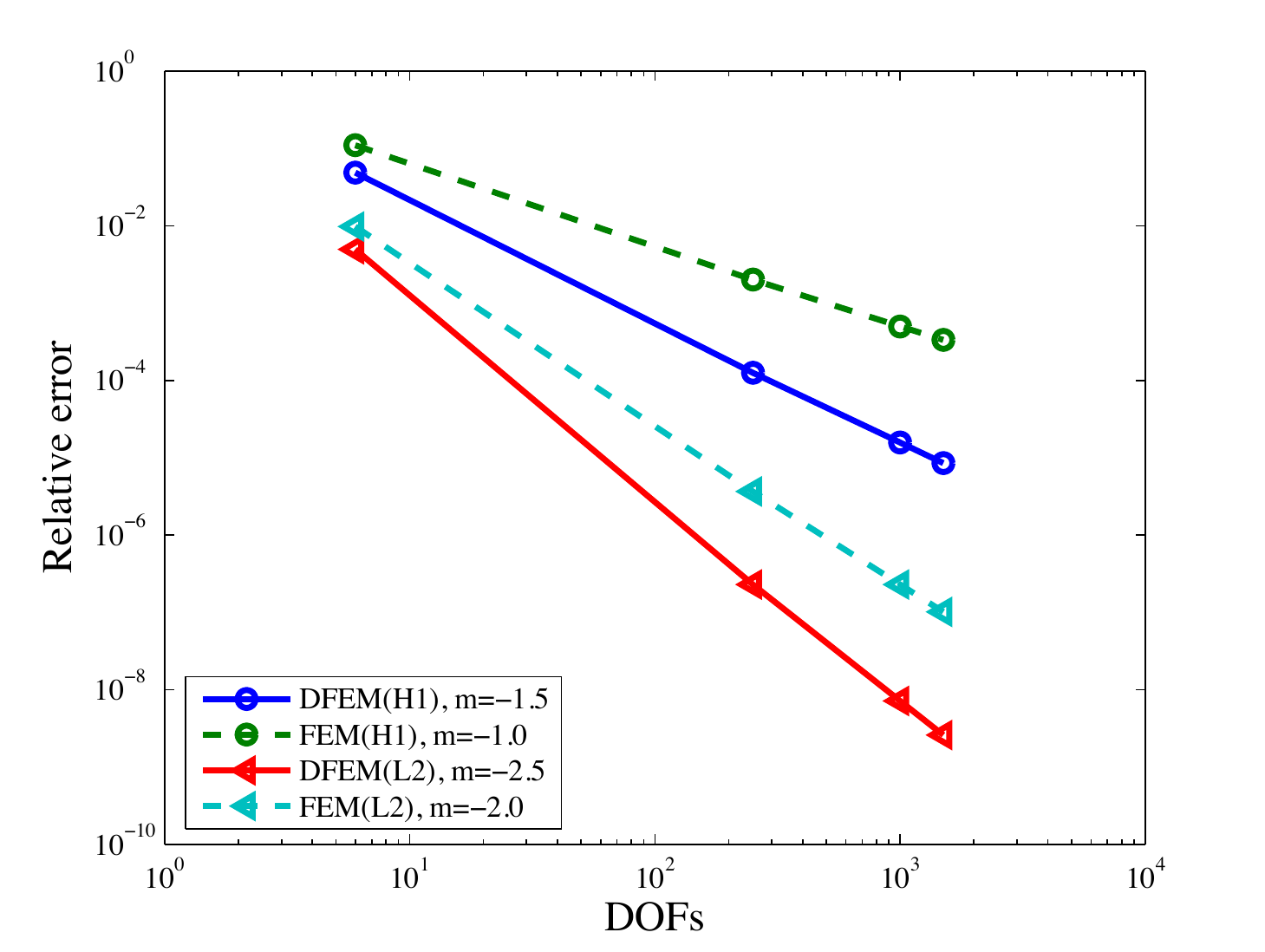}}
  \caption{Relative error in displacement and energy norm of 1D bar}  \label{1dnorm}
\end{figure}

Fig.\ref{1dresult} compares the displacement and stress values obtained by both FEM and DFEM. It can be observed from the figure that the results of DFEM captures the exact quadratic solution much better than the piecewise linear curve obtained by the FEM for the displacement field. In addition, apart from the end or boundary points, the numerical results of DFEM for stress distribution agrees well with the exact solution. The deterioration of the DFEM solution near the boundary nodes are attributed to the automatic recover of the calculation of nodal gradients at the end points, which will be explained in the following section.
Fig.\ref{1dnorm} plots the relative error in displacement and energy norm of the 1D bar problem (The definitions of these norms can be referred from section 4). It is clearly illustrated that DFEM approximation in 1D can achieve the convergence rate between the linear and quadratic order.
\subsection{2D approximation by double interpolation}

\begin{figure}
\begin{center}
\includegraphics[width=0.9\textwidth]{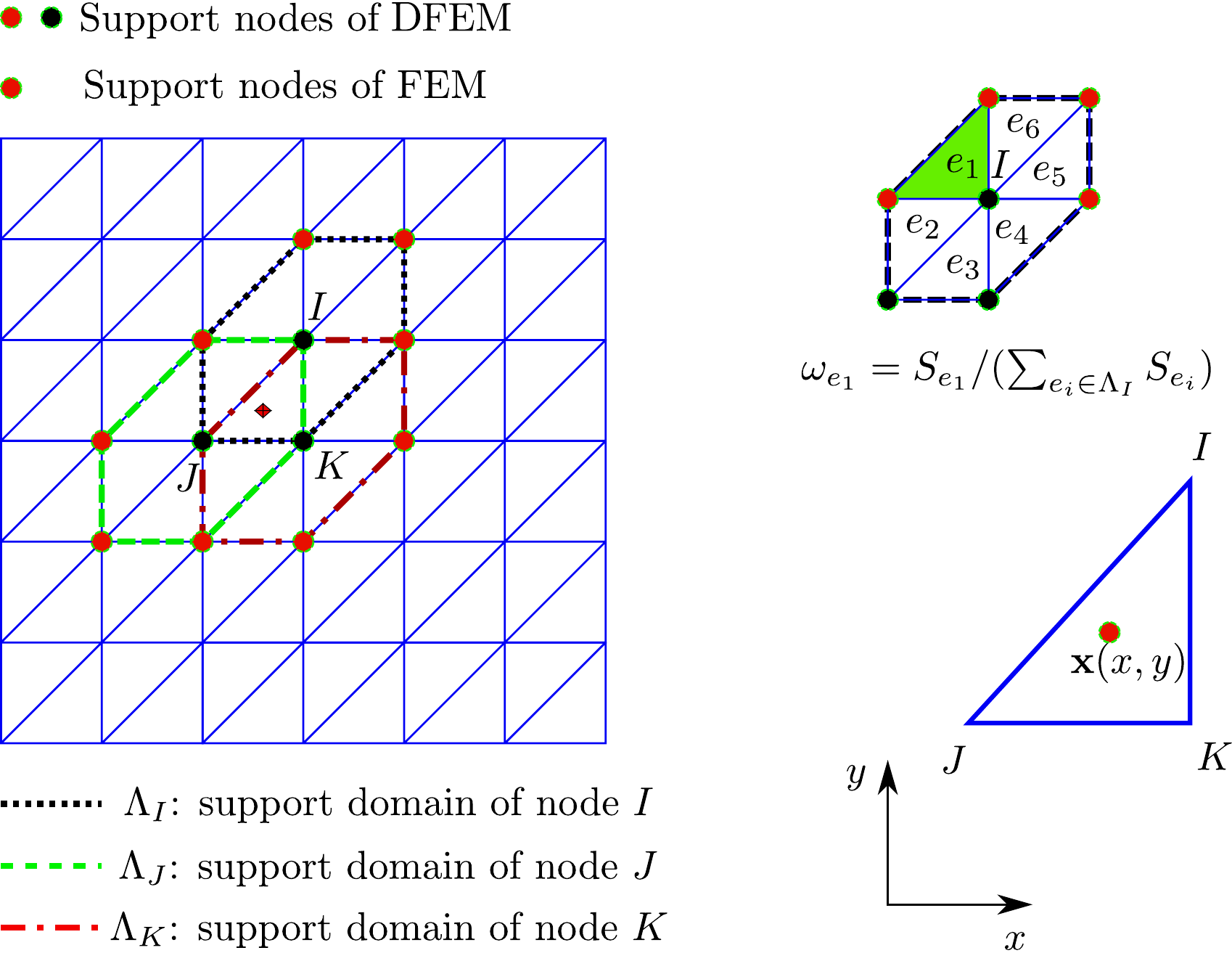}
\end{center}
\caption{Illustration for the support domain of DFEM} \label{TFEM}
\end{figure}

As illustrated in Fig.\ref{TFEM}, $\mathbf{x}$ denotes the coordinate of the point of interest in triangle $IJK$. Similar to 1D formulation, the 2D double-interpolation approximation in a triangular elemental mesh can be expressed as follows:

\begin{equation}
\mathbf{u}^h(\mathbf{x})=\sum_{L=1}^{N_S}{\hat{N}_L(\mathbf{x})\mathbf{u}^L}
\end{equation}

\begin{equation}\label{eq:09}
\begin{aligned}
\hat{N}_L(\mathbf{x})=&\phi_I(\mathbf{x}) N_L(\mathbf{x}_I)+\psi_I(\mathbf{x})\bar{N}_{L,x}(\mathbf{x}_I)+\varphi_I(\mathbf{x})\bar{N}_{L,y}(\mathbf{x}_I)+\\
&
\phi_J(\mathbf{x})N_L(\mathbf{x}_J)+\psi_J(\mathbf{x})\bar{N}_{L,x}(\mathbf{x}_J)+\varphi_J(\mathbf{x})\bar{N}_{L,y}(\mathbf{x}_J)+\\
&
\phi_K(\mathbf{x}) N_L(\mathbf{x}_K)+\psi_K(\mathbf{x})\bar{N}_{L,x}(\mathbf{x}_K)+\varphi_K(\mathbf{x})\bar{N}_{L,y}(\mathbf{x}_K)
\end{aligned}
\end{equation}
where $L=1,2,...,\mathscr{N}_S$. $\mathbf{u}^L$ is the nodal displacement vector. In the following discussion the evaluation of the average derivative of the shape function at node $I$ is considered. The average derivative of the shape function at node $I$ can be written as:

\begin{subequations} \label{eq:14}
\begin{equation}
\bar{N}_{L,x}(\mathbf{x}_I)=\sum_{e_{i,I}\in{\Lambda_I}}{\omega_{e_{i,I}} N_{L,x}^{e_i}(\mathbf{x}_I)}
\end{equation}
\begin{equation}
\bar{N}_{L,y}(\mathbf{x}_I)=\sum_{e_{i,I}\in{\Lambda_I}}{\omega_{e_{i,I}} N_{L,y}^{e_i}(\mathbf{x}_I)}
\end{equation}
\end{subequations}
$\omega_{e_{i,I}}$ is the weight of element $e_{i}$ in $\Lambda_I$ and is computed by:
\begin{equation}  \label{eq:10}
\omega_{e_{i,I}}={meas}(e_{i})/\sum_{e_i \in \Lambda_I}{{meas}(e_i)}
\end{equation}
An example of how to evaluate the weight function of an element is presented in Fig.\ref{TFEM}.
$\phi_I, \psi_I$ and $\varphi_I$ form the polynomial basis (at node $I$), which satisfies the following interpolating conditions:
\begin{equation}
\begin{array}{cccc}
\phi_I(\mathbf{x}_L)=\delta_{IL}, & \phi_{I,x} (\mathbf{x}_L)=0\quad,           & \phi_{I,y} (\mathbf{x}_L)=0\quad \\
\psi_I(\mathbf{x}_L)=0\quad,           & \psi_{I,x} (\mathbf{x}_L)=\delta_{IL}, & \psi_{I,y} (\mathbf{x}_L)=0\quad \\
\varphi_I(\mathbf{x}_L)=0\quad,        & \varphi_{I,x} (\mathbf{x}_L)=0\quad,           & \varphi_{I,y} (\mathbf{x}_L)=\delta_{IL}
\end{array}
\end{equation}

And these polynomial basis functions are given by:
\begin{subequations}
\begin{equation}
\begin{aligned}
\phi_I(\mathbf{x})=&L_I(\mathbf{x})+\left(L_I(\mathbf{x})\right)^2L_J(\mathbf{x})+\left(L_I(\mathbf{x})\right)^2 L_K(\mathbf{x}) \\
& -L_I(\mathbf{x})\left(L_J(\mathbf{x})\right)^2 -L_I(\mathbf{x})\left(L_K(\mathbf{x})\right)^2
\end{aligned}
\end{equation}
\begin{equation}
\begin{aligned}
\psi_I(\mathbf{x})=&-c_J\left(L_K(\mathbf{x})\left(L_I(\mathbf{x})\right)^2+\frac{1}{2}L_I(\mathbf{x})L_J(\mathbf{x})
L_K(\mathbf{x})\right)+  \\
&  c_K\left(\left(L_I(\mathbf{x})\right)^2L_J(\mathbf{x})+\frac{1}{2}L_I(\mathbf{x})L_J(\mathbf{x})L_K(\mathbf{x})\right)
\end{aligned}
\end{equation}
\begin{equation}
\begin{aligned}
\varphi_I(\mathbf{x})=&b_J\left(L_K(\mathbf{x})\left(L_I(\mathbf{x})\right)^2+
\frac{1}{2}L_I(\mathbf{x})L_J(\mathbf{x})L_K(\mathbf{x})\right)-\\
&  b_K\left(\left(L_I(\mathbf{x})\right)^2L_J(\mathbf{x})
+\frac{1}{2}L_I(\mathbf{x})L_J(\mathbf{x})L_K(\mathbf{x})\right)
\end{aligned}
\end{equation}
\end{subequations}

Note that the polynomial basis functions $\phi_J, \psi_{J}, \varphi_J,\phi_K, \psi_K$ and $\varphi_K$ can be obtained by the above definitions via cyclic permutation of indices $I,J$ and $K$. In the above equations, $L_I(\mathbf{x}), L_J(\mathbf{x})$ and $L_K(\mathbf{x})$ are the area coordinates of point of interest $\mathbf{x}$ in triangle $IJK$. For point of interest $\mathbf{x}$ in Fig.\ref{TFEM}, the $L_I(\mathbf{x}), b_I$ and $c_I$ are presented as follows.

\begin{subequations}
\begin{equation}
L_I(\mathbf{x})=\frac{1}{2\bigtriangleup}(a_I+b_Ix+c_Iy)
\end{equation}
\begin{equation}
a_I=x_Jy_K-x_Ky_J
\end{equation}
\begin{equation}
b_I=y_J-y_K
\end{equation}
\begin{equation}
c_I=x_K-x_J
\end{equation}
\end{subequations}
where $\bigtriangleup$ is the area of triangle $IJK$. Further, $L_J$, $L_K$, $b_J$, $b_K$, $c_J$ and $c_K$ can be obtained using the above definitions via cyclic permutation of indices $I,J$ and $K$.

When the point of interest lies on one of the edges, for example on edge $IJ$, the basis functions will boil down to 1D basis functions and will be consistent with the 1D form presented in the preceding section.

The DFEM shape functions inherit the properties such as the linear completeness, partition of unity, Dirac delta properties as that of FEM shape functions \cite{Zheng2009}. In addition, the 2D DFEM possesses $C^1$ continuity at nodes and $C^0$ continuity on edges. Compared to 3-noded triangular element, the DFEM basis functions can reach higher order without the introduction of additional nodes. However, this attractive feature comes with the price of an increased bandwith as the neighboring nodes are used to obtain the nodal gradient in the second interpolation. The details of such computational costs will be discussed in the following sections.

\begin{figure}[htbp]
  \centering
  \subfigure[Interior of the 2D domain]{\includegraphics[width=0.35\textwidth]{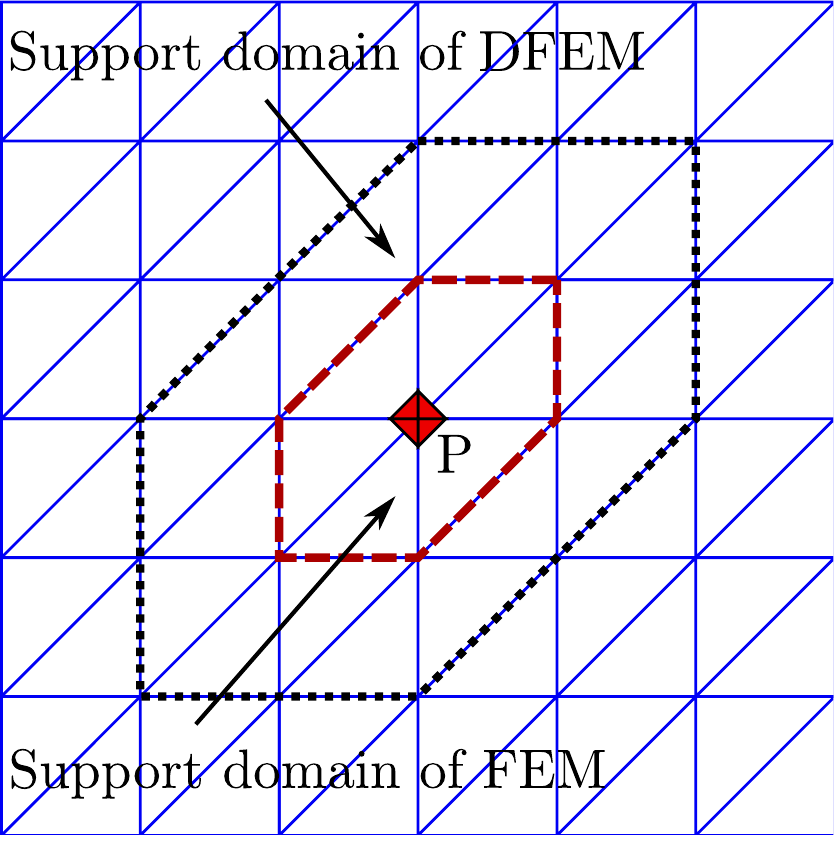}}
  \subfigure[3D plot]{\includegraphics[width=0.55\textwidth]{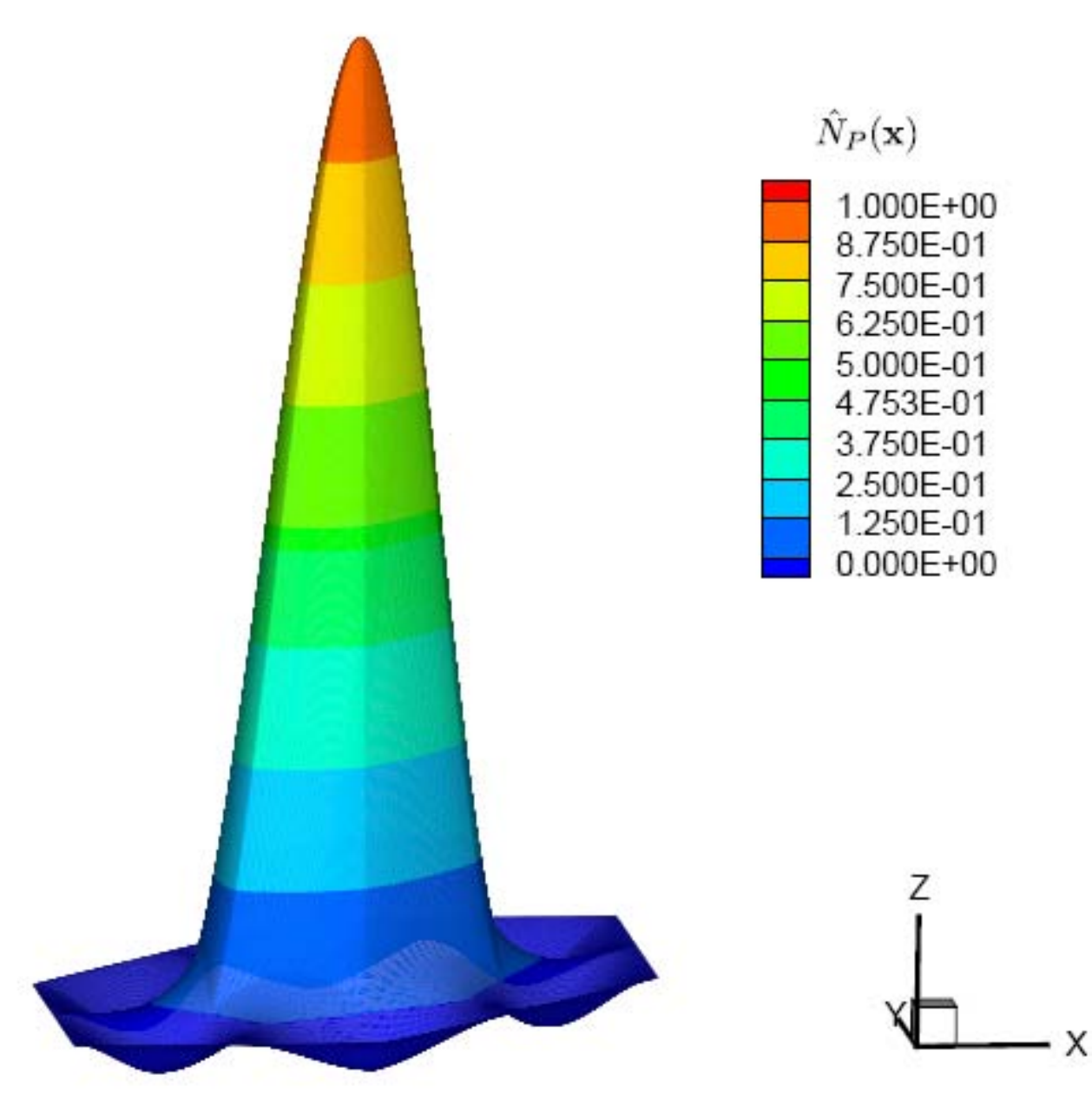}}
  \subfigure[Boundary of the 2D domain]{\includegraphics[width=0.35\textwidth]{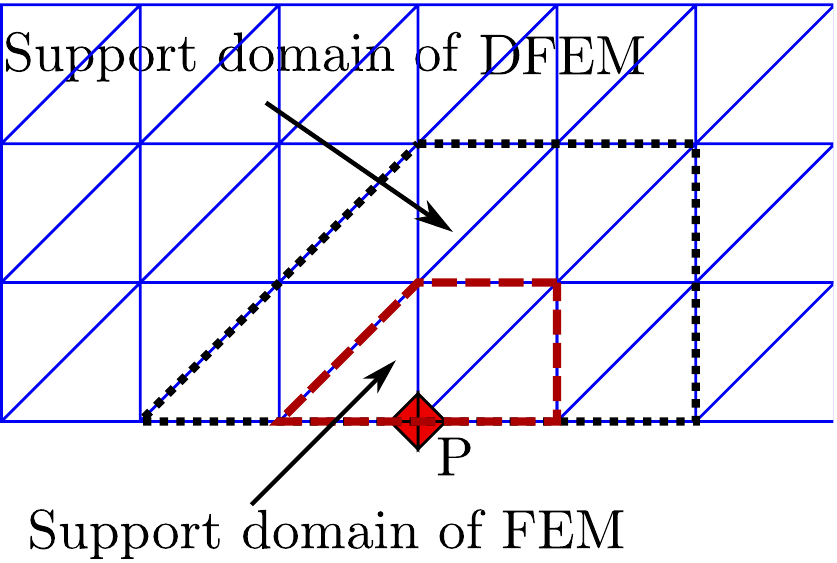}}
  \subfigure[3D plot]{\includegraphics[width=0.55\textwidth]{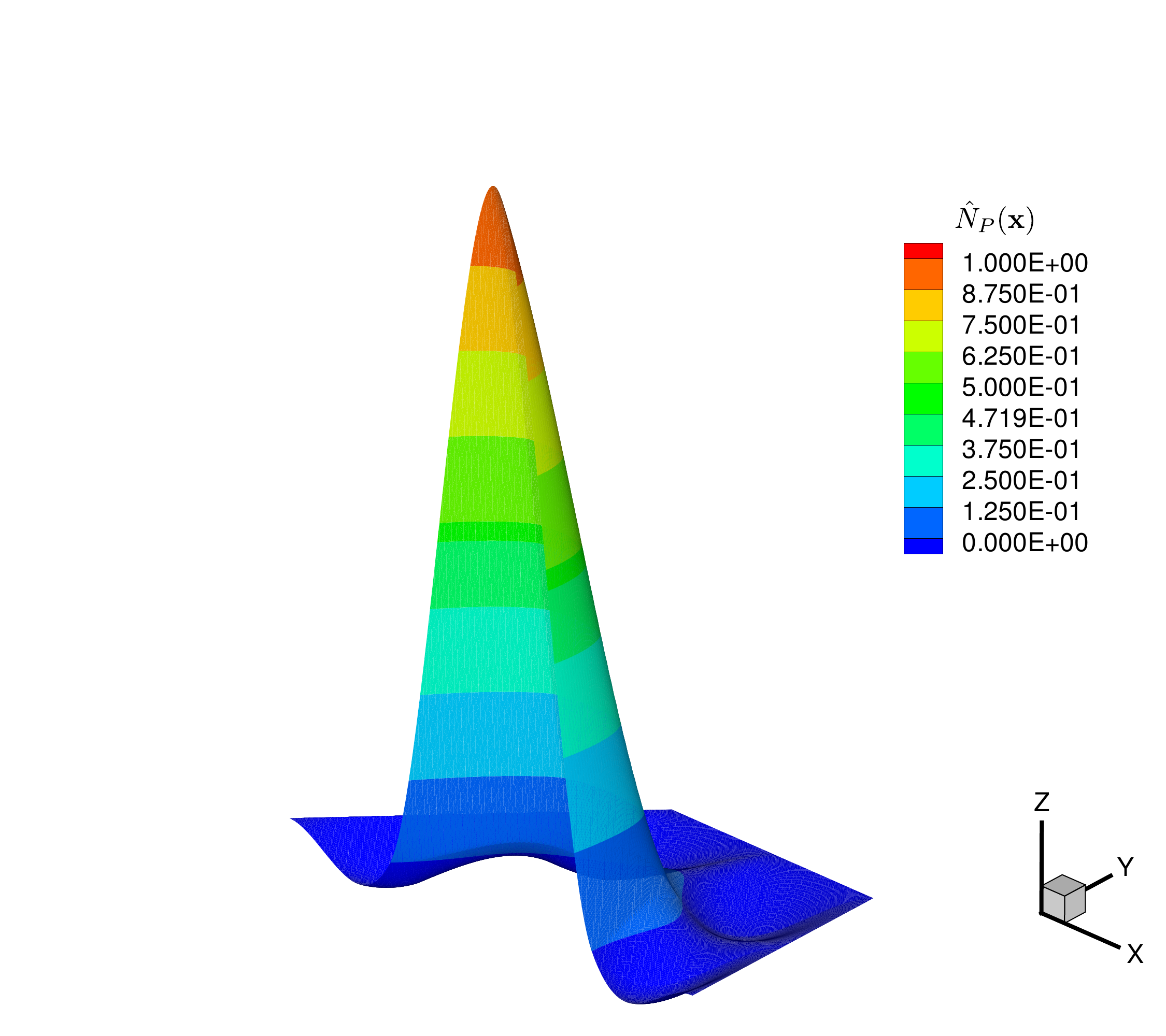}}
  \caption{The shape functions of DFEM in 2D}  \label{TFEMs}
\end{figure}

\subsection{Modification of the nodal gradients}
\begin{figure}
\centering
  {\includegraphics[width=0.5\textwidth]{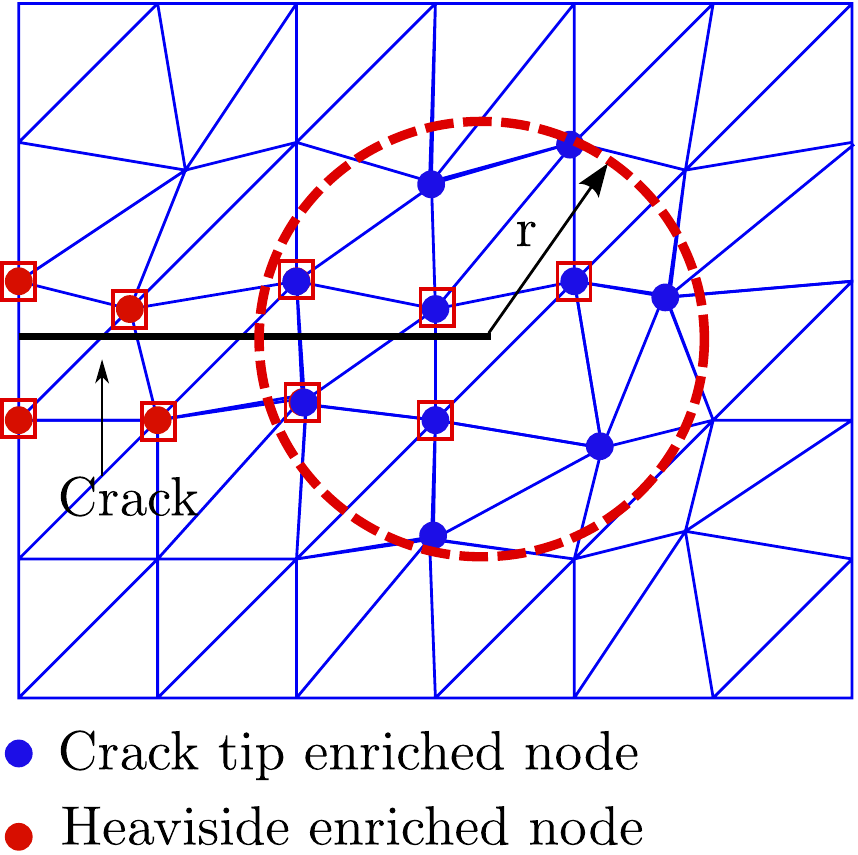}}
  \caption{nodal enrichment in XDFEM}
  \label{enrichD}
\end{figure}

When a $C^0$ node is needed, for instance on the nodes at material boundary, it is useful to modify
the calculation of the average nodal gradient as discussed below. The calculation of the nodal gradient can be performed as follows:
\begin{equation}
\bar{N}_{L,x}({\mathbf{x}_I})=N_{L,x}^{e}(\mathbf{x}_I)
\end{equation}
The right hand side is the derivative of $N_{L}$ computed in element $e$, in which the point of interest $\mathbf{x}$ is located. It can be easily observed that nodes at the endpoint of 1D bar automatically satisfy the above equation. Numerical implementation of XDFEM also indicates that all the enriched nodes in standard XFEM (the nodes circled by red boxes in Fig.\ref{enrichD}) have to be degenerated to $C^0$. One possible cause for this scenario relies on the fact that during the average calculation of gradients in equation (\ref{eq:14}), the contribution from split element cannot be computed directly from continuous FEM due to the discontinuity.

\subsection{The enriched 2D double-interpolation approximation}

\begin{figure}
\centering
  {\includegraphics[width=0.5\textwidth]{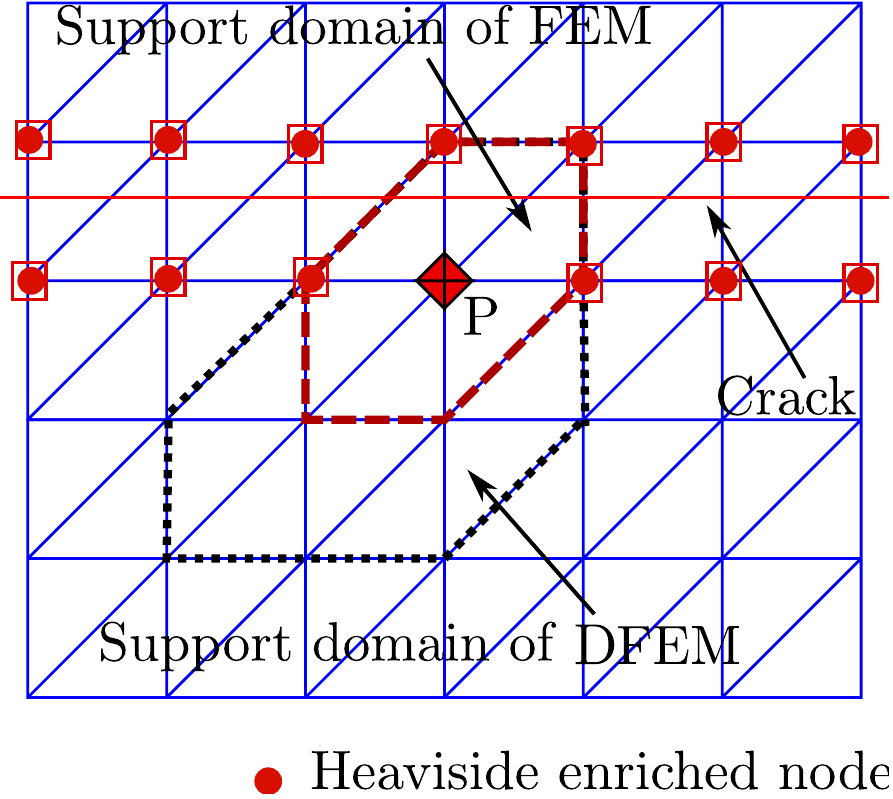}}
  \caption{The support domain of enriched DFEM. The support domain of the enriched nodes in XDFEM is decreased due to the nodal generation}
  \label{enrichshapdomain}
\end{figure}

\begin{figure}
\centering
  {\includegraphics[width=0.5\textwidth]{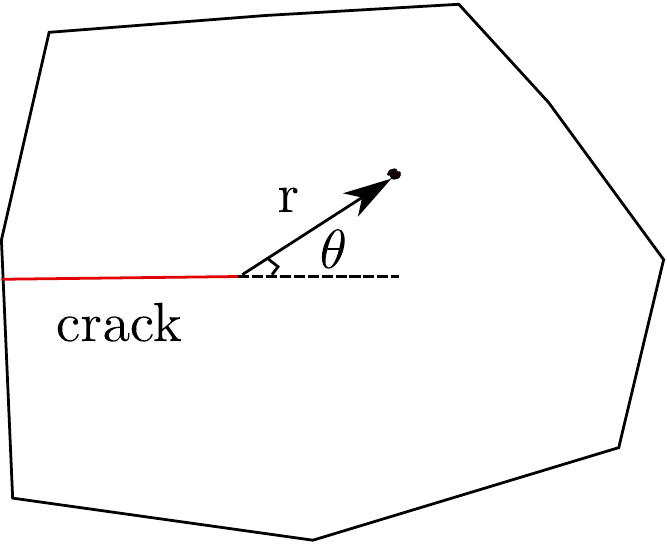}}
  \caption{Polar coordinate of crack tip}
  \label{polartip}
\end{figure}

The extended finite element method uses the local partition of unity which allows for the addition of a priori knowledge about the solution of boundary value problems into the approximation space of the numerical solution. The crack can be described in XFEM by enriching the standard displacement approximation as follows:

\begin{equation}
\mathbf{u}^h(\mathbf{x})=\sum_{I \in \mathscr{N}_I}\hat{N}_I(\mathbf{x})\mathbf{u}^I+
\sum_{J \in \mathscr{N}_J}\hat{N}_J(\mathbf{x})H(\mathbf{x})\mathbf{a}^J+
\sum_{K\in \mathscr{N}_K}\hat{N}_K(\mathbf{x})\sum_{\alpha=1}^{4}f_{\alpha}(\mathbf{x})\mathbf{b}^{K\alpha}
\end{equation}
where $\mathbf{u}^I$ are the regular DOFs. $\mathbf{a}^J$ are the additional Heaviside enriched DOFs. $\mathbf{b}^{K\alpha}$ are the additional crack tip enriched DOFs. $\mathscr{N}_I$,$\mathscr{N}_J$ and $\mathscr{N}_K$ are the collections of regular nodes, Heaviside enriched nodes and crack tip enriched nodes, respectively. $H(\cdot)$ is the Heaviside function. The crack tip enrichment function is defined as:

\begin{equation}
\left\{f_{\alpha}(r,\theta),\alpha=1,4\right\}=\left\{\sqrt{r}\text{sin}\frac{\theta}{2},\sqrt{r}\text{cos}\frac{\theta}{2},
\sqrt{r}\text{sin}\frac{\theta}{2}\text{sin}\theta,\sqrt{r}\text{cos}\frac{\theta}{2}\text{sin}\theta
\right\}
\end{equation}
where $(r,\theta)$ are the polar coordinates of the crack tip(Fig.\ref{polartip}). Fig.\ref{enrichshapdomain} compare the Heaviside enriched shape functions between XFEM and XDFEM.

\begin{figure}[htbp]
  \centering
  \subfigure[XDFEM]{\includegraphics[width=0.49\textwidth]{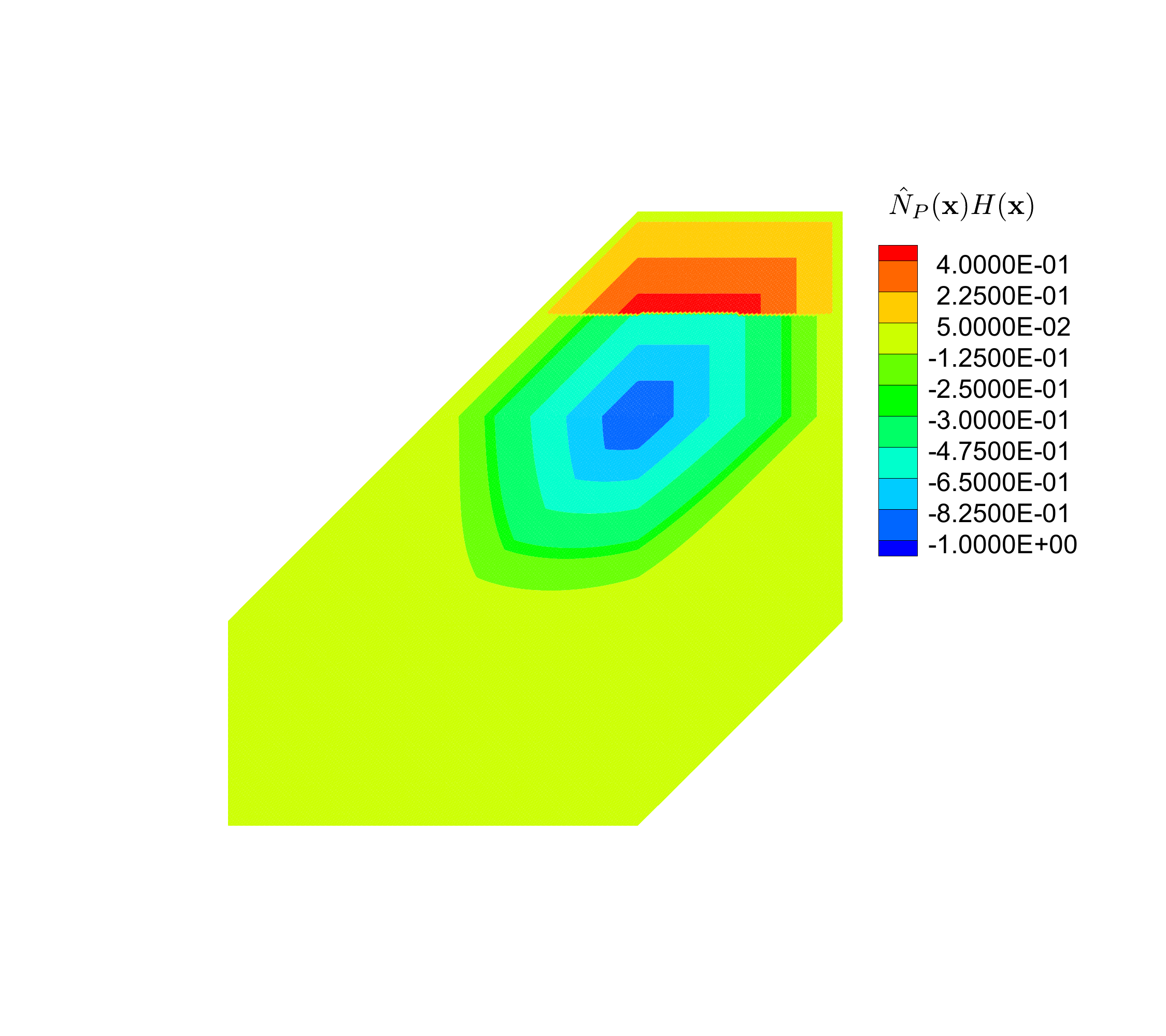}}
  \subfigure[XDFEM]{\includegraphics[width=0.49\textwidth]{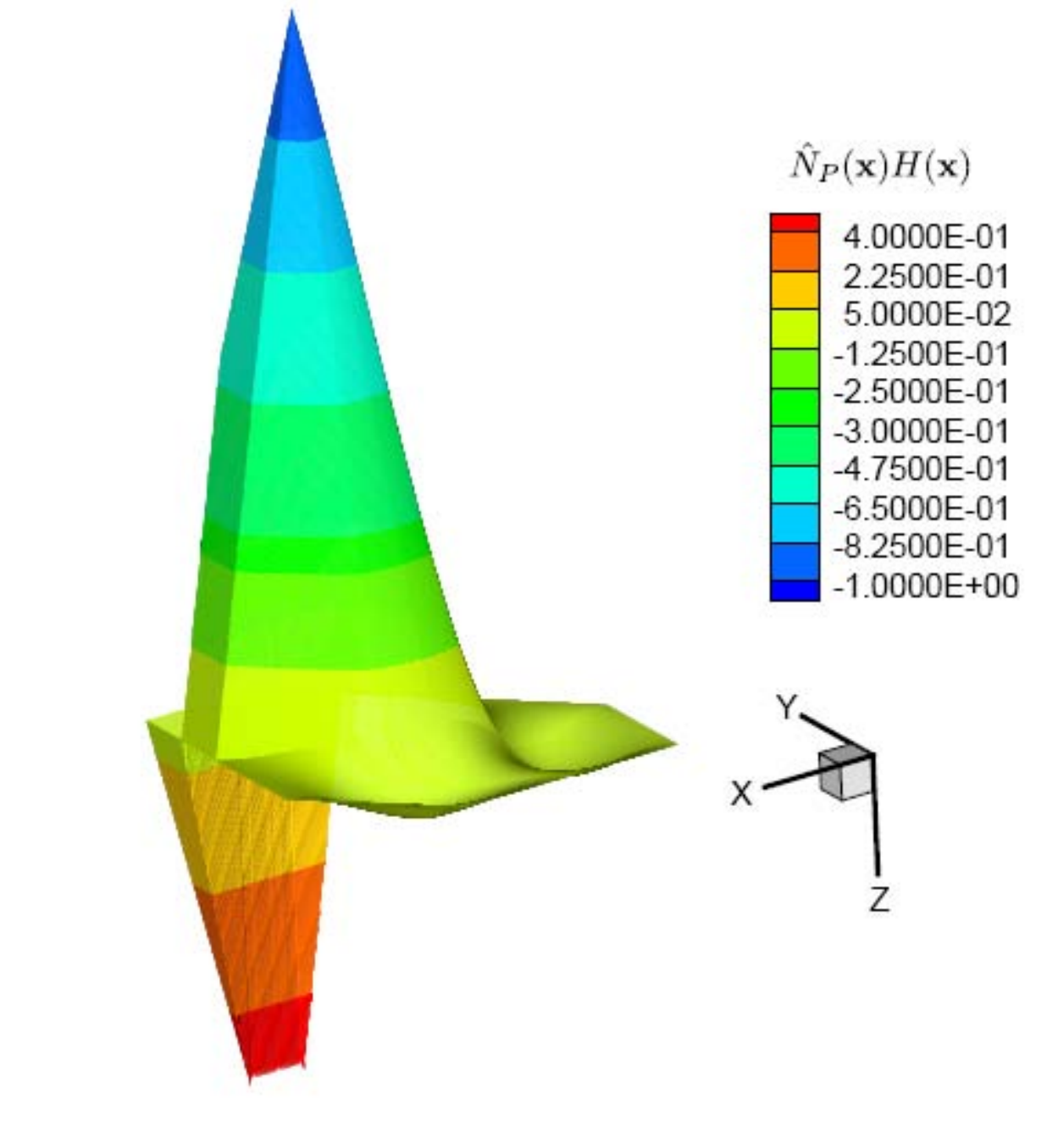}}
  \subfigure[XFEM]{\includegraphics[width=0.49\textwidth]{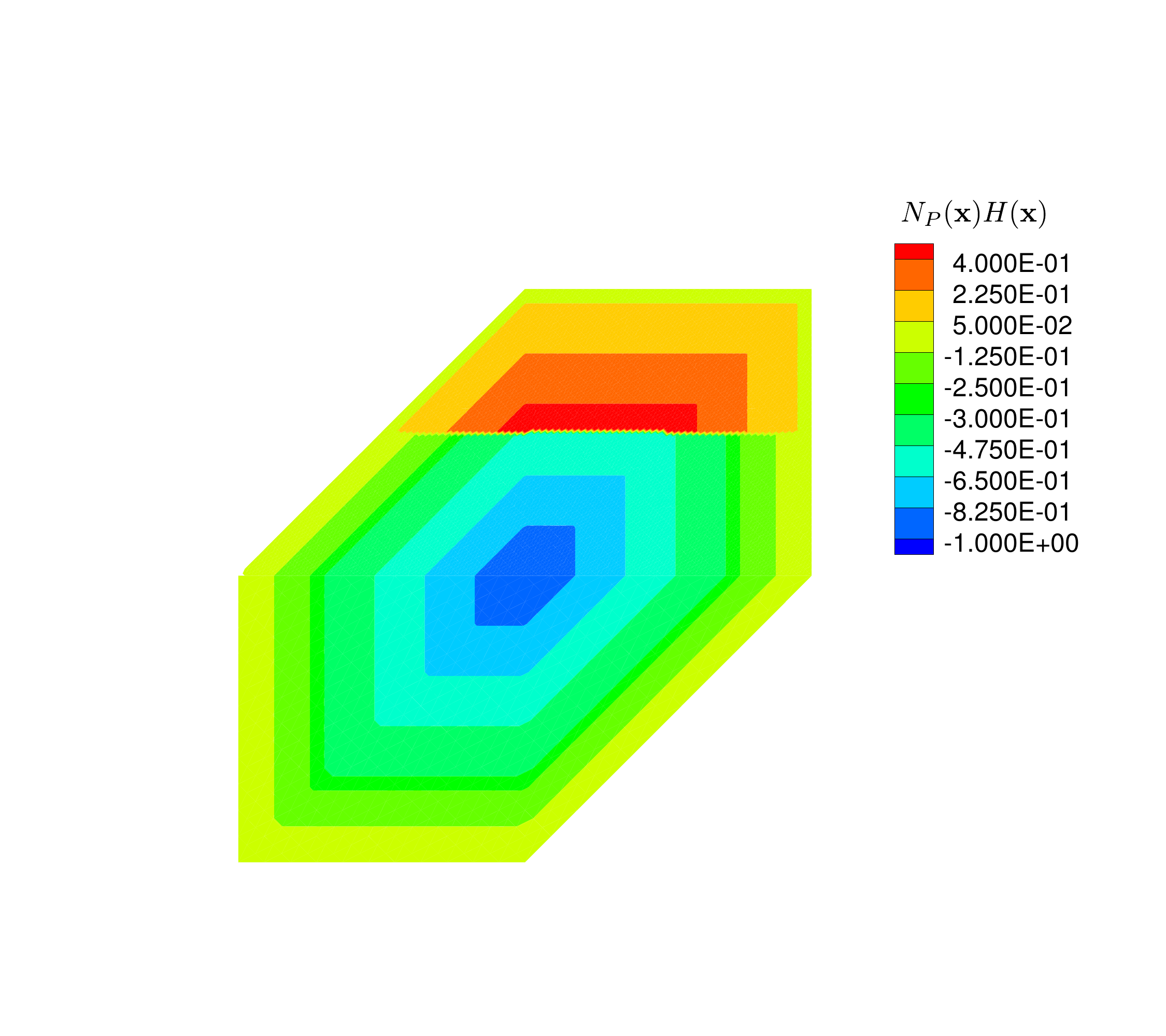}}
  \subfigure[XFEM]{\includegraphics[width=0.49\textwidth]{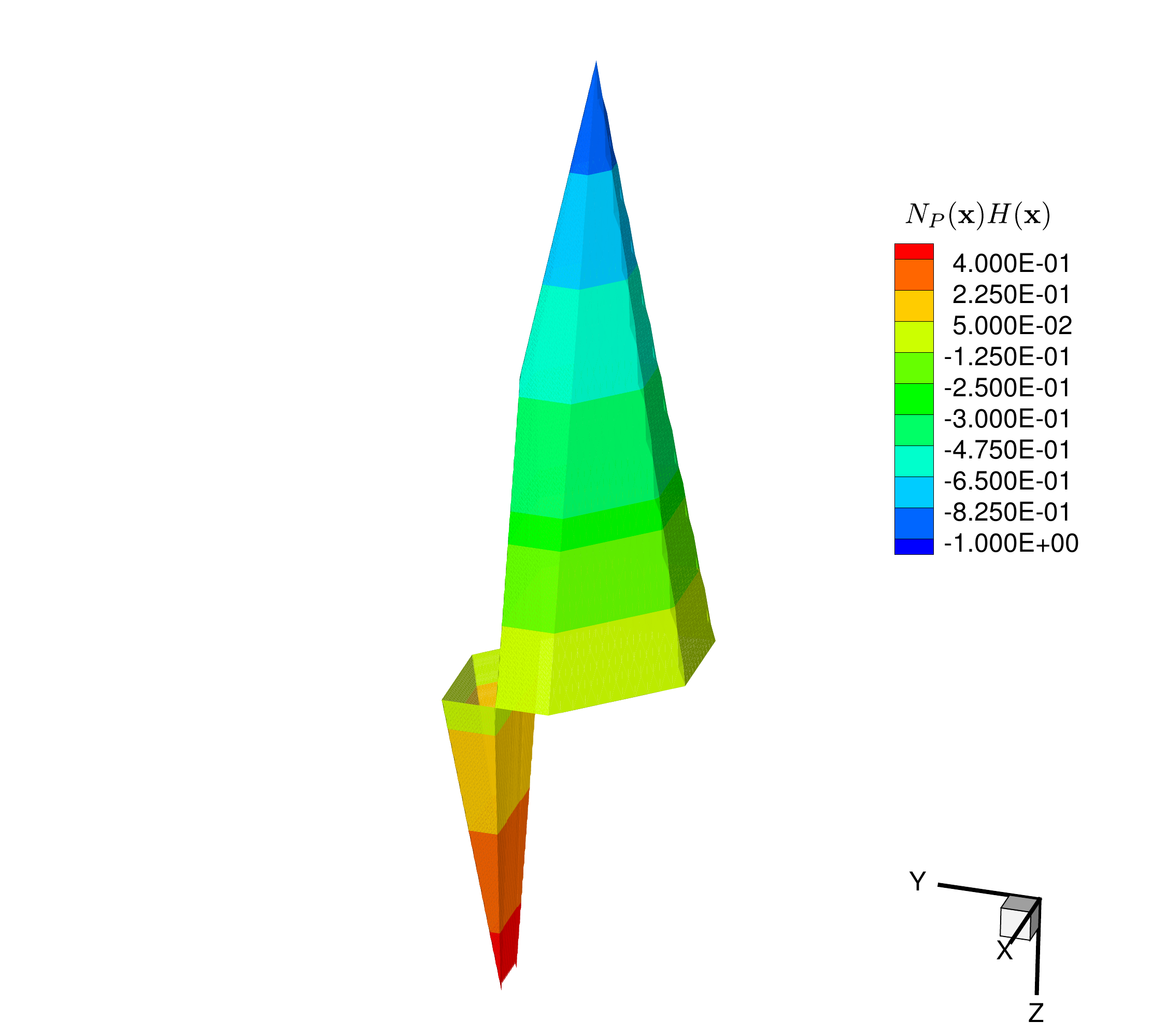}}
  \caption{Contour plot of Heaviside enriched shape functions}
  \label{enrichshape}
\end{figure}

\section{Weak form and discretized formulations}
For elastostatic problems with the assumptions of small strains and small displacements, the equilibrium equations and boundary conditions are given by:

\begin{equation}
\begin{array}{cccc}
\nabla\cdot\boldsymbol{\sigma}=
\boldsymbol{0}\qquad & \textrm{on}\quad &\Omega \\
\boldsymbol{\sigma}\cdot\mathbf{n}=\bar{\mathbf{t}}\qquad & \textrm{on}\quad &\Gamma_t \\
\mathbf{u}=\mathbf{\bar{\mathbf{u}}}\qquad &\textrm{on}\quad &\Gamma_u \\
\end{array}
\end{equation}

Here $\bar{\mathbf{t}}$ is the traction imposed on boundary $\Gamma_t$. Further, assuming traction free crack faces:

\begin{equation}
\boldsymbol{\sigma}\cdot\mathbf{n}=\boldsymbol{0} \qquad\textrm{on}\quad\Gamma_{{cr}^+},\Gamma_{{cr}^-}
\end{equation}

where $\Gamma_{{cr}^+},\Gamma_{{cr}^-}$ are the upper and lower crack surfaces respectively. The strain-displacement relation and the constitutive law are respectively as:
\begin{subequations}
\begin{equation}
\boldsymbol{\epsilon}=\frac{1}{2}\left(\nabla \mathbf{u}+({\nabla \mathbf{u}})^{\textrm{T}}\right)
\end{equation}
\begin{equation}
\boldsymbol{\sigma}=\mathbf{C}:\boldsymbol{\epsilon}
\end{equation}
\end{subequations}

Substituting the trial and test functions into the weak form, and using the arbitrary nodal variations, the discretized equations can be simplified as,
\begin{equation}
\mathbf{Kd}=\mathbf{f}
\end{equation}
where, $\mathbf{d}$ is the nodal vector and $\mathbf{K}$ is the stiffness matrix. The elemental form of $\mathbf{K}$ is given by;
\begin{equation} \label{eq:11}
\mathbf{K}^e_{IJ}=
\left[\begin{array}{cccc}
\mathbf{K}_{IJ}^{uu} & \mathbf{K}_{IJ}^{ua} & \mathbf{K}_{IJ}^{ub} \\
\mathbf{K}_{IJ}^{au} & \mathbf{K}_{IJ}^{aa} & \mathbf{K}_{IJ}^{ab} \\
\mathbf{K}_{IJ}^{bu} & \mathbf{K}_{IJ}^{ba} & \mathbf{K}_{IJ}^{bb} \\
\end{array}\right]
\end{equation}

The external force vector $\mathbf{f}$ is defined as
\begin{equation} \label{eq:12}
\mathbf{f}_{I}=\{\mathbf{f}^u_I \quad \mathbf{f}^a_I \quad \mathbf{f}^{b^1}_I \quad \mathbf{f}^{b^2}_I \quad \mathbf{f}^{b^3}_I \quad \mathbf{f}^{b^4}_I\}.
\end{equation}
The submatrices and vectors in equations (\ref{eq:11}) and (\ref{eq:12}) are presented as follows:
\begin{subequations}
\begin{equation}\label{eq:13}
\hspace{5mm}  \mathbf{K}^{rs}_{IJ}=\int_{\Omega^e}(\mathbf{B}_I^r)^T\mathbf{C}\mathbf{B}_J^s \textrm{d}\Omega
\qquad (r,s=u,a,b)
\end{equation}
\begin{equation}
\mathbf{f}_I^u=\int_{\partial \Omega_t^h\cap \partial \Omega^e} \hat{N}_I \bar{\mathbf{t}}\textrm{d}\Gamma
\end{equation}
\begin{equation}
\mathbf{f}_I^a=\int_{\partial \Omega_t^h\cap \partial \Omega^e} \hat{N}_I H \bar{\mathbf{t}}\textrm{d}\Gamma
\end{equation}
\begin{equation}
\mathbf{f}_I^{b^\alpha}=\int_{\partial \Omega_t^h\cap \partial \Omega^e} \hat{N}_I f_\alpha \bar{\mathbf{t}}\textrm{d}\Gamma
\qquad (\alpha=1, 2, 3, 4)
\end{equation}
\end{subequations}

In equations (\ref{eq:13}), $\mathbf{B}^u_I,\mathbf{B}^a_I$ and $\mathbf{B}^{b^\alpha}_I$ are given by
\begin{subequations}
\begin{equation}
\mathbf{B}^u_I=
\left[\begin{array}{cccc}
\hat{N}_{I,x} & 0      \\
0     & \hat{N}_{I,y}  \\
\hat{N}_{I,y} & \hat{N}_{I,x}
\end{array}\right],
\end{equation}
\begin{equation}
\mathbf{B}^a_I=
\left[\begin{array}{cccc}
{(\hat{N}_I (H-H_I))}_{,x} & 0               \\
0              & {(\hat{N}_I (H-H_I))}_{,y}  \\
{(\hat{N}_I (H-H_I))}_{,y} & {(\hat{N}_I (H-H_I)}_{,x}
\end{array}\right]
\end{equation}
\begin{equation}
\mathbf{B}^b_I=\left[\mathbf{B}^{b^1}_I\quad \mathbf{B}^{b^2}_I\quad\mathbf{B}^{b^3}_I\quad\mathbf{B}^{b^4}_I\right]
\end{equation}
\begin{equation}
\mathbf{B}^{b^\alpha}_I=
\left[\begin{array}{cccc}
{(\hat{N}_I (f_\alpha-f_{\alpha I}))}_{,x} & 0                         \\
0                        & {(\hat{N}_I (f_\alpha-f_{\alpha I}))}_{,y}  \\
{(\hat{N}_I (f_\alpha-f_{\alpha I}))}_{,y} & {(\hat{N}_I (f_\alpha-f_{\alpha I}))}_{,x}
\end{array}\right]\quad(\alpha=1-4)
\end{equation}
\end{subequations}

In order to obtain the nodal displacements in a more straightforward manner, the shifted-basis is adopted in the above equations.

\section{Numerical examples}

A set of numerical examples are performed to compare the efficiency of the double-interpolation and its enriched form. In order to assess the convergence rate of each model, the relative error of the $L_2$ displacement norm and $H_1$ energy norm are defined, respectively, as:

\begin{subequations}
\begin{equation}
R_d=\sqrt{\frac{\int_\Omega(\mathbf{u}^h-\mathbf{u})^{\text{T}}(\mathbf{u}^h-\mathbf{u})\text{d}\Omega}
{\int_\Omega\mathbf{u}^{\text{T}}\mathbf{u}\text{d}\Omega}}
\end{equation}
\begin{equation}
R_e=\sqrt{\frac{\int_\Omega(\boldsymbol{\sigma}^h-\boldsymbol{\sigma})^{\text{T}}(\boldsymbol{\epsilon}^h-
\boldsymbol{\epsilon})\text{d}\Omega}
{\int_\Omega\boldsymbol{\sigma}^{\text{T}}\boldsymbol{\epsilon}\text{d}\Omega}}
\end{equation}
\end{subequations}

where, the fields with superscript $"h"$ refer to the approximation, and $\boldsymbol{\sigma}, \boldsymbol{\epsilon}, \mathbf{u}$ are exact fields.
Unless specified otherwise, the \emph{Young's modulus} $E$ and \emph{Possion's ratio} $\nu$ are assumed to be $1000$ and $0.3$ respectively. The constants $\mu$ and $\kappa$ are given by;

\begin{subequations}
\begin{equation}
\mu=\frac{E}{2(1+\nu)}
\end{equation}
\begin{displaymath}
\kappa=\left\{
\begin{array}{ll}
3-4\nu, & \text{Plane strain} \\
(1-\nu)/(3+\nu), & \text{Plane stress}  \\
\end{array}\right.
\end{displaymath}
\end{subequations}

\subsection{Infinite plate with a center hole}
\begin{figure}
\center
  \subfigure[]{\includegraphics[width=0.49\textwidth]{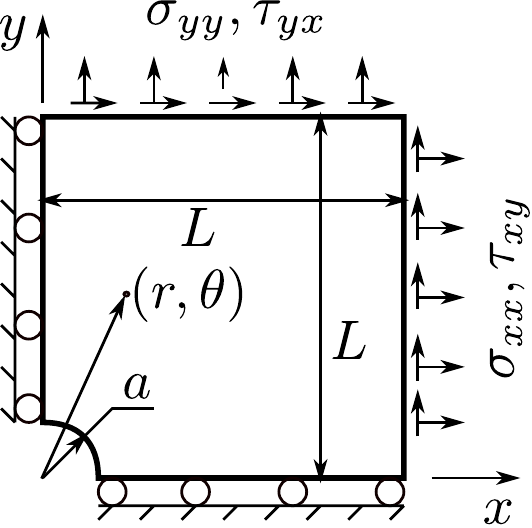}}
  \subfigure[]{\includegraphics[width=0.49\textwidth]{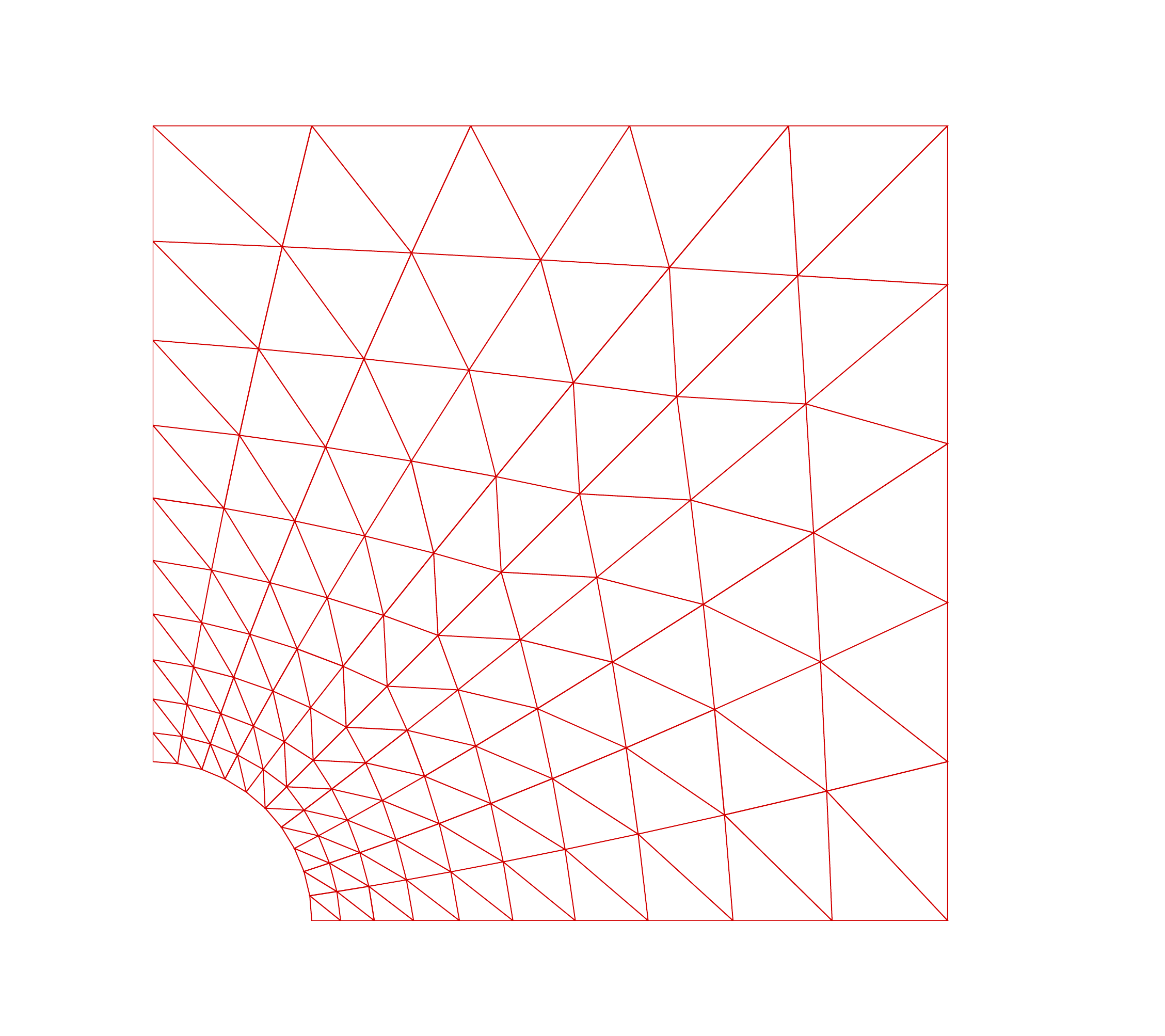}}
  \caption{a:$1/4$ model of the infinite plate with a center hole; b: The typical mesh division}
  \label{modplate}
\end{figure}

Fig.\ref{modplate} presents the upper right quadrant model of an infinite plate with a center hole subjected to remote tensile loads. Here the geometry parameters $L=5$ and $a=1$. The analytical solutions for stress and displacement fields are given as \cite{timoshenko1972theory}:
\begin{subequations}
\begin{equation}
\sigma_{xx}(r,\theta)=1-\frac{a^2}{r^2}\left(\frac{3}{2}\text{cos}2\theta+\text{cos}4\theta\right)+
\frac{3a^4}{2r^4}\text{cos}4\theta
\end{equation}
\begin{equation}
\sigma_{yy}(r,\theta)=-\frac{a^2}{r^2}\left(\frac{1}{2}\text{cos}2\theta-\text{cos}4\theta\right)-
\frac{3a^4}{2r^4}\text{cos}4\theta
\end{equation}
\begin{equation}
\tau_{xy}(r,\theta)=1-\frac{a^2}{r^2}\left(\frac{1}{2}\text{sin}2\theta+\text{sin}4\theta\right)+
\frac{3a^4}{2r^4}\text{sin}4\theta
\end{equation}
\begin{equation}
u_r(r,\theta)=\frac{a}{8\mu}\left[\frac{r}{a}(\kappa+1)\text{cos}\theta+\frac{2a}{r}((1+\kappa)
\text{cos}\theta+\text{cos}3\theta)
-\frac{2a^3}{r^3}\text{cos}3\theta\right]
\end{equation}
\begin{equation}
u_{\theta}(r,\theta)=\frac{a}{8\mu}\left[\frac{r}{a}(\kappa-1)\text{sin}\theta+\frac{2a}{r}((1-\kappa)
\text{sin}\theta+\text{sin}3\theta)
-\frac{2a^3}{r^3}\text{sin}3\theta\right]
\end{equation}
\end{subequations}
where $(r,\theta)$ are the polar coordinates. The exact stress is imposed on the top and right boundary of the model. The number of nodes used in the four models are $121$, $441$, $1681$ and $6561$.

\begin{figure}[htbp]
  \centering
  \subfigure{\includegraphics[width=1.\textwidth]{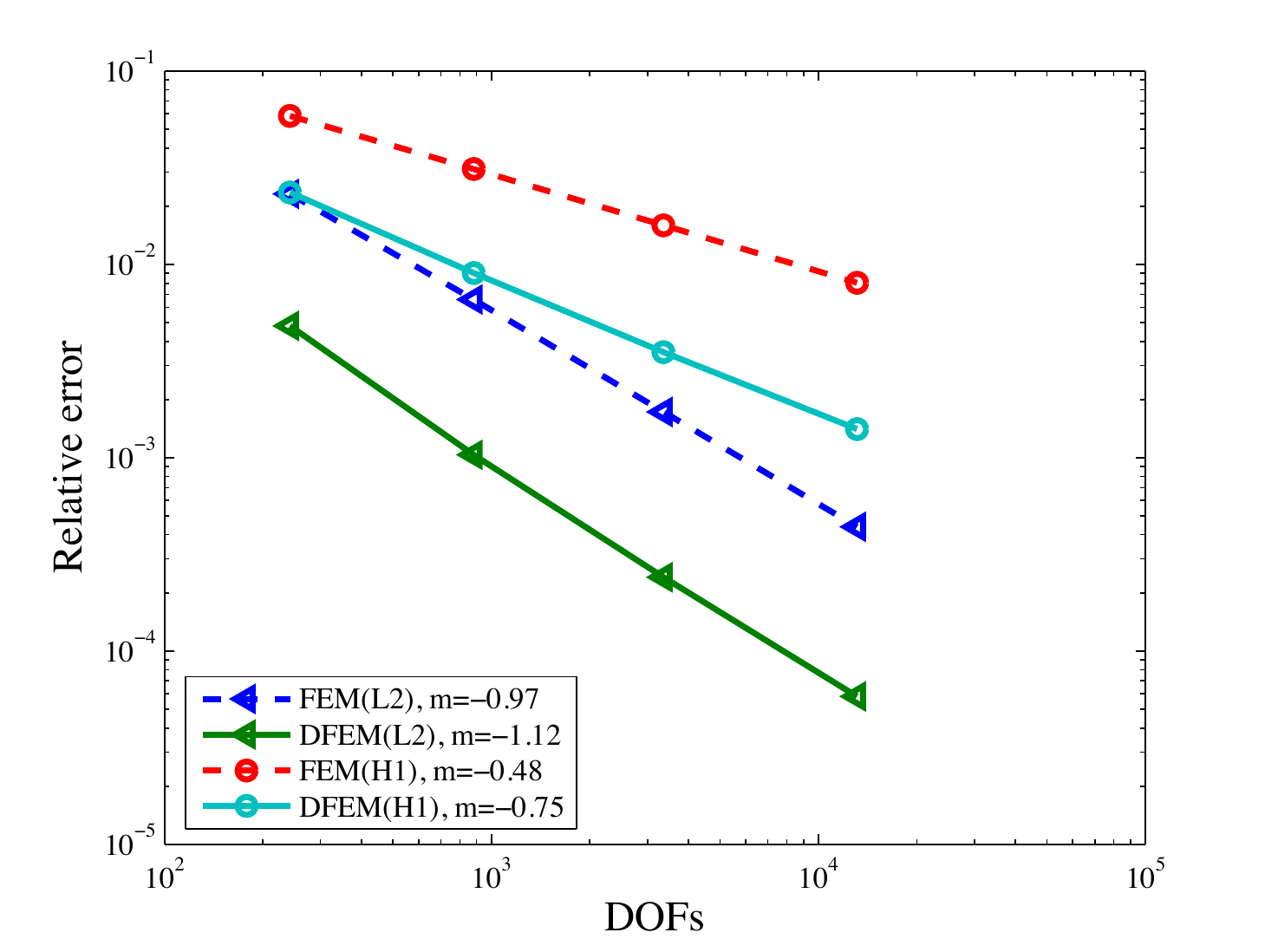}}
  \caption{Relative error in displacement and energy norm of in infinite plate with a hole}
  \label{plhole}
\end{figure}

In this example, the numerical results obtained using DFEM and 3-noded triangle (T3) finite element are compared for the same mesh model. The relative error in displacement and energy norm for this example are plotted in Fig.\ref{plhole}. The DFEM shows a significantly better precision and higher convergence rate than the standard 3-noded FEM.

\subsection{Timoshenko Beam}
\begin{figure}
\center
{\includegraphics[width=1.\textwidth]{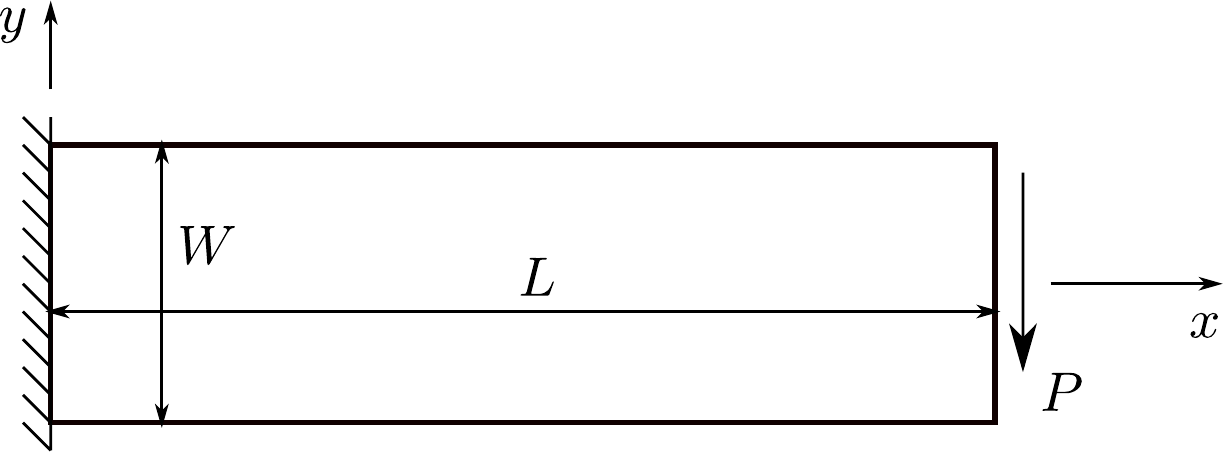}}
  \caption{Physical model of cantilever beam }
  \label{modbeam}
\end{figure}

Fig.\ref{modbeam} illustrates the physical model of a cantilever beam. In this example, plane stress conditions are assumed. The geometric parameters are taken as $L=48$ and $W=12$. The analytical solution for displacement and stress fields is given in \cite{timoshenko1972theory} as:

\begin{subequations}
\begin{equation}
u_x(x,y)=\frac{Py}{6EI}\left[(6L-3x)x+(2+\nu)(y^2-\frac{W^2}{4})\right]
\end{equation}
\begin{equation}
u_y(x,y)=-\frac{P}{6EI}\left[3\nu y^2(L-x)+(4+5\nu)\frac{W^2x}{4}+(3L-x)x^2\right]
\end{equation}
\begin{equation}
\sigma_{xx}(x,y)=\frac{P(L-x)y}{I}
\end{equation}
\begin{equation}
\sigma_{yy}(x,y)=0
\end{equation}
\begin{equation}
\tau_{xy}(x,y)=-\frac{P}{2I}\left(\frac{W^2}{4}-y^2\right)
\end{equation}
\end{subequations}
where $P=1000.$ and $I=W^3/12$. The exact displacement is applied to the left boundary and the exact traction is applied to the right boundary.

\begin{figure}
\center
{\includegraphics[width=1.\textwidth]{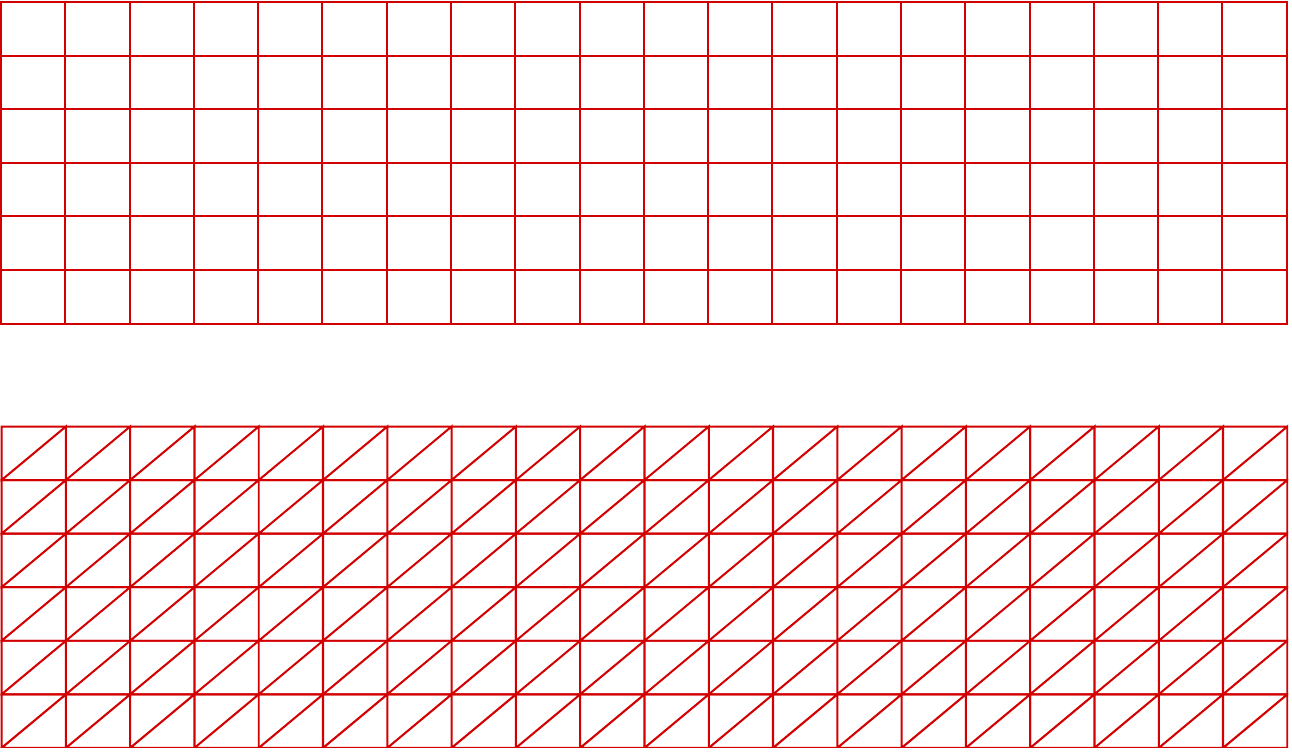}}
  \caption{Mesh discretization using regular quadrilateral and triangular element }
  \label{mesh3beam}
\end{figure}

\begin{figure}[htbp]
  \centering
  \subfigure[]{\includegraphics[width=1.\textwidth]{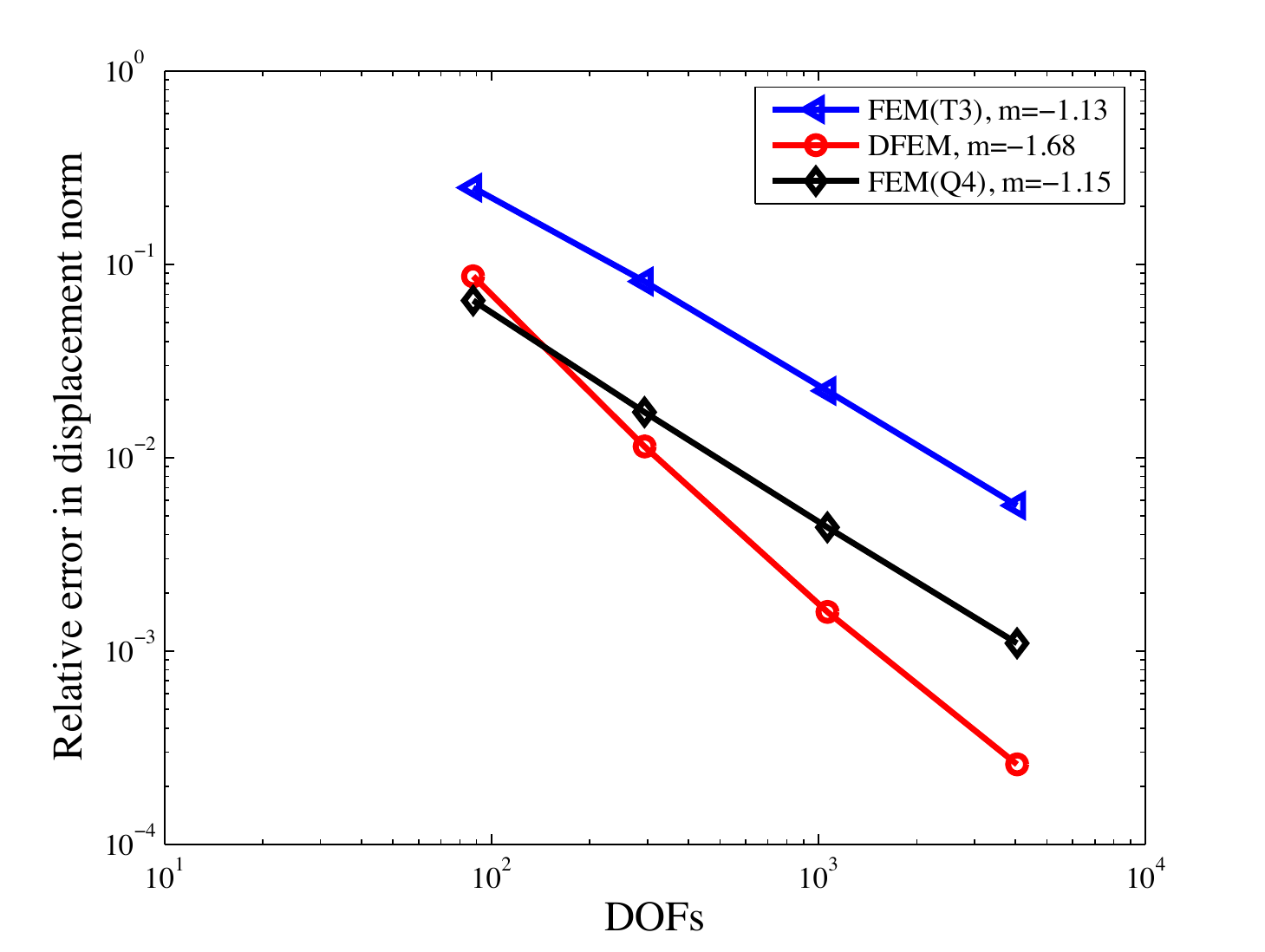}}
  \subfigure[]{\includegraphics[width=1.\textwidth]{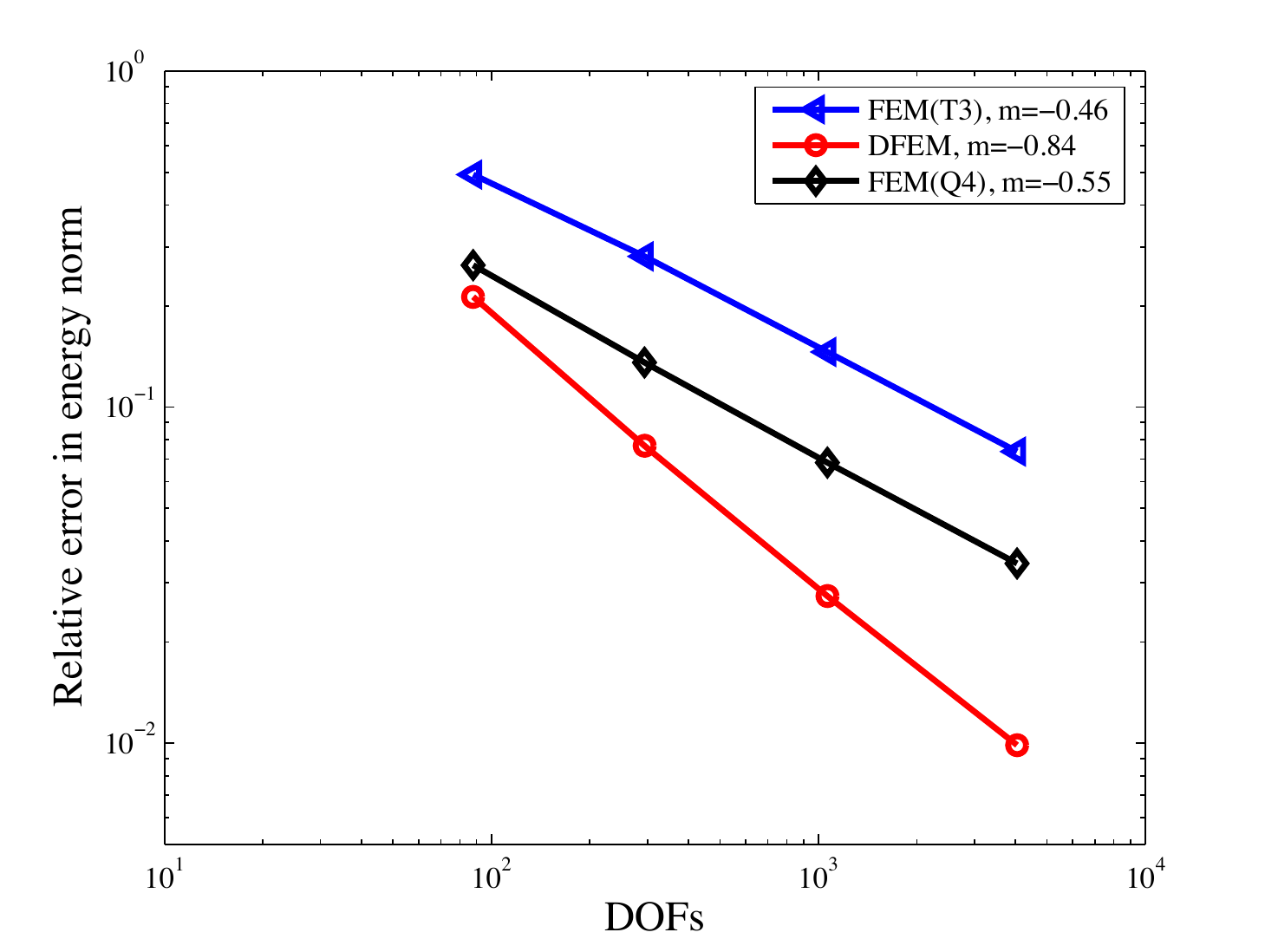}}
  \caption{Relative error in displacement and energy norm of Timoshenko beam}
  \label{beam}
\end{figure}

Structured meshes are used in this example to ensure regular node location and to enable easier comparison among the 3-noded triangle (T3) finite element, 4-noded quadrilateral (Q4) element and DFEM (Fig.\ref{mesh3beam}). Four mesh sizes, $3\times 10$, $6\times 20$, $12\times 40$ and $24\times 80$ , are used. It can be verified that the bilinear Q4 FE is more accurate than the constant strain T3 FE. However, their convergence rates do not demonstrate any significant difference as is known from a priori error estimation. The DFEM solution, however, demonstrates better accuracy and super-convergence compared to Q4 and T3. It can be concluded that for the same number of DOFs, DFEM is more accurate than Q4 FE.

\subsection{Griffith crack}

\begin{figure}
  {\includegraphics[width=1.\textwidth]{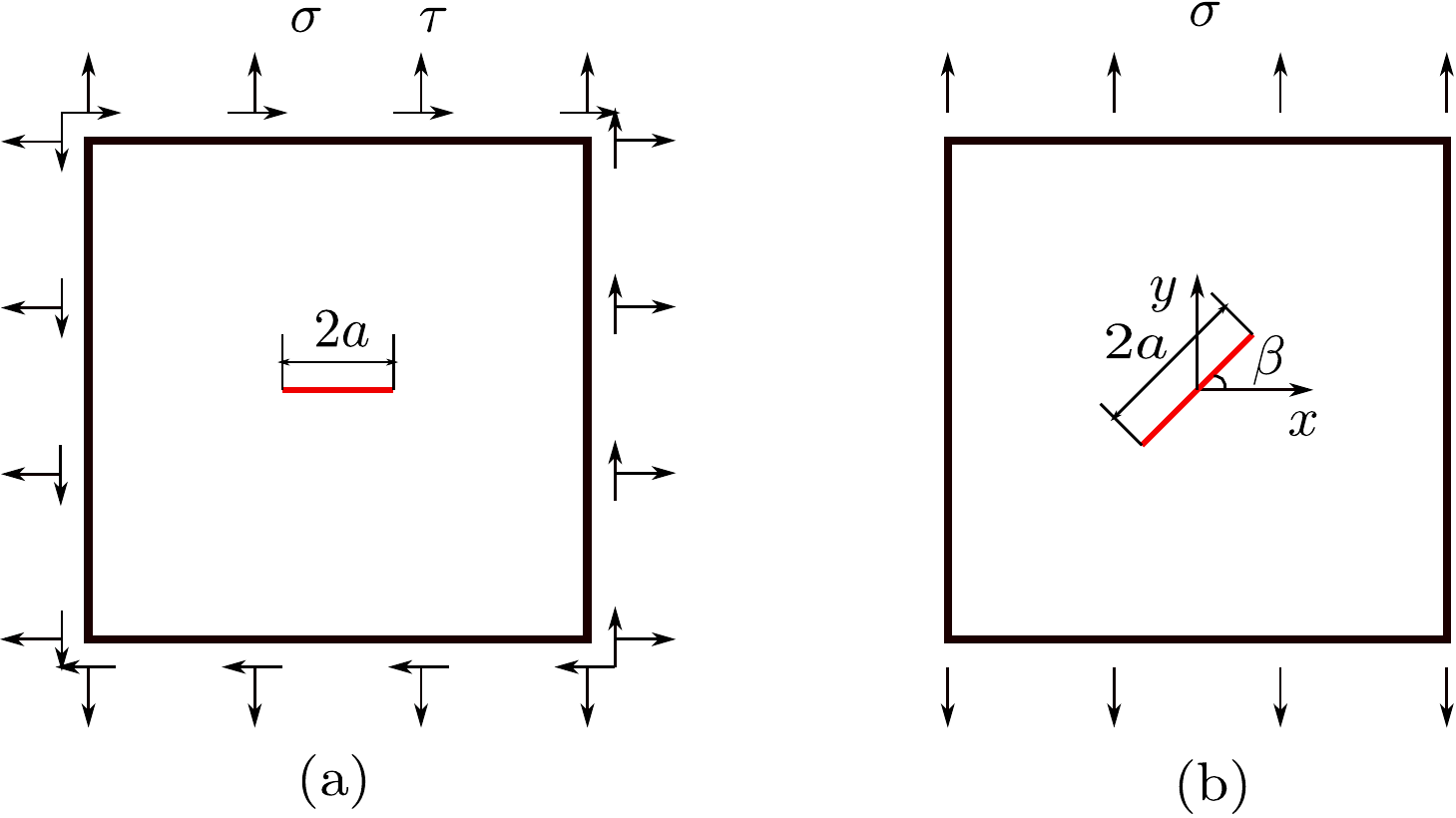}}
  \caption{(a)Griffith crack; (b) inclined crack}
  \label{griff}
\end{figure}

A Griffith crack problem is shown in Fig.\ref{griff}(a). An infinite plate with a crack  segment (a=1.) subjected to remote tensile loads is considered here. A square domain ($10\times 10$) is selected in the vicinity of crack tip. The analytical displacement and stress field near the crack tip are given by \cite{Westergaard1939}:

\begin{subequations}
\begin{equation}
\begin{aligned}
\sigma_{xx}(r,\theta)=&\frac{K_I}{\sqrt{2\pi r}}\text{cos}\frac{\theta}{2}\left(1-
\text{sin}\frac{\theta}{2}\text{sin}\frac{3\theta}{2}\right) \\
&
-\frac{K_{II}}{\sqrt{2\pi r}}\text{sin}\frac{\theta}{2}\left(2+
\text{cos}\frac{\theta}{2}\text{cos}\frac{3\theta}{2}\right)
\end{aligned}
\end{equation}
\begin{equation}
\begin{aligned}
\sigma_{yy}(r,\theta)=\frac{K_I}{\sqrt{2\pi r}}\text{cos}\frac{\theta}{2}\left(1+
\text{sin}\frac{\theta}{2}\text{sin}\frac{3\theta}{2}\right)
+\frac{K_{II}}{\sqrt{2\pi r}}\text{sin}\frac{\theta}{2}\text{cos}\frac{\theta}{2}
\text{cos}\frac{3\theta}{2}
\end{aligned}
\end{equation}
\begin{equation}
\begin{aligned}
\tau_{xy}(r,\theta)=\frac{K_I}{\sqrt{2\pi r}}\text{sin}\frac{\theta}{2}\text{cos}\frac{\theta}{2}
\text{cos}\frac{3\theta}{2}+\frac{K_{II}}{\sqrt{2\pi r}}\text{cos}\frac{\theta}{2}\left(1-
\text{sin}\frac{\theta}{2}\text{sin}\frac{3\theta}{2}\right)
\end{aligned}
\end{equation}
\begin{equation}
\begin{aligned}
u_x(r,\theta)=&\frac{K_I}{2\mu}\sqrt{\frac{r}{2\pi}}
\text{cos}\frac{\theta}{2}\left(\kappa-1+2\text{sin}^2\frac{\theta}{2}\right)\\
&
+\frac{(1+\nu)K_{II}}{E}\sqrt{\frac{r}{2\pi}}\text{sin}\frac{\theta}{2}\left(\kappa+1
+2\text{cos}^2\frac{\theta}{2}\right)
\end{aligned}
\end{equation}
\begin{equation}
\begin{aligned}
u_y(r,\theta)&=\frac{K_I}{2\mu}\sqrt{\frac{r}{2\pi}}
\text{sin}\frac{\theta}{2}\left(\kappa+1-2\text{cos}^2\frac{\theta}{2}\right)\\
&\quad
+\frac{(1+\nu)K_{II}}{E}\sqrt{\frac{r}{2\pi}}\text{cos}\frac{\theta}{2}\left(1-\kappa
+2\text{sin}^2\frac{\theta}{2}\right)
\end{aligned}
\end{equation}
\end{subequations}
where $K_I$ and $K_{II}$ are the stress intensity factors (SIFs) for mode-I and mode-II, respectively. $(r,\theta)$ are the polar coordinates used to define the crack geometry.

\subsubsection{Convergence study}

\begin{figure}[htbp]
  \centering
  \subfigure[]{\includegraphics[width=1.\textwidth]{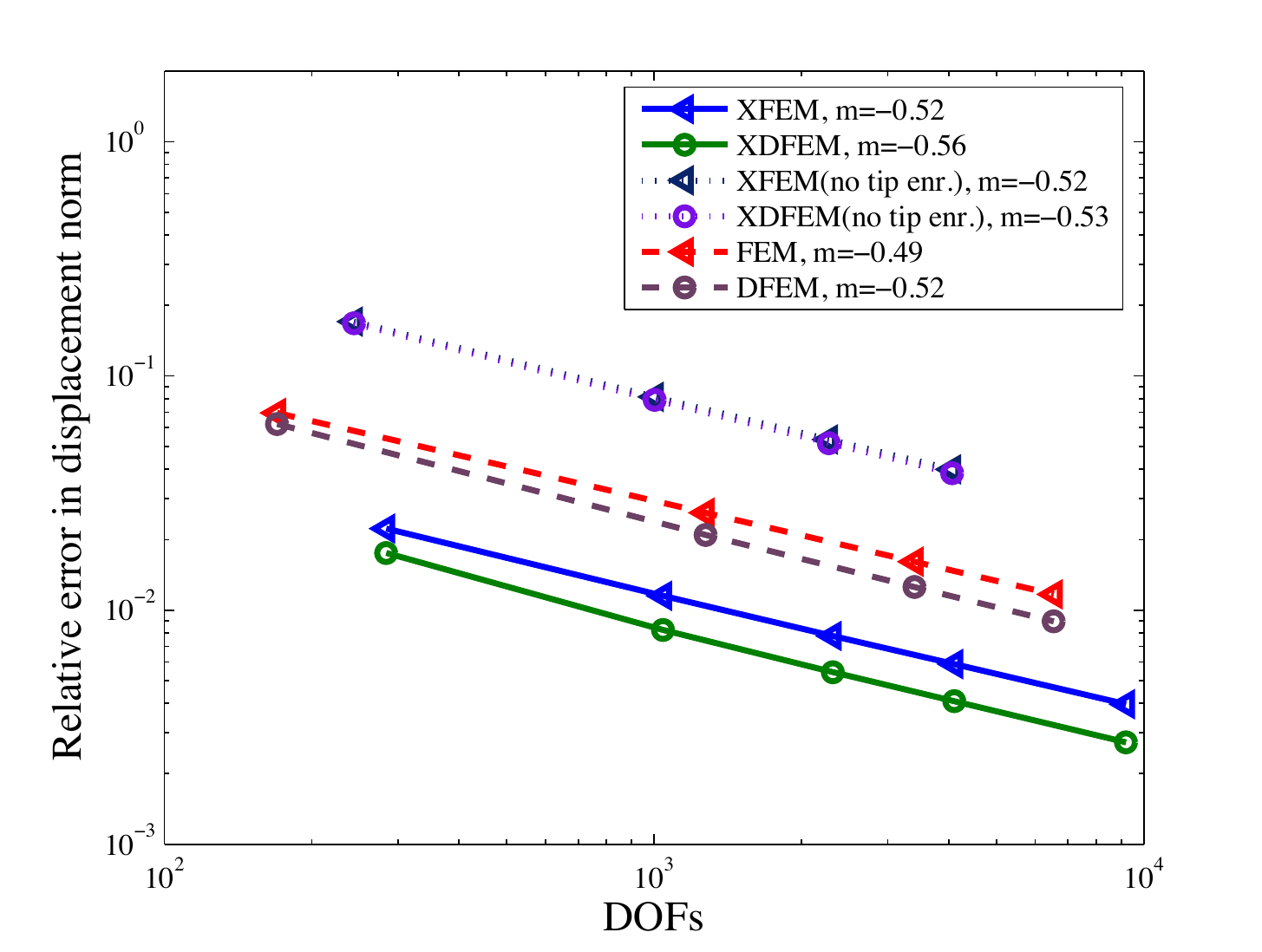}}
  \subfigure[]{\includegraphics[width=1.\textwidth]{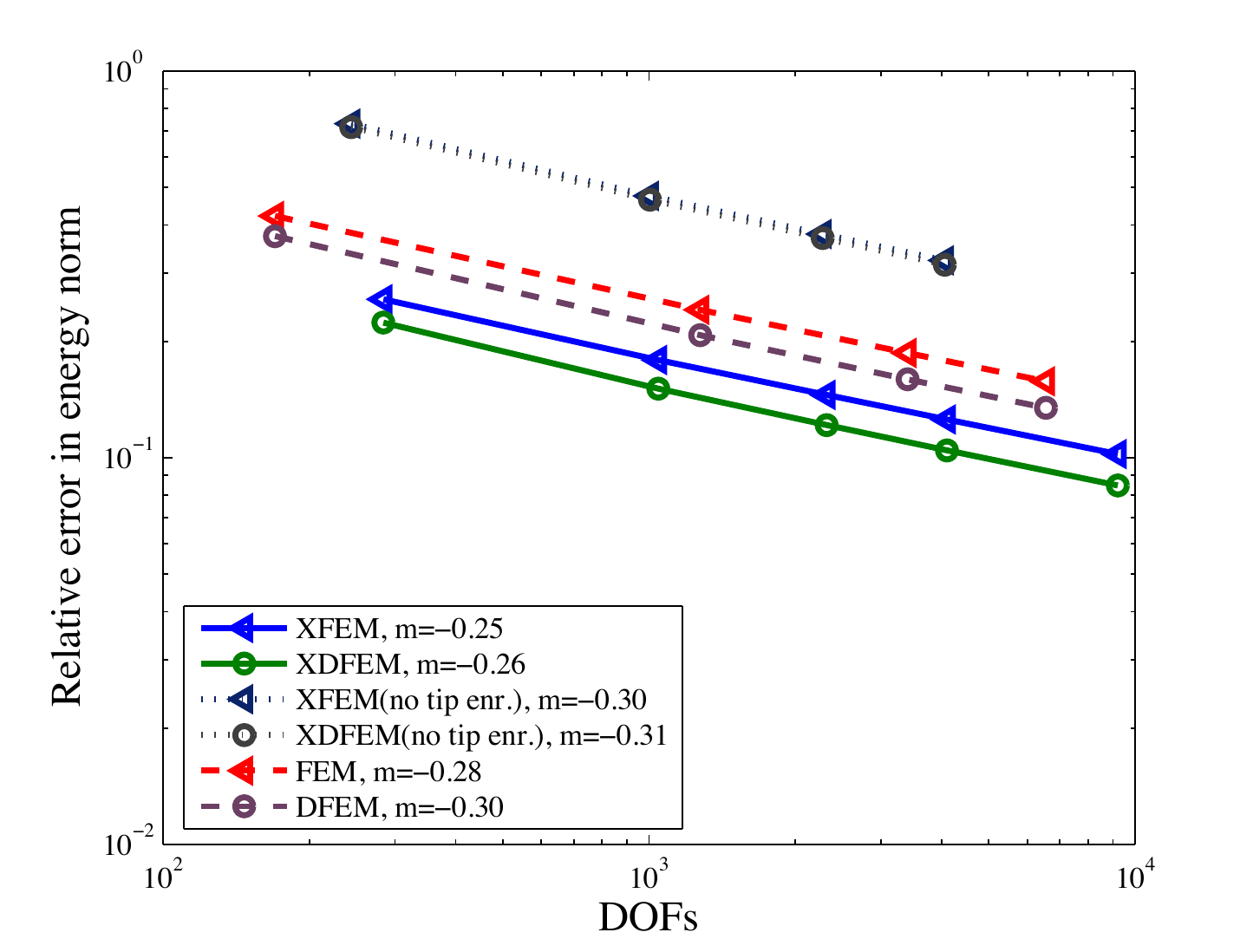}}
  \caption{Relative error in displacement and energy norm of in Griffith crack}
  \label{cra}
\end{figure}

Griffith crack problem is used here to investigate the enrichment effects of DFEM. The convergence rate in XDFEM is studied in three folds: explicit crack model, heaviside enrichment only and full enrichment model. These results are plotted in Fig.\ref{cra}. From Fig.\ref{cra}, it can be concluded that the DFEM yields better accuracy and slightly improves the convergence rate compared to FEM for all the cases considered. It also transpires from the results that the full enrichment of DFEM produces better accuracy than modelling the crack explicitly.

\begin{figure}[htbp]
  \centering
  \subfigure[Relative error in displacement and energy norm]{\includegraphics[width=1.\textwidth]{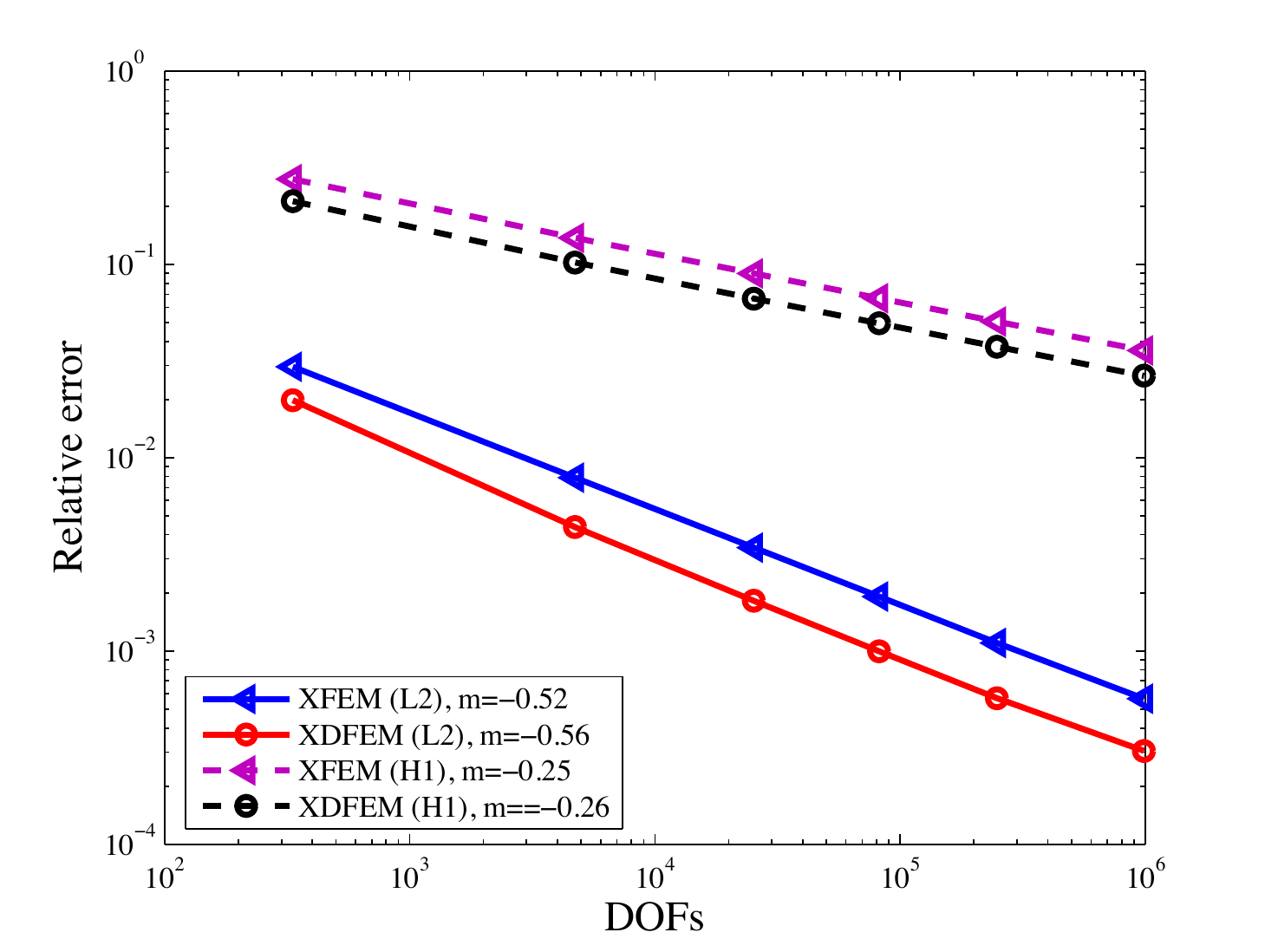}}
  \subfigure[Relative error of SIF]{\includegraphics[width=1.\textwidth]{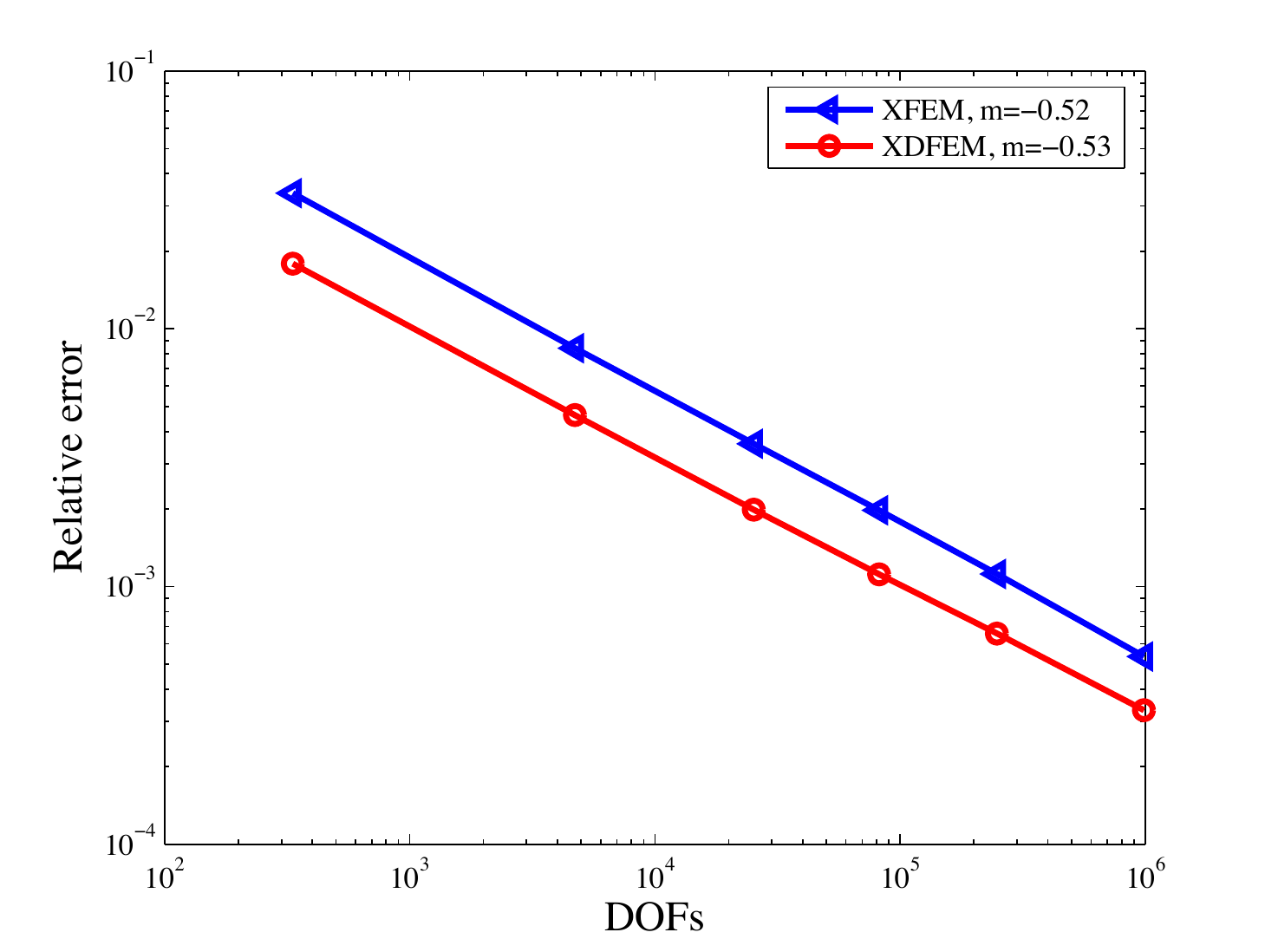}}
  \caption{Convergence results of Griffith (mode I)}
  \label{MODEI}
\end{figure}

\begin{figure}[htbp]
  \centering
  \subfigure[Relative error in displacement and energy norm]{\includegraphics[width=1.\textwidth]{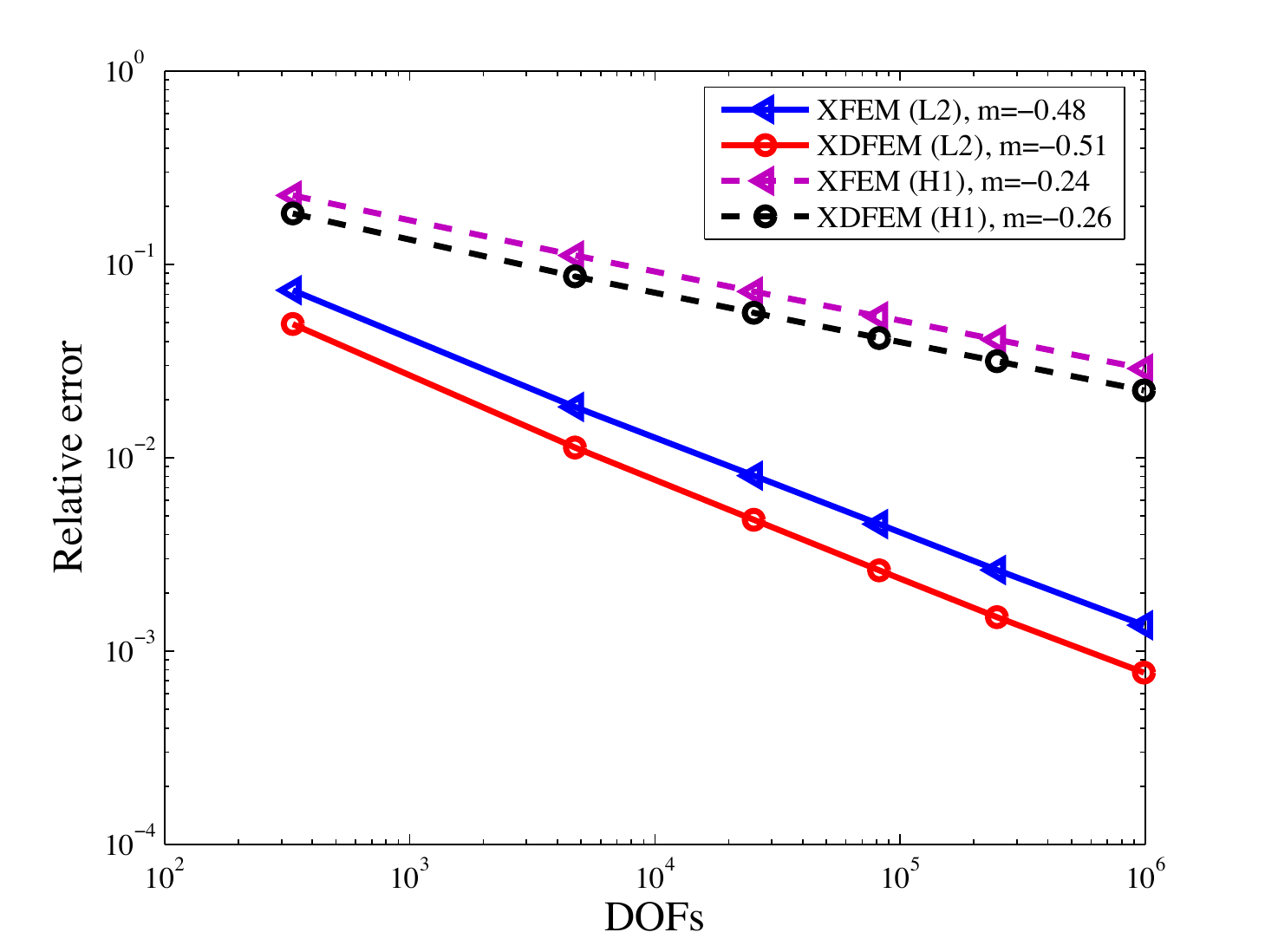}}
  \subfigure[Relative error of SIF]{\includegraphics[width=1.\textwidth]{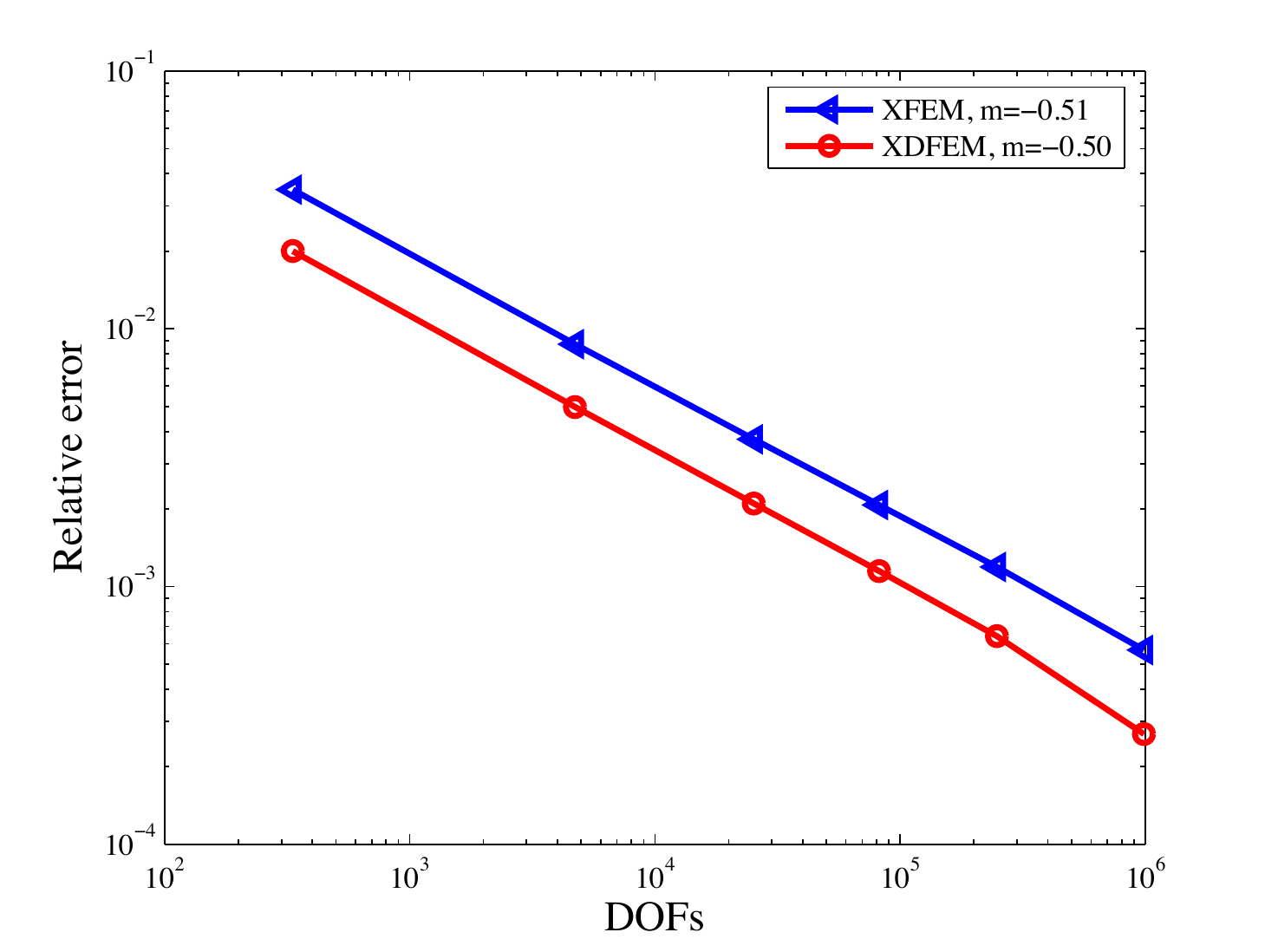}}
  \caption{Convergence results of Griffith crack (mode II)}
  \label{MODEII}
\end{figure}

Approximately 1 million DOFs models of mode I and mode II crack have been simulated to assess the convergence rate (see Fig.\ref{MODEI} and \ref{MODEII}). The relative errors of SIF are also depicted in the plots. From the results it can be inferred that the high order shape functions of XDFEM should be the main reason to achieve high precision and it again confirms the fact that increasing the order of the shape functions can not lead to the optimal convergence rate for XFEM.


\begin{figure}[htbp]
  \centering
  \subfigure[]{\includegraphics[width=1.\textwidth]{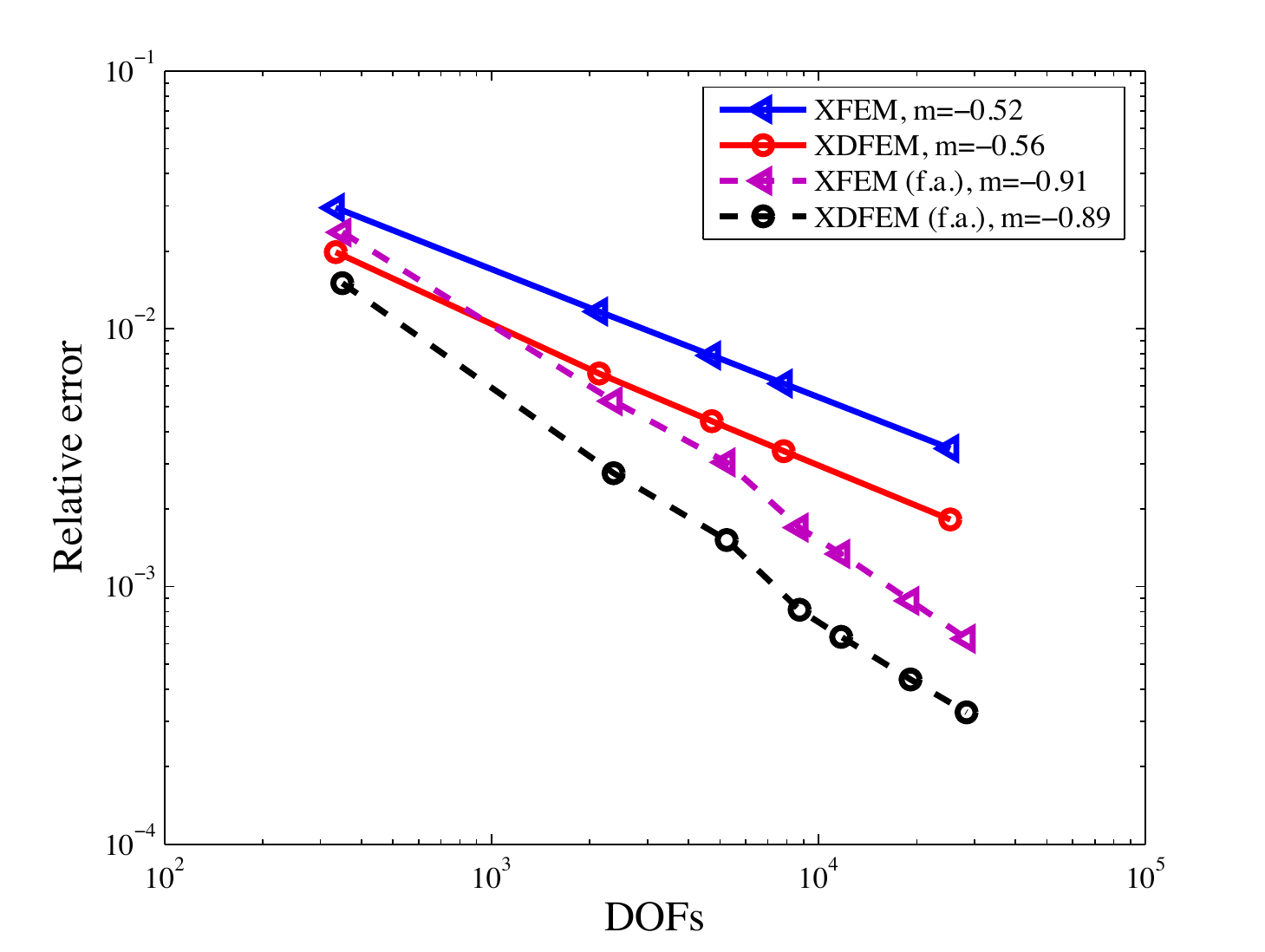}}
  \subfigure[]{\includegraphics[width=1.\textwidth]{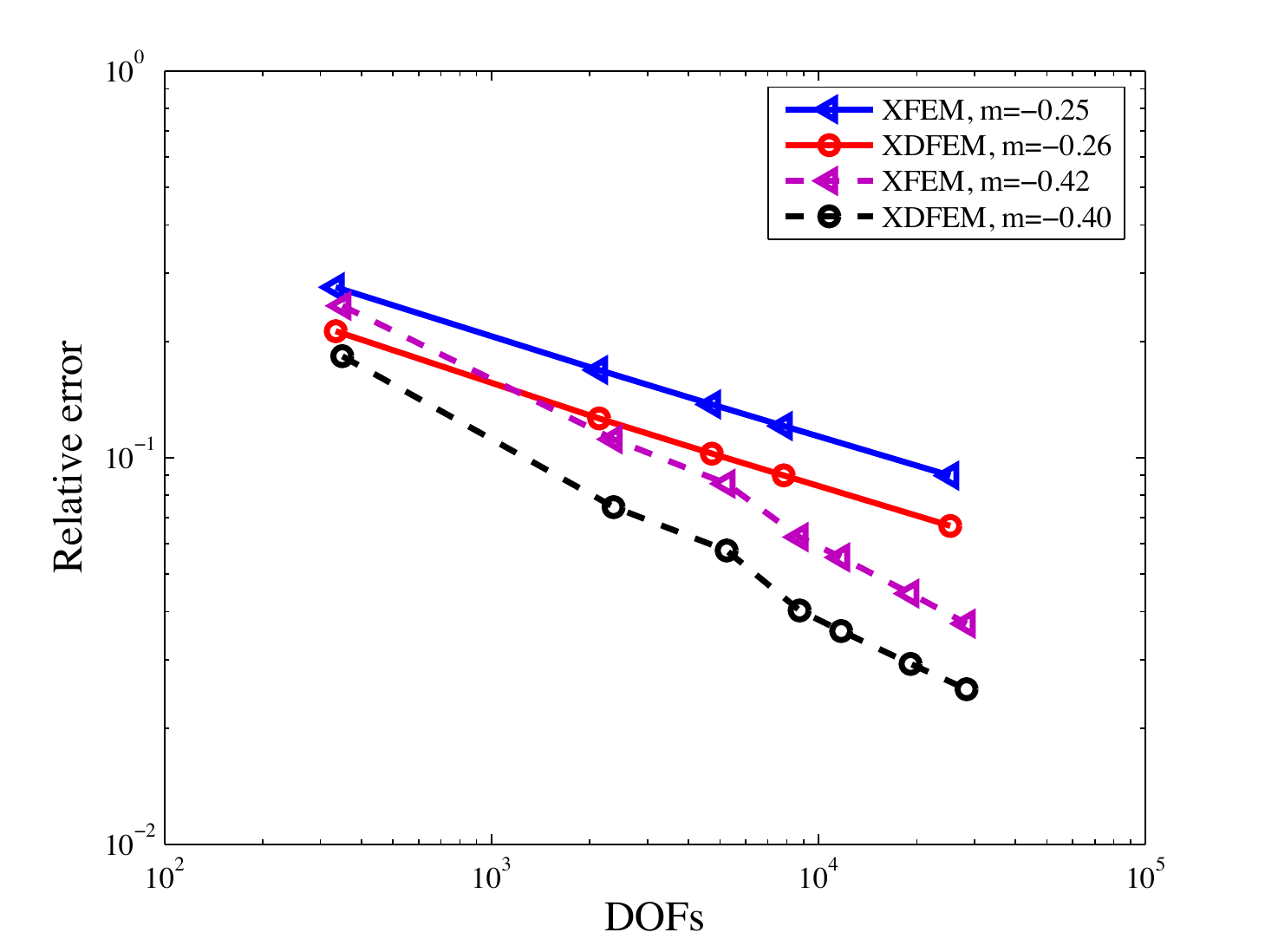}}
  \caption{Convergence results of standard enrichment and fixed area enrichment (mode I): (a) displacement norm error; (b) energy norm error}
  \label{fix}
\end{figure}

\begin{figure}
  {\includegraphics[width=1.\textwidth]{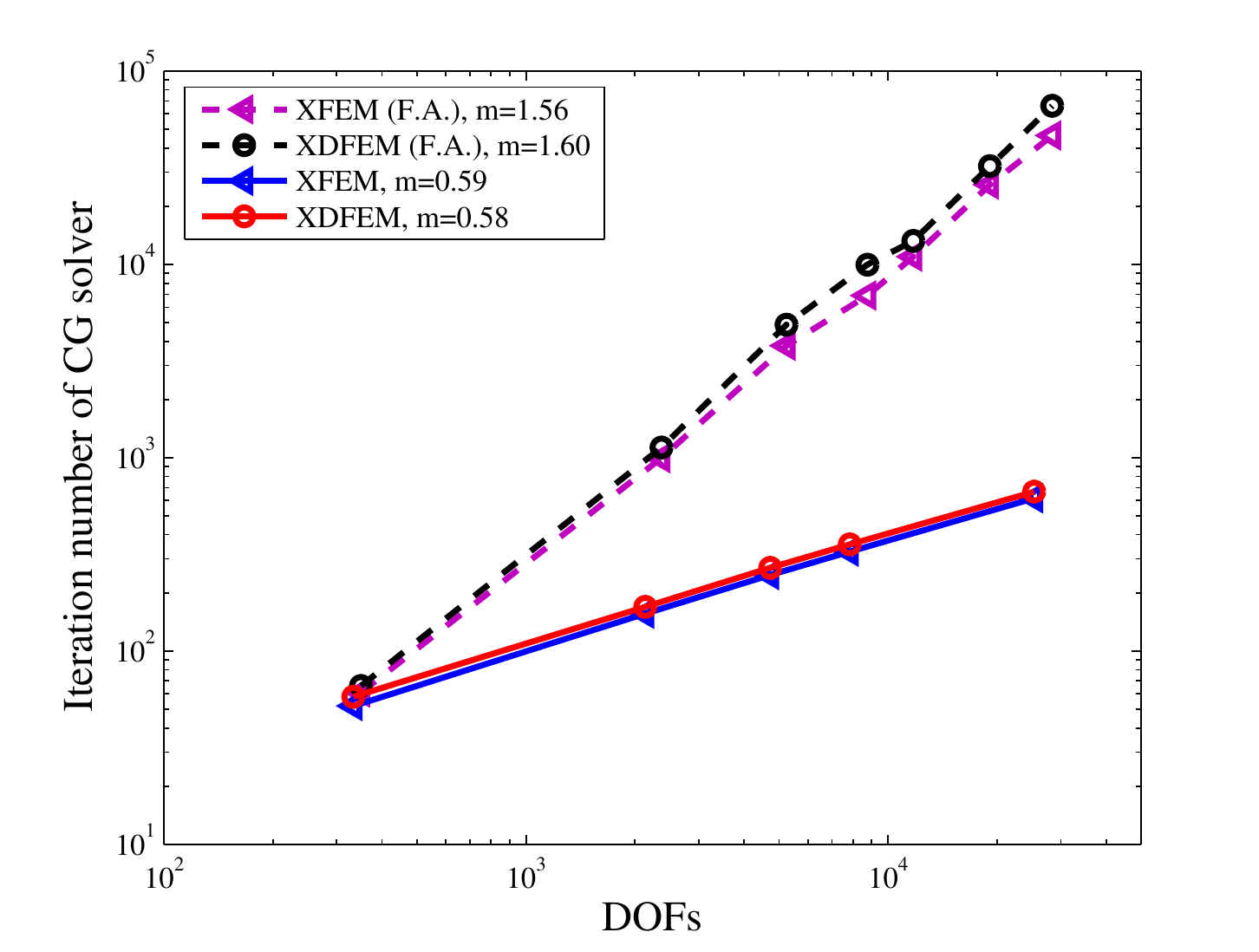}}
  \caption{Iteration number of CG solver}
  \label{numcg}
\end{figure}

The lack of  optimal convergence rate and the aforementioned degeneration of the enriched nodes observed in XDFEM motivated the authors to investigate the effect of  fixed area enrichment in the standard XDFEM. For the purpose of comparison,as in the case of XDFEM, the XFEM with fixed area enrichment was also investigated to evaluate the L2 and H2 norms (Fig. \ref{fix}). The radius parameter $r$ is selected as $1/5$ of the crack length. It can be seen that both two methods achieve the optimal convergence rate with the fixed area enrichment. And XDFEM results are still more accurate than that of XFEM. Fig.\ref{numcg} illustrates the iterative number of the Conjugate Gradient (CG) solver, which can be regarded as an indication to the condition number of the stiffness matrix. It is observed that XDFEM performs slightly worse than XFEM in terms of the condition number. And the deterioration rate of fixed area enrichment is much higher than standard enrichment in both methods. These conclusions are in agreement with the investigation reported by Laborde \emph{et al}.\cite{NME:NME1370}.


\subsubsection{Accuracy study}

\begin{figure}[htbp]
  \centering
  \subfigure[XDFEM]{\includegraphics[width=0.8\textwidth]{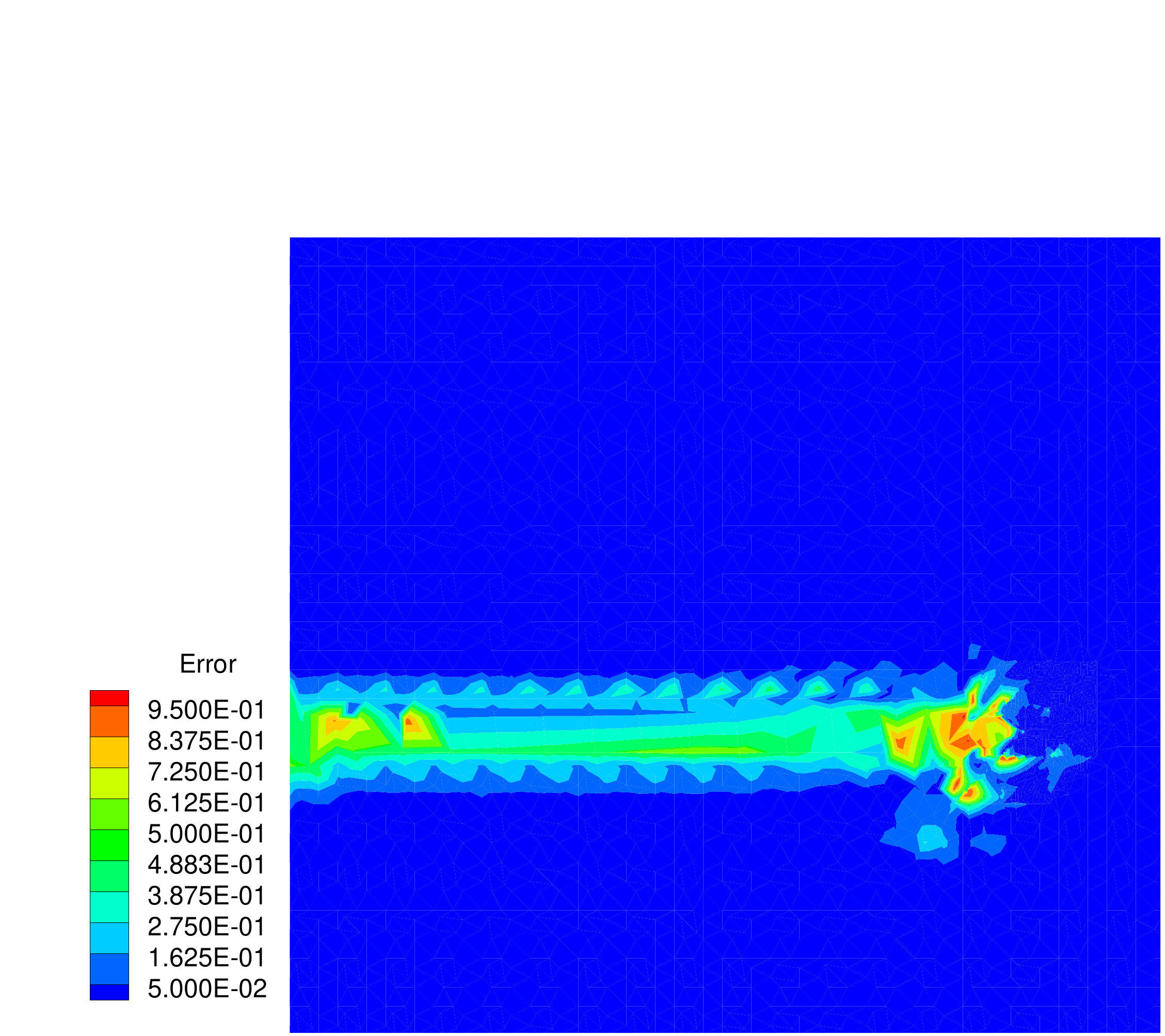}}
  \caption{The contour of relative error in equivalent stress of Mode I crack}
  \label{serror}
\end{figure}

\begin{figure}
  \centering
  {\includegraphics[width=0.8\textwidth]{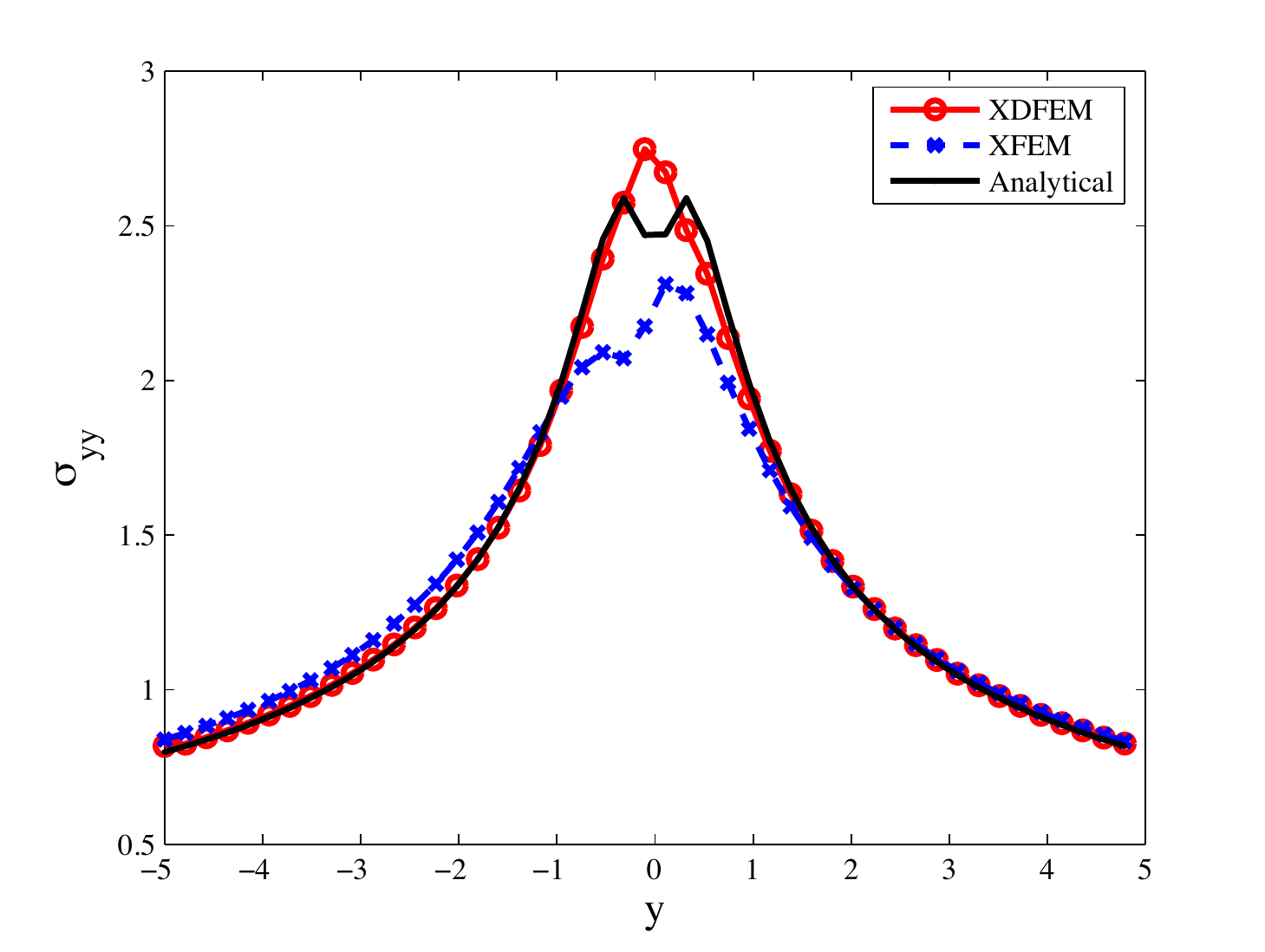}}
  \caption{The comparison of $\sigma_{y}$ at the sample points on the specified line}
  \label{ssline}
\end{figure}

Though it has been already established from the convergence curves that the XDFEM is more accurate than XFEM, it is necessary to investigate whether XDFEM can improve the precision along the crack discontinuity and in the vicinity of the crack tip. Fig.\ref{serror} plots the error contour of the equivalent stress for mode-I crack in a very coarse mesh ($11\times11$). It can be observed clearly that XDFEM improves precision of the stress field significantly near the crack. The stress component $\sigma_{yy}$ are selected at the sample points which are on the line perpendicular to the crack in front of the tip. The results are plotted in Fig.\ref{ssline}. It can be noted from Fig.\ref{ssline}, that the XDFEM results ae much closer to the analytical solution than that of XFEM. Especially, near the crack tip, the XDFEM performs better due to the inclusion of nodal gradients in the approximation.

\begin{figure}[htbp]
  \centering
  \subfigure[XFEM]{\includegraphics[width=0.8\textwidth]{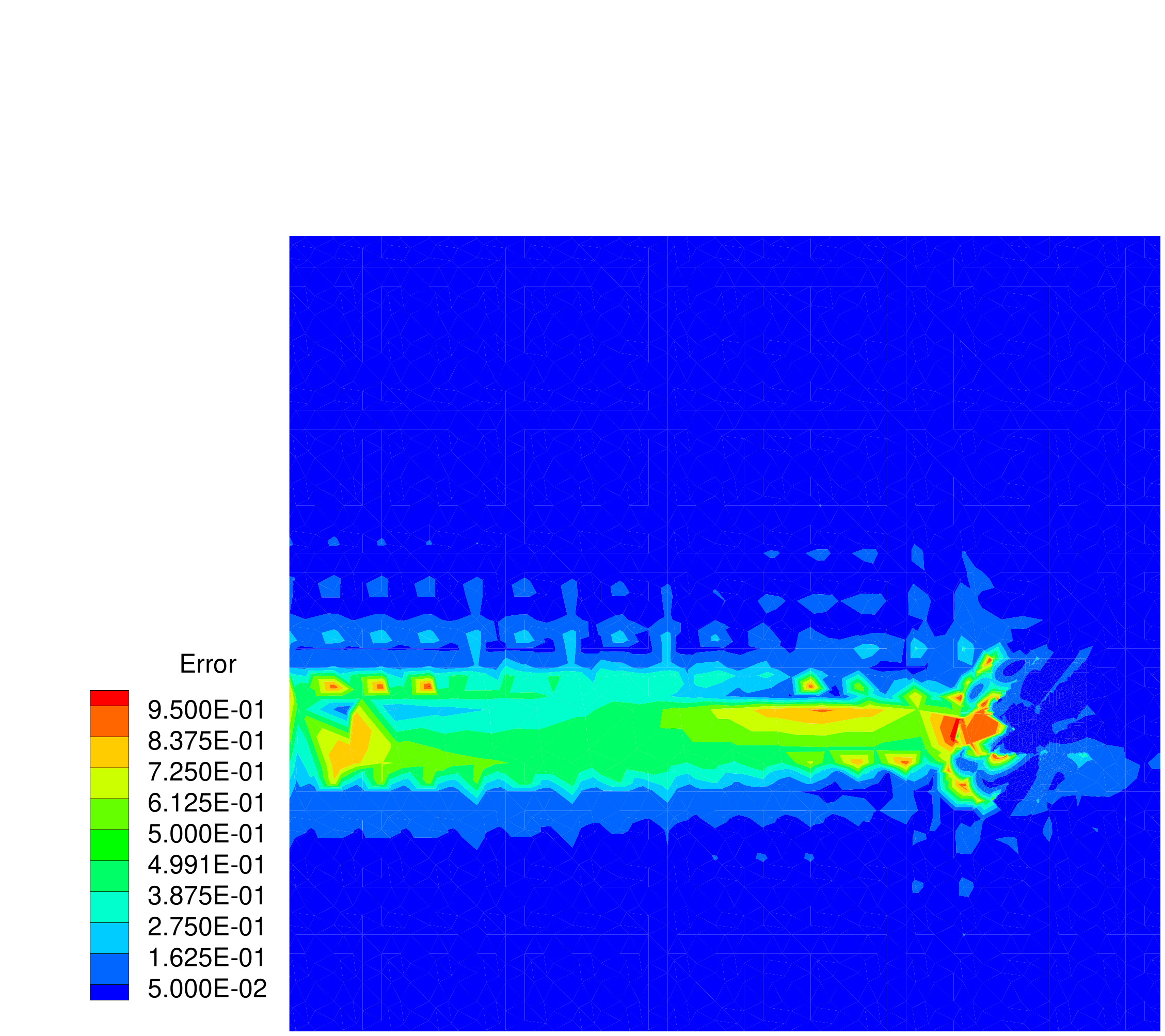}}
\end{figure}

\begin{figure}[htbp]
  \centering
  \subfigure[]{\includegraphics[width=0.49\textwidth]{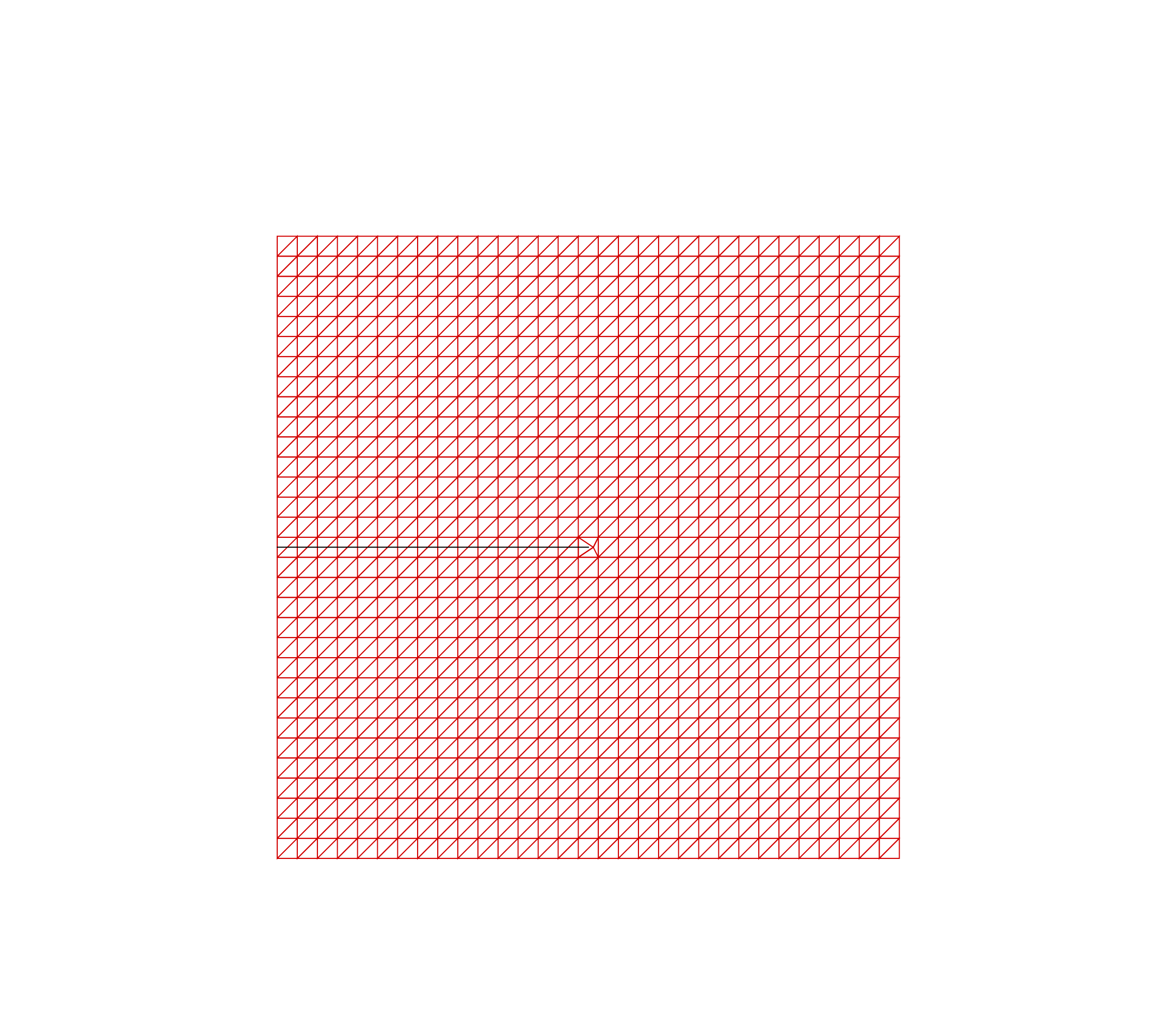}}
  \subfigure[]{\includegraphics[width=0.49\textwidth]{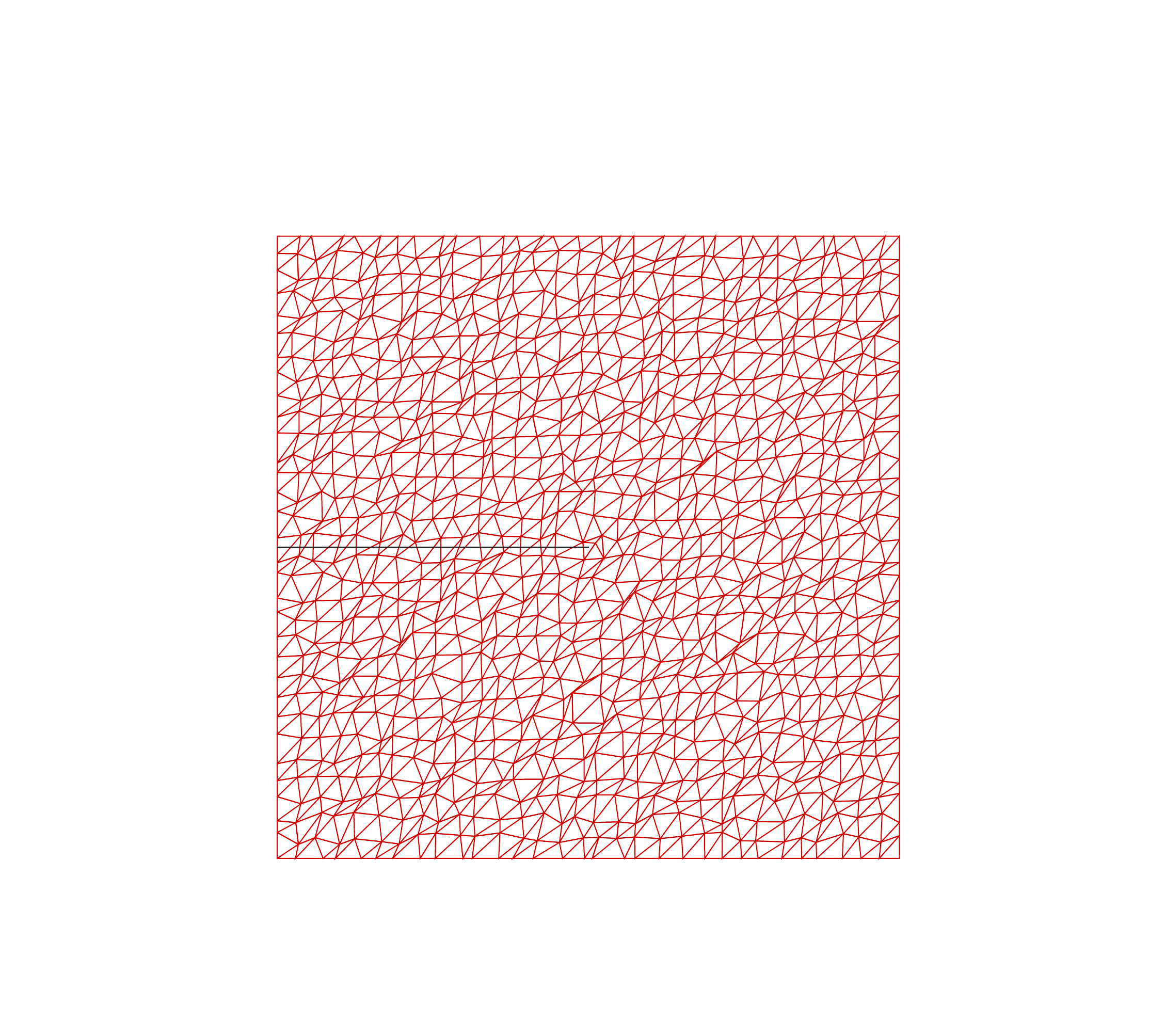}}
  \caption{Mesh design to check the mesh distortion effect:(a)structured mesh;(b) distorted mesh}
  \label{distort}
\end{figure}

The mesh distortion effect is also investigated in this example. The strutured mesh and typical distorted mesh are shown in Fig.\ref{distort}. The results are listed in Table \ref{distorttable}. And the precision of XDFEM in distorted mesh is still better than that of XFEM.
\begin{table}[!hcbp]
\center
\begin{tabular}{|c|c|c|c|c|}
\hline
\multicolumn{1}{|c|}{}& \multicolumn{2}{|c|}{Strutured mesh} & \multicolumn{2}{|c|}{distorted mesh} \\
\hline
DOFs & XFEM & XDFEM & XFEM & XDFEM \\
\hline
$334$ & $0.2272$ & $0.1832$ & $0.2313$ & $0.1882$  \\
\hline
$4726$ & $0.1112$ & $0.08672$ & $0.1132$ & $0.08863$  \\
\hline
$7834$ & $0.09769$ & $0.07600$ & $0.1016$ & $0.08261$ \\
\hline
$17134$ & $0.08006$ & $0.06212$ & $0.08223$ & $0.06215$  \\
\hline
\end{tabular}
\caption{Relative error in energy norm in strutured and distorted mesh}
\label{distorttable}
\end{table}

\subsubsection{Investigation of computational time}
\begin{figure}[htbp]
  \centering
  \subfigure[]{\includegraphics[width=1.\textwidth]{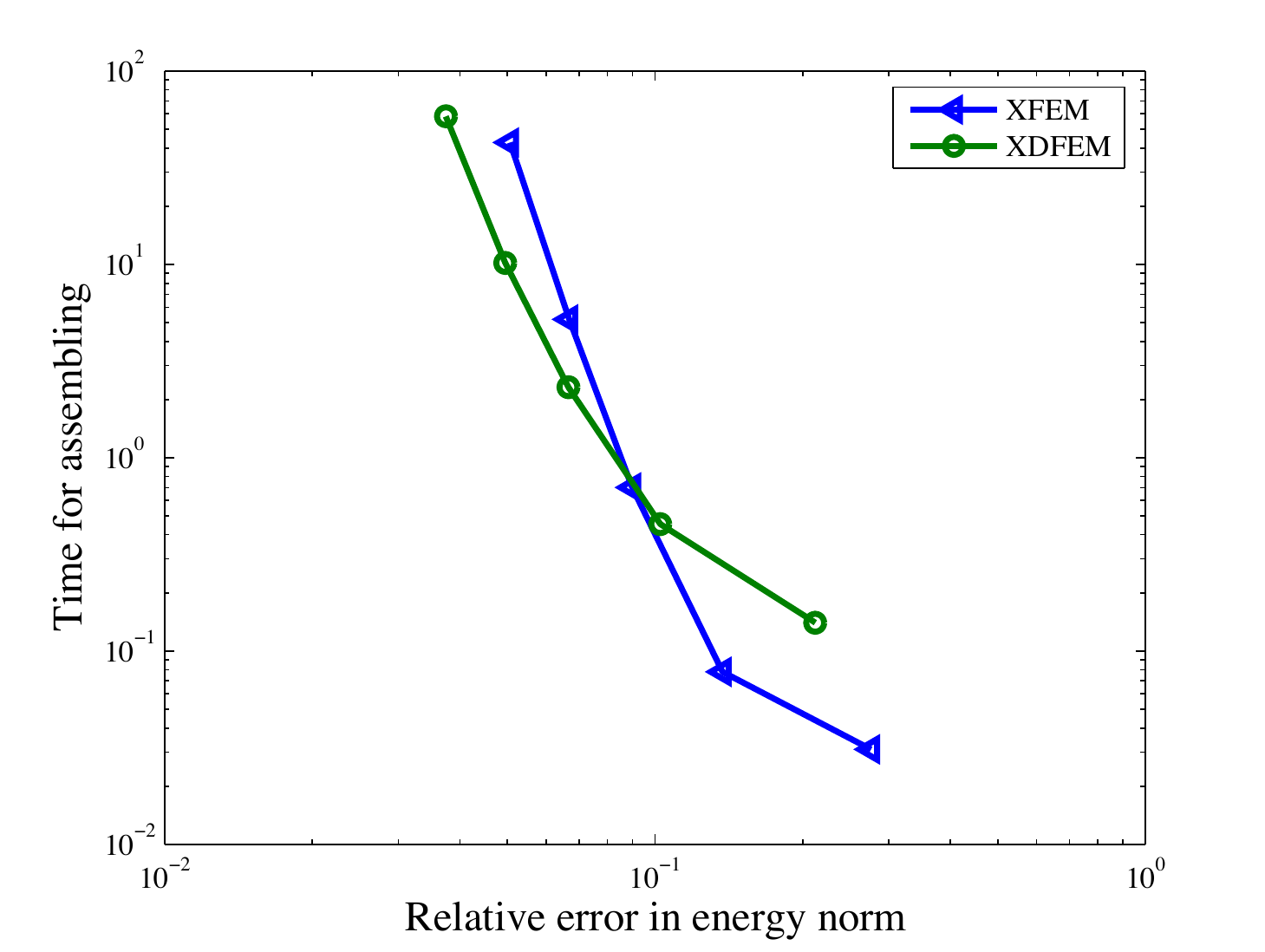}}
  \subfigure[]{\includegraphics[width=1.\textwidth]{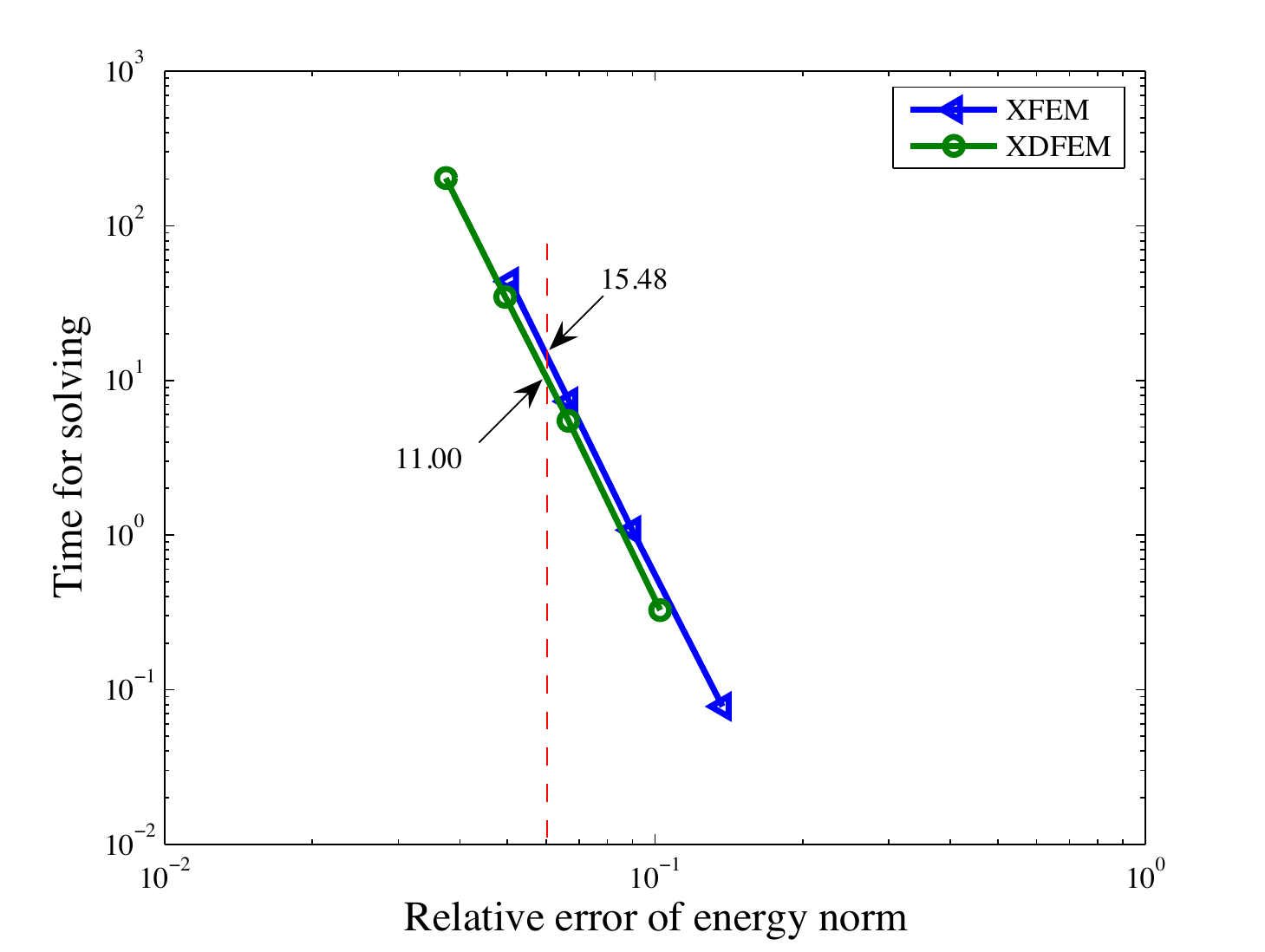}}
  \caption{Comparison of time costs for XFEM and XDFEM in Griffith crack problem }
  \label{times}
\end{figure}

\begin{figure}[htbp]
  \centering
{\includegraphics[width=1.\textwidth]{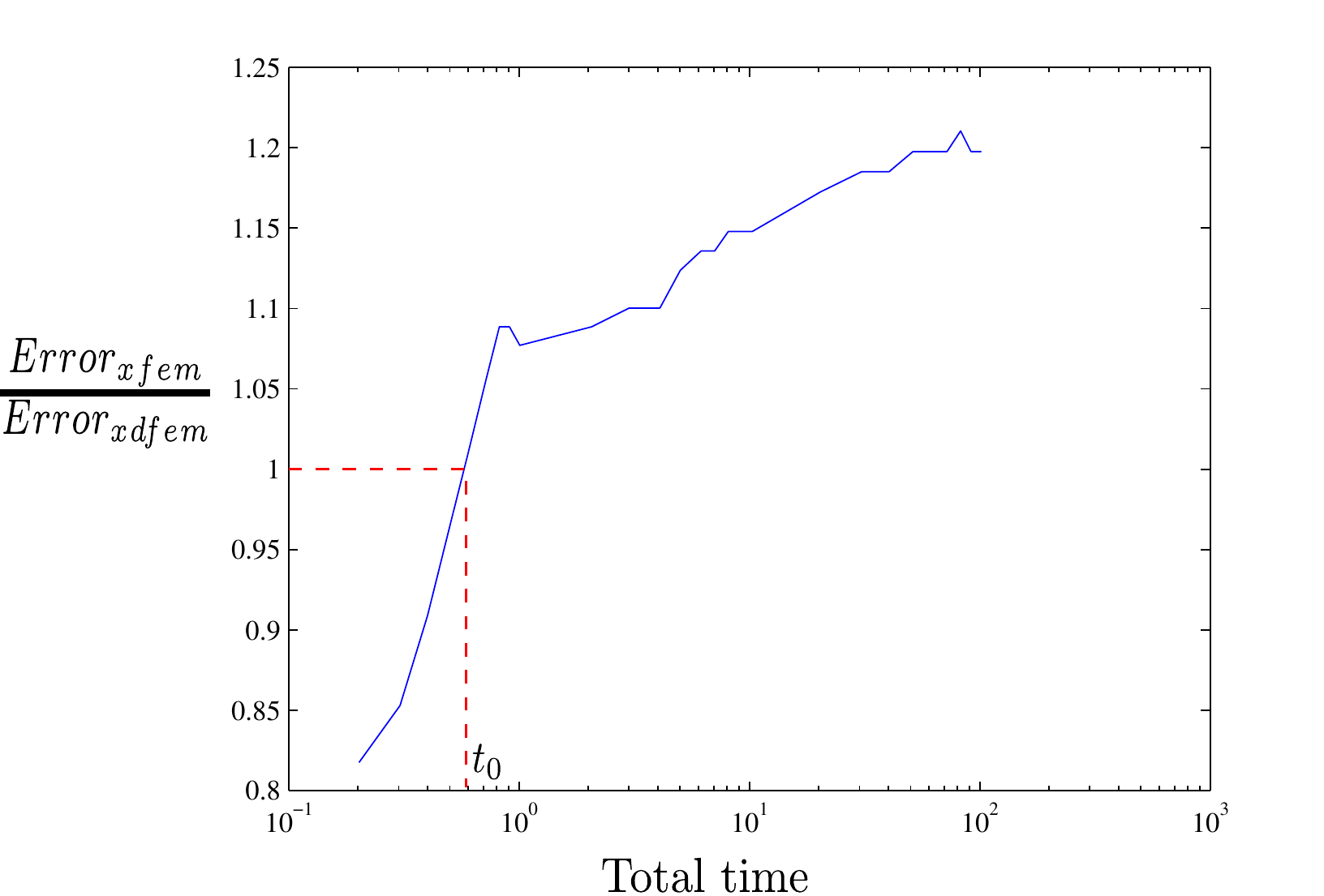}}
  \caption{Total time comparison XFEM and XDFEM in energy norm error in Griffith crack problem }
  \label{time}
\end{figure}

It should be highlighted that the support domain of DFEM element is much bigger than that of FEM element due to the introduction of the nodal gradient into the approximation (see Fig.\ref{domain1D}, Fig.\ref{TFEM}). This directly results in increased bandwidth of the stiffness matrix in DFEM, even though it has the same DOFs as FEM. Hence the solution time is also increased in DFEM or XDFEM, when the same mesh model is adopted. So it is necessary to find a balance between the DOFs, the solution time and precision for both methods. Fig.\ref{times} and \ref{time} show the comparison of the time cost in assembling the stiffness matrix, solving the linear equations and the total time of the two processes (in XFEM, the post-processing of stress smoothing is also included). Here the time for assembling the stiffness matrix implies the time for element stiffness matrix calculation and it does not include the time for searching the support triangles for one node because this process is regarded as a pre-processing procedure. It can be seen that with the model size increasing, XDFEM requires less time to obtain the same precision. For the solving process, XFEM produces error 1.4 times higher ($\frac{XFEM 15.48}{XDFEM 11.}=1.4$) than the XDFEM at a same solving time. The total time comparison shows that after $t_0=0.6$, the XFEM tends to produce higher and higher error than XDFEM. The main factor for this process is the bandwidth. It should be mentioned that the program is solved by CG method, with the storage in the strategy of varied half bandwidth method. Therefore it can be inferred that if the element-by-element approach is used when combining the CG solver, the XDFEM should be much more efficient as the bandwidth effect "disappears".


\subsection{Inclined center crack}
The inclined center crack is investigated in this section. The model is presented in Fig.\ref{griff}(b). The infinite plate is subjected to remote tensile load in $y$ direction and the inclined angle $\beta$ is measured in counter-clockwised direction from the $x$ direction. The half crack length $a=1000$. A square domain area ($10\times 10$) encircling the crack tip is selected and the exact displacement is applied as in the previous example. The analytical SIFs are given as
\begin{subequations}
\begin{equation}
K_I=\sigma \sqrt{\pi a}\text{cos}^2\beta
\end{equation}
\begin{equation}
K_{II}=\sigma \sqrt{\pi a}\text{cos}\beta\text{sin}\beta
\end{equation}
\end{subequations}
Table \ref{siftable} shows that the $K_I$ and $K_{II}$ vary with the inclined angle in inclined crack. It can be observed from Table \ref{siftable} that both XDFEM and XFEM results agree well with the analytical solution. The precision of SIFs in XDFEM are better than that of XFEM. This example demonstrates that XDFEM performs well in the mixed mode crack problem.

\begin{table}[!hcbp]
\center
\begin{tabular}{|c|c|c|c|c|c|c|}
\hline
\multicolumn{1}{|c|}{}& \multicolumn{3}{|c|}{$K_I$} & \multicolumn{3}{|c|}{$K_{II}$} \\
\hline
$\beta$  & exact & XFEM & XDFEM & exact & XFEM & XDFEM \\
\hline
$0$ & $1.0000$ & $1.0058$ & $1.0030$ & $0.0000$ & $-0.0010$ & $-0.0003$ \\
\hline
$\frac{\pi}{12}$ & $0.9330$ & $0.9381$ & $0.9357$ & $0.2500$ & $0.2502$ & $0.2503$ \\
\hline
$\frac{\pi}{6}$ & $0.7500$ & $0.7537$ & $0.7520$ & $0.4330$ & $0.4343$ & $0.4339$ \\
\hline
$\frac{\pi}{4}$ & $0.5000$ & $0.5022$ & $0.5012$ & $0.5000$ & $0.5018$ & $0.5011$ \\
\hline
$\frac{\pi}{3}$ & $0.2500$ & $0.2508$ & $0.2505$ & $0.4330$ & $0.4348$ & $0.4340$ \\
\hline
$\frac{5\pi}{12}$ & $0.0670$ & $0.0670$ & $0.0671$ & $0.2500$ & $0.2511$ & $0.2506$ \\
\hline
$\frac{\pi}{2}$ & $0.0000$ & $0.0000$ & $0.0000$ & $0.0000$ & $0.0000$ & $0.0000$ \\
\hline
\end{tabular}
\caption{Normalized SIFs in inclined center crack ($47\times47$)}
\label{siftable}
\end{table}

\subsection{XDFEM for crack propagation}

\begin{figure}
  {\includegraphics[width=1.\textwidth]{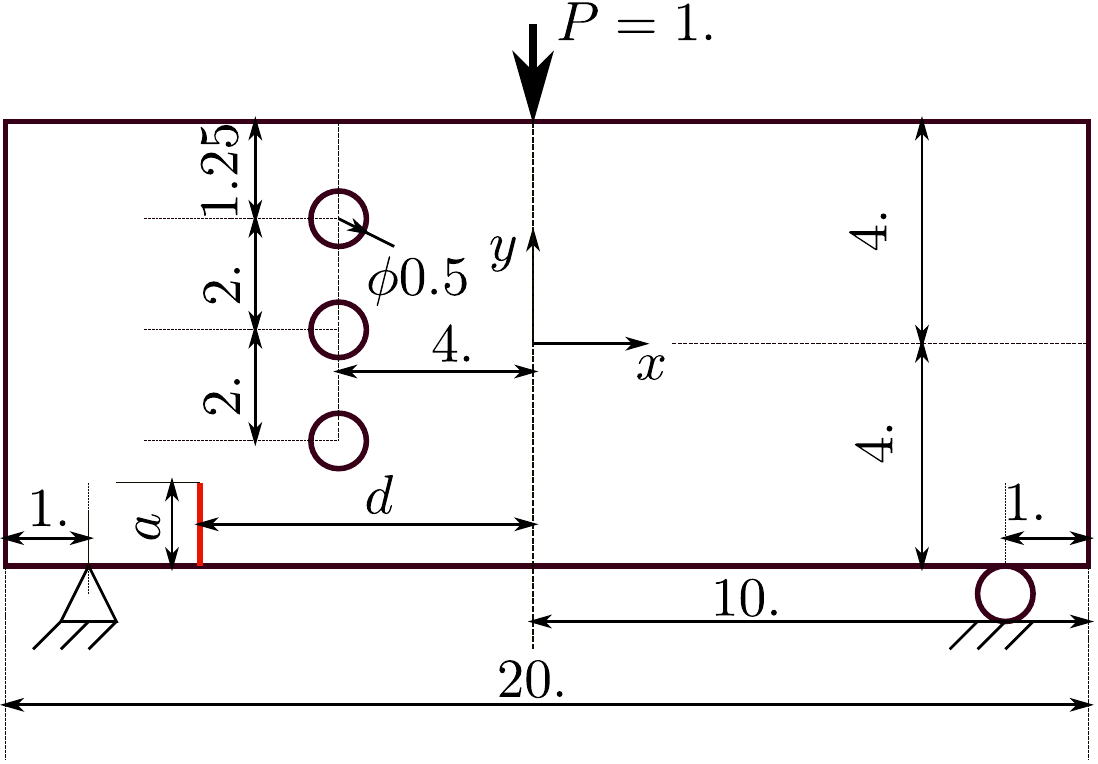}}
  \caption{Physical model of 3 points bending beam with 3 holes}
  \label{3hole}
\end{figure}

\begin{figure}[htbp]
\centering
  \subfigure[XFEM, the 94th step]{\includegraphics[width=0.8\textwidth]{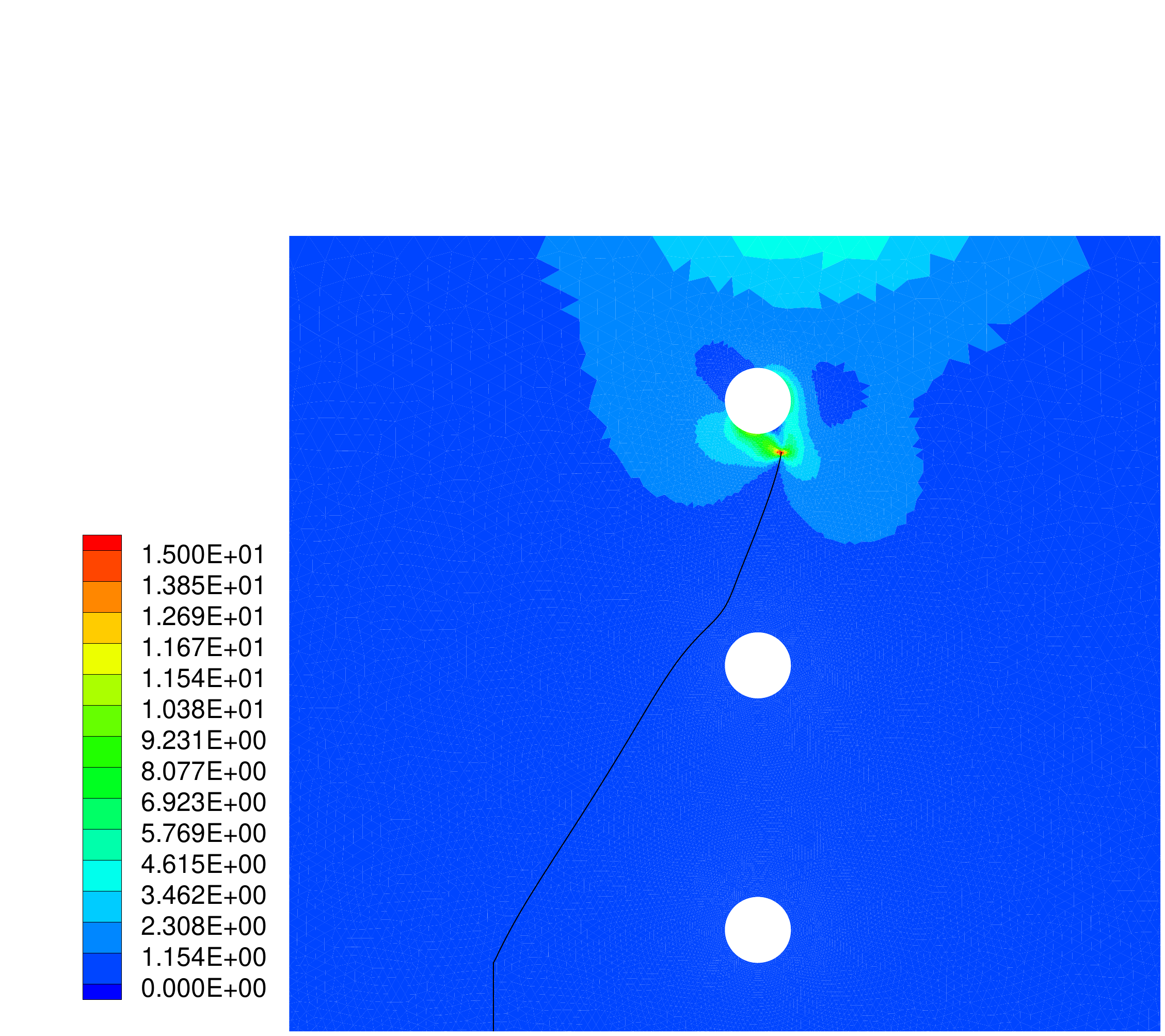}}
  \subfigure[XDFEM, the 94th step]{\includegraphics[width=0.8\textwidth]{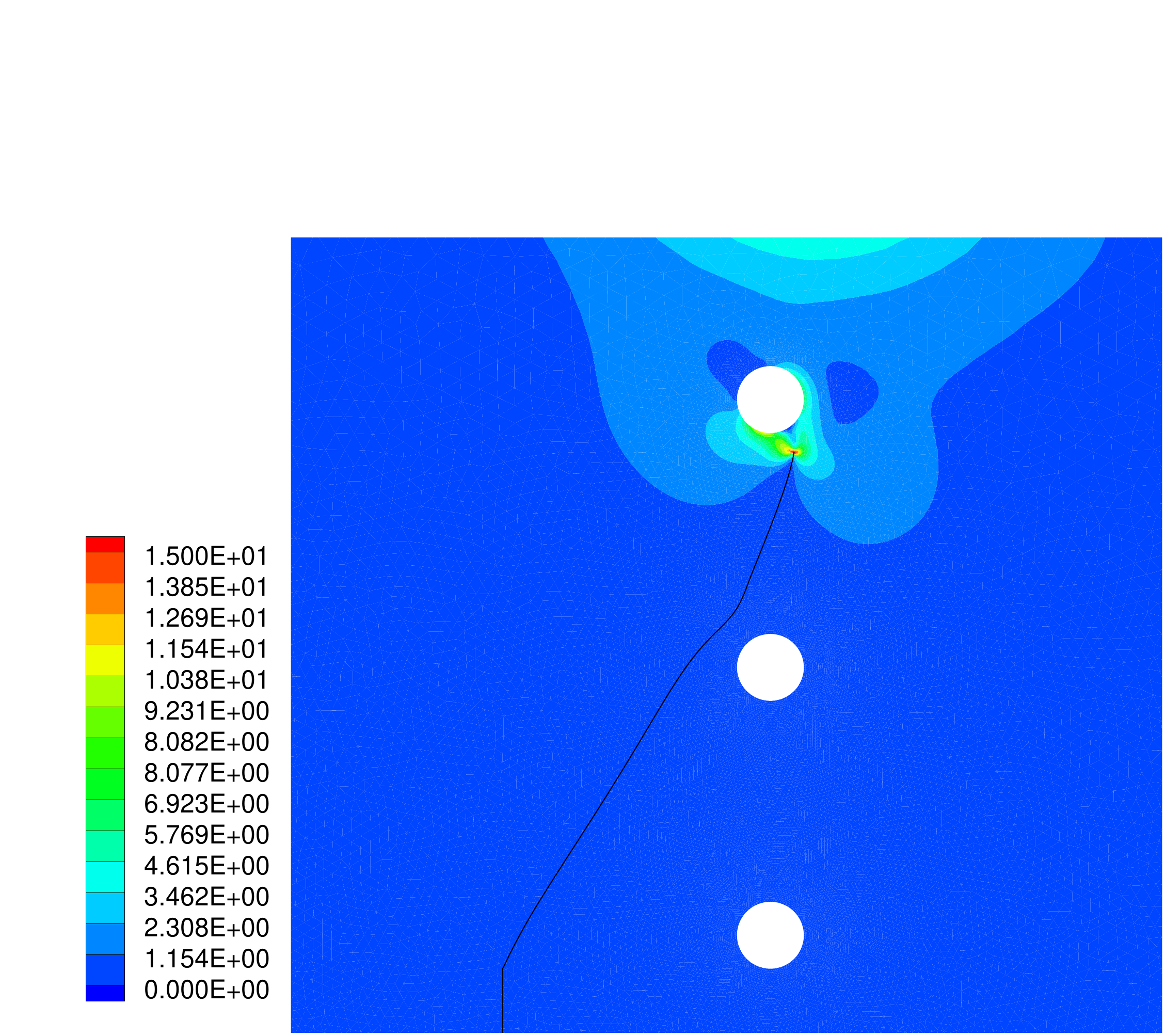}}
  \caption{Contour plots of equivalent stress in case 3}
  \label{s13hole}
\end{figure}

\begin{figure}[htbp]
\centering
  \subfigure[]{\includegraphics[width=0.8\textwidth]{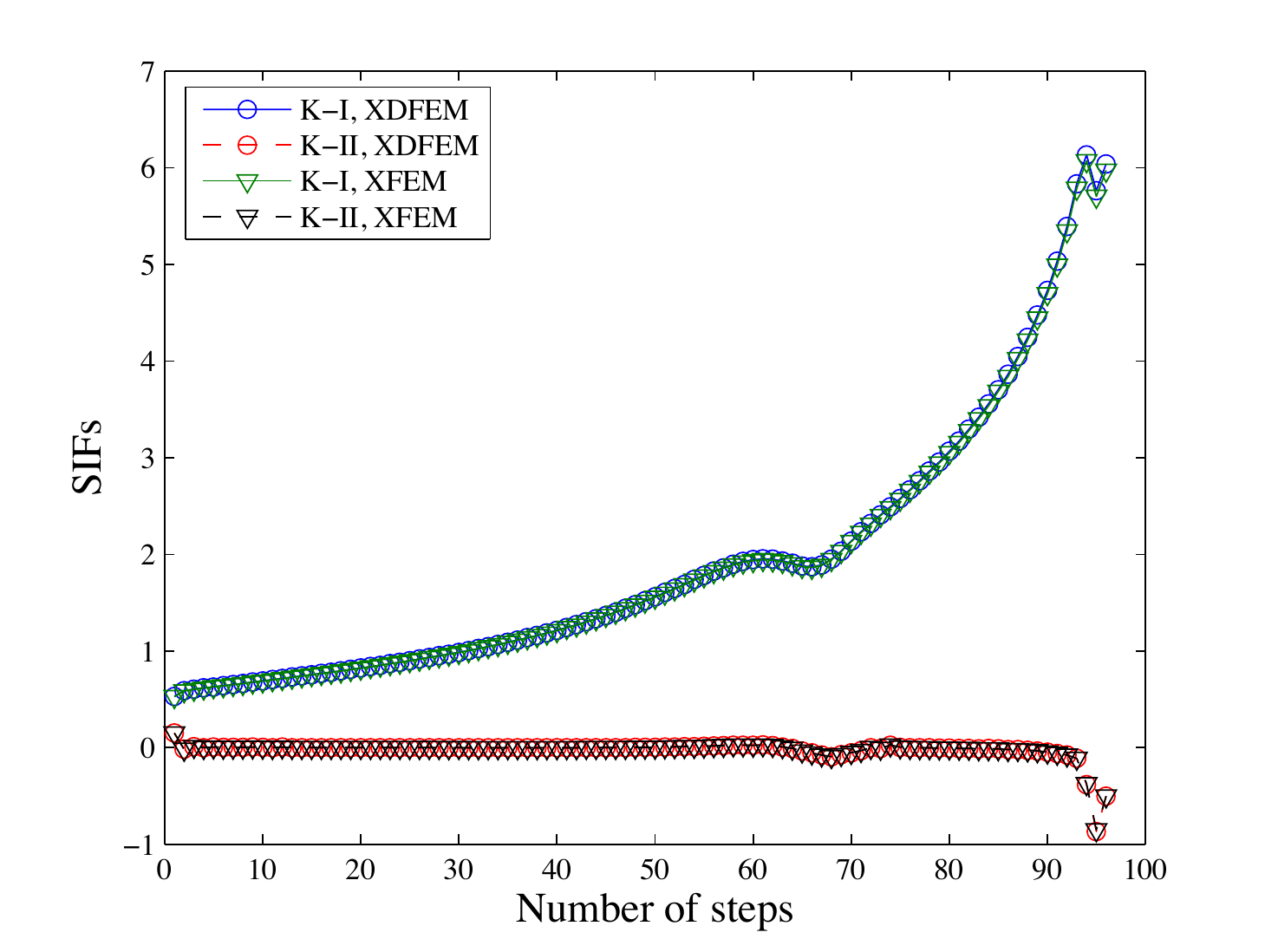}}
  \caption{SIFs variation in three cases}
  \label{sif3hole}
\end{figure}

\begin{figure}
\centering
{\includegraphics[width=1.\textwidth]{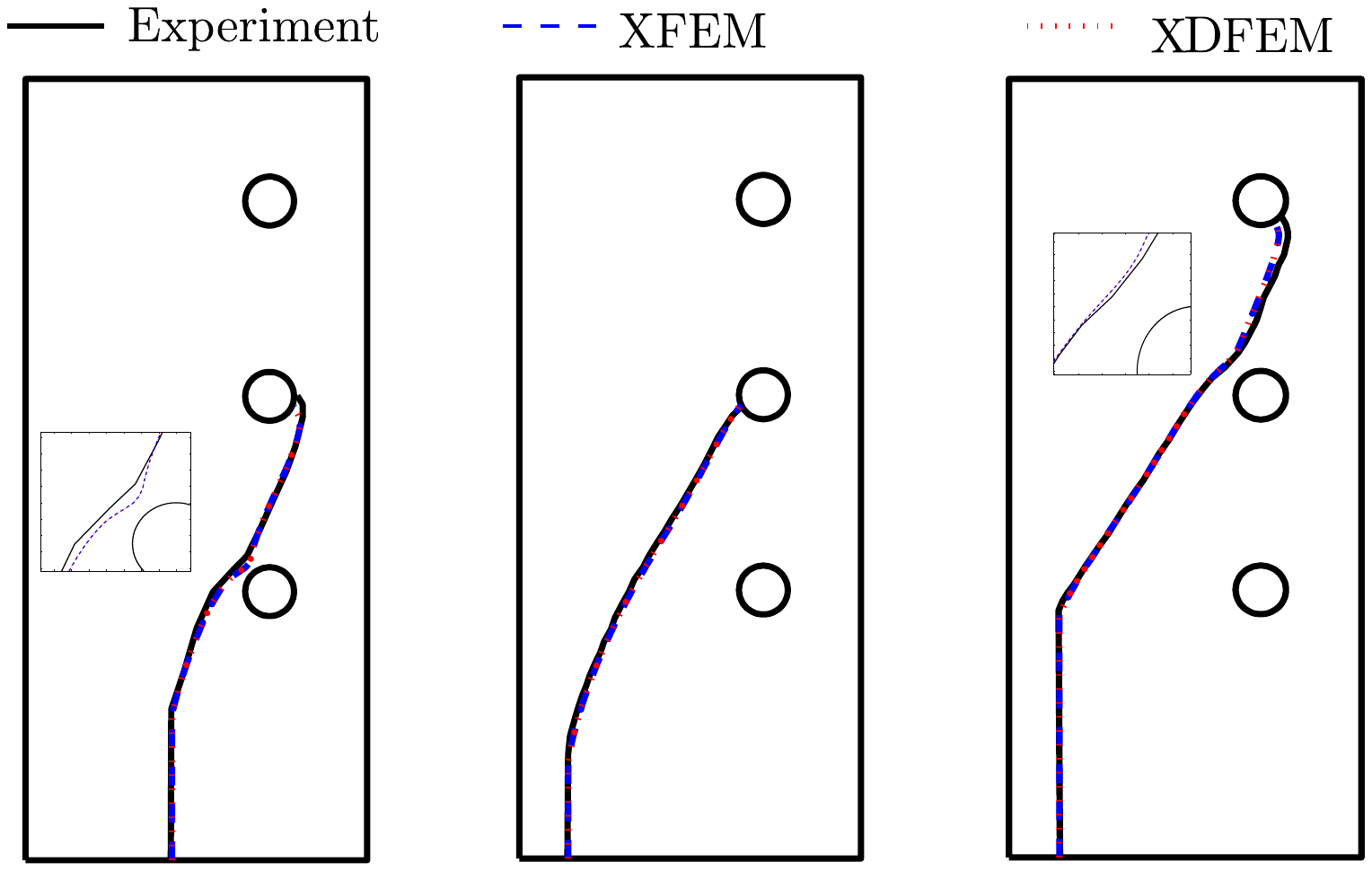}}
  \caption{Crack evolution of the three cases}
  \label{case3hole}
\end{figure}

The example of three point bending beam with three holes is simulated in this section to test the versatility of XDFEM in the application of crack propagation. Holes tend to be an important source of uncertainty to the crack growth in mechanical structures. This experiment is designed to explore the effect of holes on the crack trajectories. The geometry and load condition are illustrated in Fig.\ref{3hole}. Plaxiglas specimens are used as the material and $\nu=0.37$ is used in the simulations. Plane strain condition is assumed. With the variation of the position of the initial crack, different crack trajectories are obtained\cite{Ingraffea1990} \cite{Bittencourt1996}. A set of test cases as listed in Table \ref{holescase} are simulated. The maximum hoop stress criterion \cite{Moes1999} is adopted to judge the orientation of crack propagation. The model is discretized by $27869$ nodes and $55604$ triangular elements. Fig.\ref{s13hole} compares the stress contours of XFEM and XDFEM. Apparently, XDFEM can obtain the smoothed stress field without any post-processing. Fig.\ref{case3hole} illustrates the crack evolution of the listed three cases. And the results show that both methods are in good agreement with the experiment. The SIFs for the three crack trajectories are plotted in Fig.\ref{sif3hole}. It can be observed that the value of each case compares well. In the numerical tests it is noted that, though, the error of energy norm decreases in XDFEM, the crack path is still very consistent with XFEM. It can also be observed from Fig.\ref{case3hole} that, there is very minor difference in the crack path trajectory between XFEM and XDFEM. This observation can be attributed to the fact that, with small increment and very fine mesh, the XFEM can readily  produce accurate crack propagation. However, the accuracy of the stress field associated with this solution may not be significantly high. Similar results are presented by smoothed XFEM \cite{Chen2012}. This example proves the robustness and stability of the XDFEM in simulating the crack propagation.

\begin{table}[!hcbp]
\center
\begin{tabular}{|c|c|c|c|c|}
\hline
      & $d$   & $a$   & crack increment & number of propagation \\
\hline
case 1& $5$   & $1.5$ & $0.052$ & $67$\\
\hline
case 2& $6$   & $1.0$  & $0.060$ & $69$\\
\hline
case 3& $6$   & $2.5$ & $0.048$ & $97$\\
\hline
\end{tabular}
\caption{Test cases for the three points bending problem}
\label{holescase}
\end{table}

\begin{figure}[htbp]
\centering
  \subfigure[XFEM, the 69th step]{\includegraphics[width=0.8\textwidth]{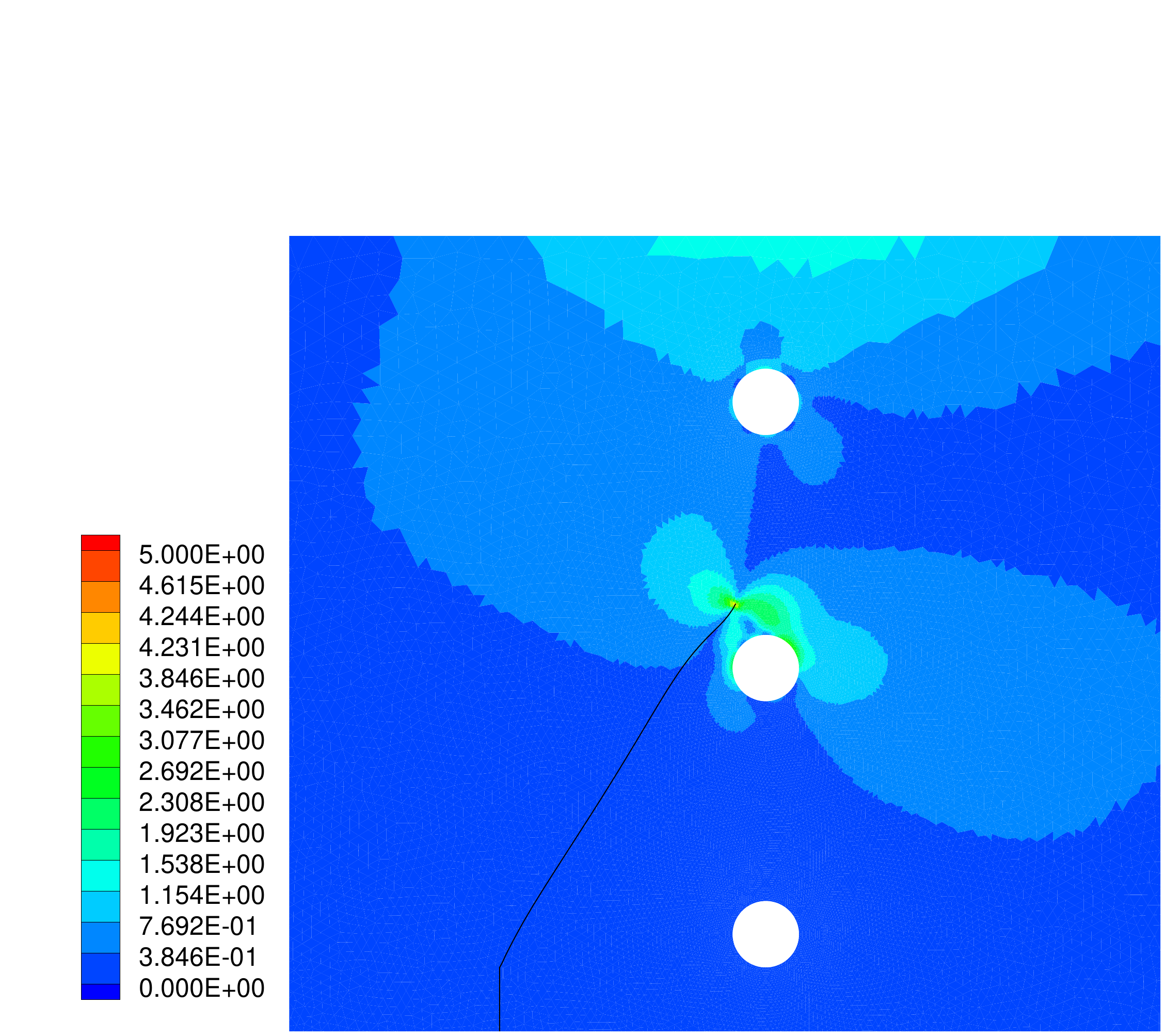}}
  \subfigure[XDFEM, the 69th step]{\includegraphics[width=0.8\textwidth]{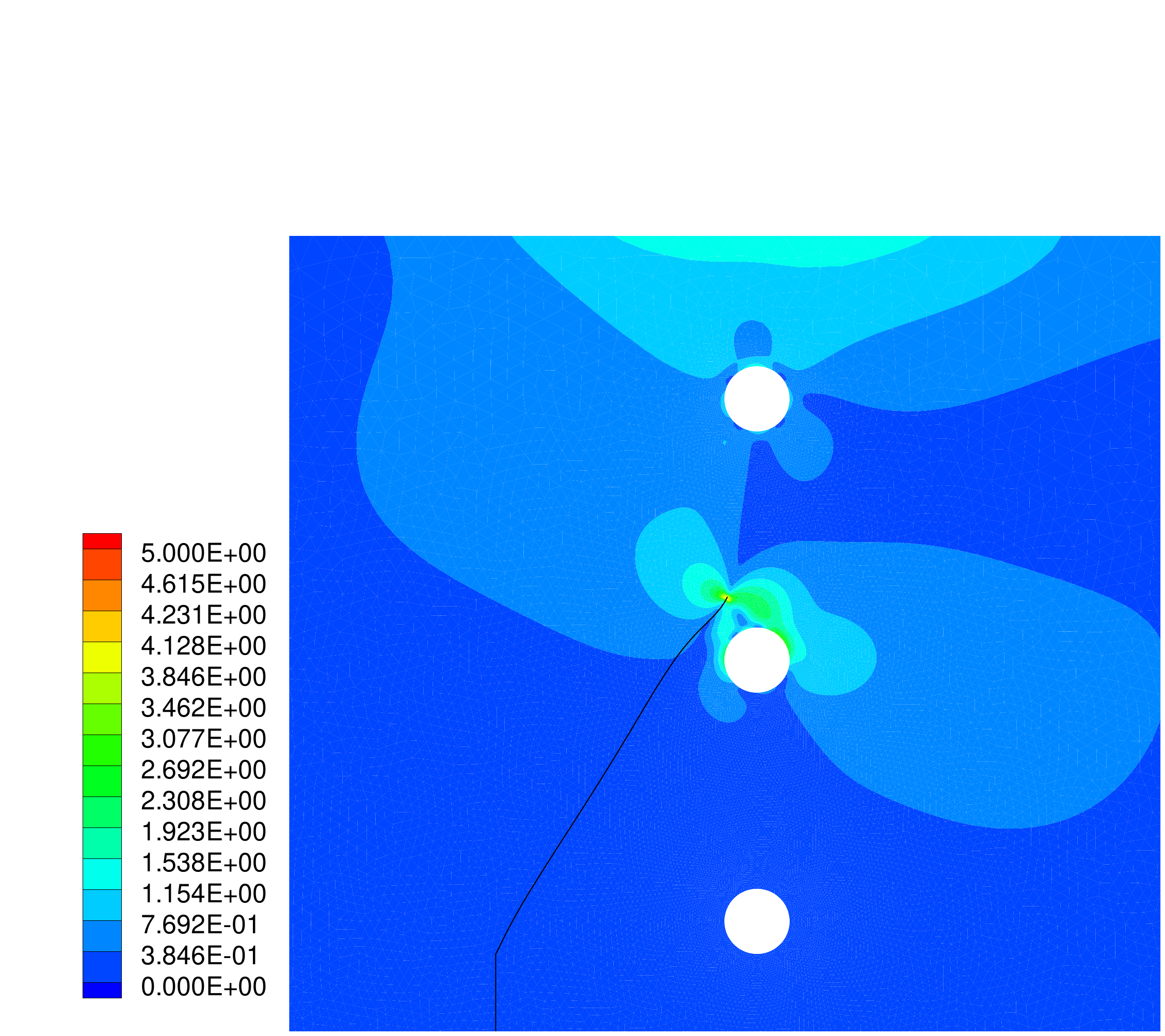}}
\end{figure}
\begin{figure}[htbp]
\centering
  \subfigure[]{\includegraphics[width=0.8\textwidth]{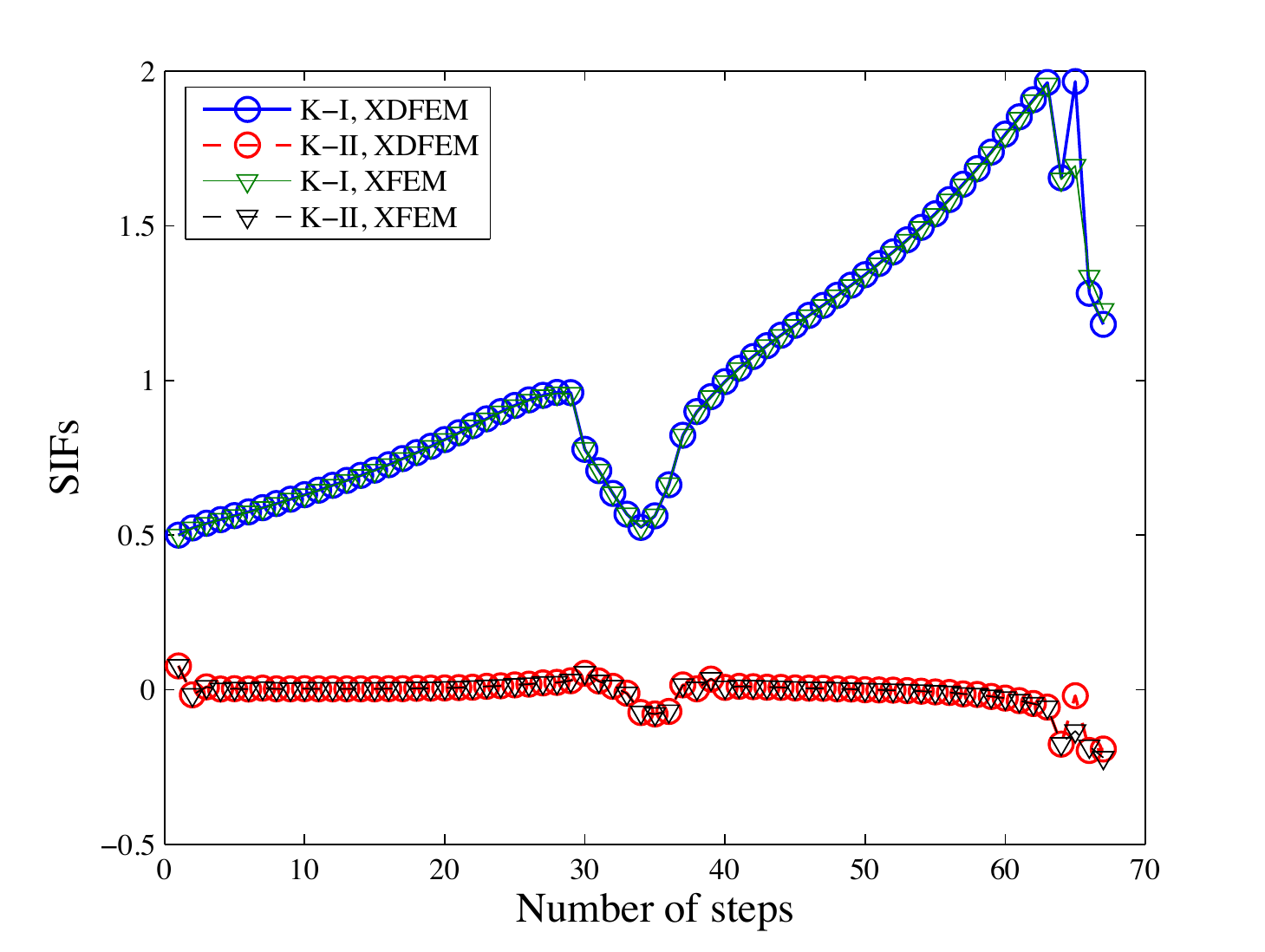}}
  \subfigure[]{\includegraphics[width=0.8\textwidth]{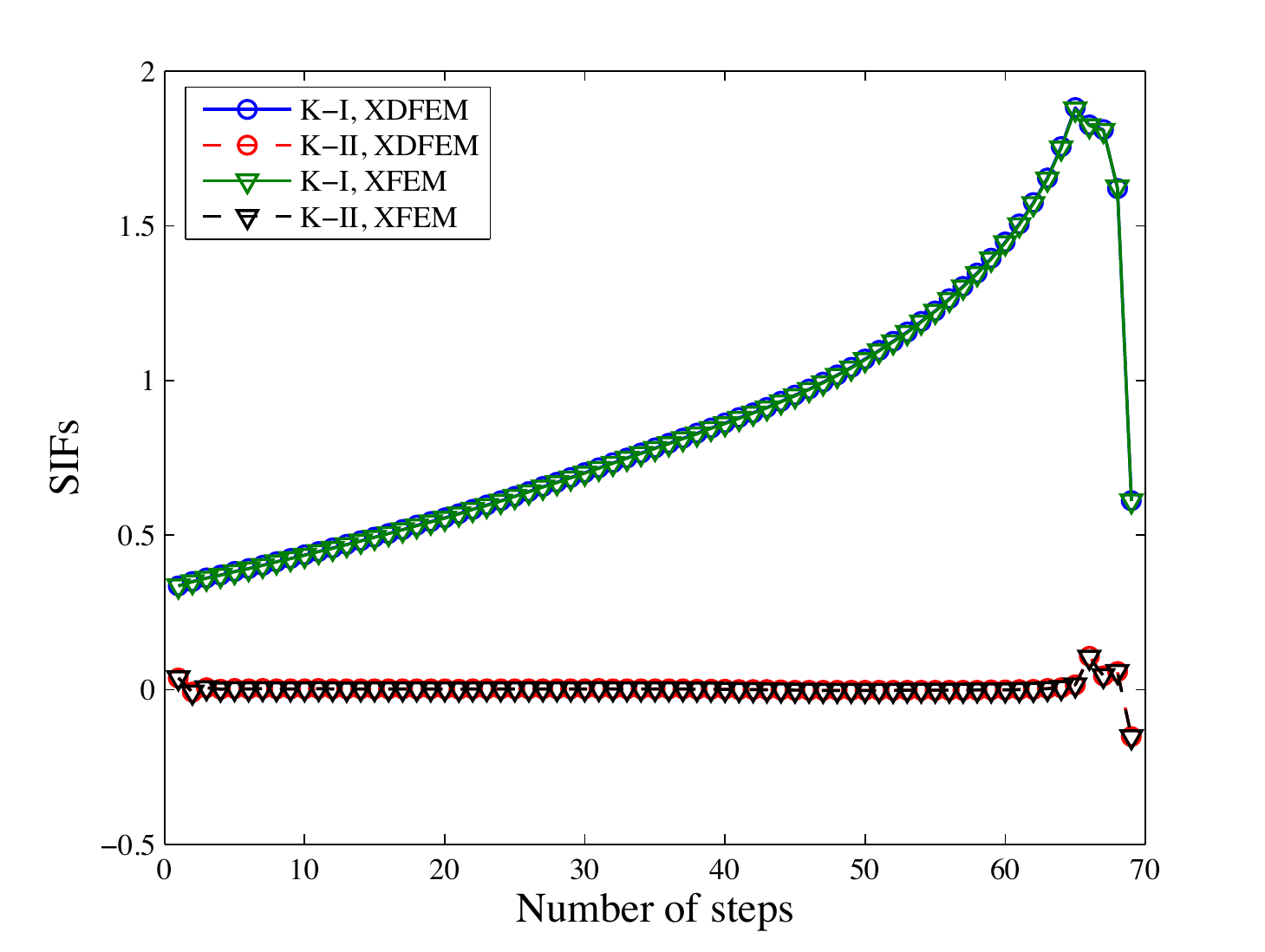}}
\end{figure}

\section{Conclusion}
A novel double-interpolation approximation based on the simplex mesh discretization is introduced and used in its natural and enriched form to simulate the problem of linear elastic fracture mechanics. Several examples are investigated to demonstrate the basic features of DFEM and XDFEM. The conclusion are listed as follows:

\begin{itemize}
\item Due to the nodal gradient aspect of the approximation, the DFEM and XDFEM can improve the precision of simulation result, without the inclusion of additional DOFs. Numerical tests have shown that the double-interpolation method is even more accurate than the Q4 element in the same model size, despite of using the simplex mesh dicretization. It is well known that the quadrilateral (hexahadronal) mesh can achieve much better accuracy while, simplex mesh is more convenient to generate and has a strong adaptivity for using with arbitrary shape. DFEM proposed here proves to unite the two factors together to provide very accurate results using the simplex mesh.
\item In contrast to common high order finite element, DFEM has $C^1$ continuity on most nodes. For continuum mechanics problem, it does not require any post-processing for recovering the nodal stress. For XDFEM, only the enriched nodes need extra post-processing. Post-processing procedure is almost omitted in DFEM and XDFEM, which improves the efficiency of the simulation and maintains all the variables in the valid space.
\item In terms of computational cost, DFEM has proved to be more efficient than XFEM in the computational efficiency study. It should be highlighted that the major factor which hampers the efficiency is the increased bandwidth issue of the double-interpolation method. When the element-by-element strategy is used, the bandwidth should not be a problem and XDFEM is expected to be much faster than XFEM when combined with the preconditioned CG (PCG) solver.
\item The application of XDFEM in crack propagation is as robust as XFEM, with a smoother stress field generated for evaluation. The investigation in the local error near the crack tip shows that XDFEM can obtain more accurate stress field. This could be more helpful in improving the accuracy in 3D XFEM, in which the precision is very low due to the low continuity of the solution.
\end{itemize}

In conclusion, XDFEM proves to be a very promising alternative to XFEM. For the future work, the 3D XDFEM should be investigated to verify the accuracy of the solution with more practical implementation such as the preconditioning technique and element-by-element storage. Further it would be beneficial to identify a procedure to maintain $C^1$ continuity at the enriched nodes as well. This will then extent the application of XDFEM to model curved discontinuities. To further analyze the advantageous of XDFEM, detailed comparison should be made between XDFEM and the popular smoothed XFEM in terms of precision, stability and efficiency.

\bibliography{XTFEM2}{}
\bibliographystyle{unsrt}

\end{document}